\documentclass[12pt,a4paper]{article}
\usepackage{epsfig,amssymb,color, graphics}
\usepackage{float}

\newtheorem{theo}{Theorem}

\newtheorem{prop}{Proposition}[section]
\newtheorem{defi}{Definition}[section]
\newtheorem{lemm}{Lemma}[section]
\newtheorem{coro}{Corollary}[section]
\newtheorem{stat}{Statement}[section]
\newtheorem{den}{Denotation}[section]
\newtheorem{rema}{Remark}[section]
\newenvironment{demo}{{\bf Proof: }}{\hfill $\diamond$\medskip}

\def\Om{\Omega}
\def\ga{\gamma}    
\def\Ga{\Gamma}
\def\de{\delta}

\def\vp{\varphi}
\def\la{\lambda}

\def\si{\sigma}
\def\Si{\Sigma}
\def\ep{\varepsilon} 
\def\nd{\noindent}

\begin{document}
\sloppy
\date{}
\title{Topological classification of Morse-Smale diffeomorphisms on 3-manifolds}
\author{Ch. Bonatti\thanks{Institut de Math\'ematiques de Bourgogne, UMR 5584 du CNRS, France, bonatti@u-bourgogne.fr}\and V.~Grines\thanks{Higher School of Economy, B. Pecherskaya, 25/12, N. Novgorod,
603155 Russia, vgrines@yandex.ru.}\and O.~Pochinka\thanks{Higher School of Economy, B. Pecherskaya, 25/12, N. Novgorod, 603155 Russia, olga-pochinka@yandex.ru.}}

\maketitle
\sloppy \tableofcontents

\begin{abstract} Topological classification of even the simplest Morse-Smale diffeomorphisms on 3-manifolds does not fit into the concept of singling out a skeleton consisting of stable and unstable manifolds of periodic orbits. The reason for this lies primarily in the possible ``wild'' behaviour of separatrices of saddle points. Another difference between Morse-Smale diffeomorphisms in dimension 3 from their surface analogues lies in the variety of heteroclinic intersections: a connected component of such an intersection may be not only a point as in the two-dimensional case, but also a curve, compact or non-compact. The problem of a topological classification of Morse-Smale cascades on 3-manifolds either without heteroclinic points (gradient-like cascades) or without heteroclinic curves was solved in a series of papers from 2000 to 2006 by Ch. Bonatti, V. Grines, F. Laudenbach, V. Medvedev, E. Pecou, O. Pochinka. The present paper is devoted to a complete topological classification of the set $MS(M^3)$ of  orientation preserving  Morse-Smale diffeomorphisms $f$ given on smooth closed orientable 3-manifolds $M^3$. A complete topological invariant for a diffeomorphism $f\in MS(M^3)$ is an equivalent class of its scheme $S_f$, which contains an information on a  periodic date and a topology of embedding of two-dimensional invariant manifolds of the saddle periodic points of $f$ into the ambient manifold. 
\end{abstract}



{\it Acknowledgement:} This work was supported by the Russian Science Foundation  (project 17-11-01041).

\section{Introduction and formulation of the results}\label{I}
\subsection{Informal statement of the results}

Morse-Smale systems are the simplest dynamical systems, and correspond to 
the evolution that one may imagine, before knowing the existence of chaotic evolutions: 
every orbit flow down to a (periodic) equilibrium point and comes (if it is not an equilibrium point himself) from another periodic point. More precisely,  a diffeomorphism is called  \emph{Morse-Smale} if it has finitely many periodic points, all of them are hyperbolic, the stable and unstable manifolds of any two periodic points are transverse, and every point in the manifold lies in both an unstable manifold and a stable manifold. 

Further this simple behaviour, the importance of Morse-Smale systems comes from the fact that they are \emph{structurally stable} (see \cite{Pa,PS}): the dynamics remains unchanged (that is, conjugated to itself) under small perturbation in the $C^1$-topology.

The time-$1$ map of the gradient flow of a generic Morse function is the typical example of Morse-Smale diffeomorphism (the genericity is needed to get the transversality condition). One could therefore  hope that  the manifold and the dynamic are characterized by  simple combinatorial informations on the periodic orbits, and the position of the invariant manifold. That is indeed the case  for Morse-Smale vector-field on compact  surfaces where a complete description and an a classification (up to topological equivalence) has been obtained by Peixoto \cite{Pe} and for 2-sphere it was a formalization of the Leontovich-Mayer's scheme \cite{LM2}. 

The problem is substantially more complicated in the case of diffeomorphisms. In particular, Morse-Smale diffeomorphisms (while dynamically as ``simple as possible'') are not necessarily embedded into a flow, even in a topological flow: the simplest obstruction is the existence of transverse intersection points between stable and unstable manifolds of complementary dimensions. This \emph{heteroclinic intersection points} 
lead to the main difficulty for classification of the Morse-Smale diffeomorphisms on surface. They admit an invariant like to Peixoto's graph in the case of finitely many heteroclinic orbits, what was proved by A.  Bezdenezhnykh and V. Grines \cite{BG1}, \cite{BG3}, \cite{G8}. Indeed a Morse-Smale diffeomorphisms on a compact surface may have infinitely many heteroclinic orbits, which cut the stable and unstable manifolds in infinitely many orbits of segments: the topological relative position of these segments in the surface are a topological conjugacy invariant of the diffeomorphisms. Nevertheless these infinitely many segments seems to follow a finite pattern: ending a long sequence of papers, Ch. Bonatti and R. Langevin  \cite{BoLa} provided a finite combinatorical complete invariant for Morse-Smale diffeomorphisms of surfaces, (as well as structurally stable diffeomorphisms with non-trivial basic sets).  

\begin{figure}[H]\centerline{\epsfig{file=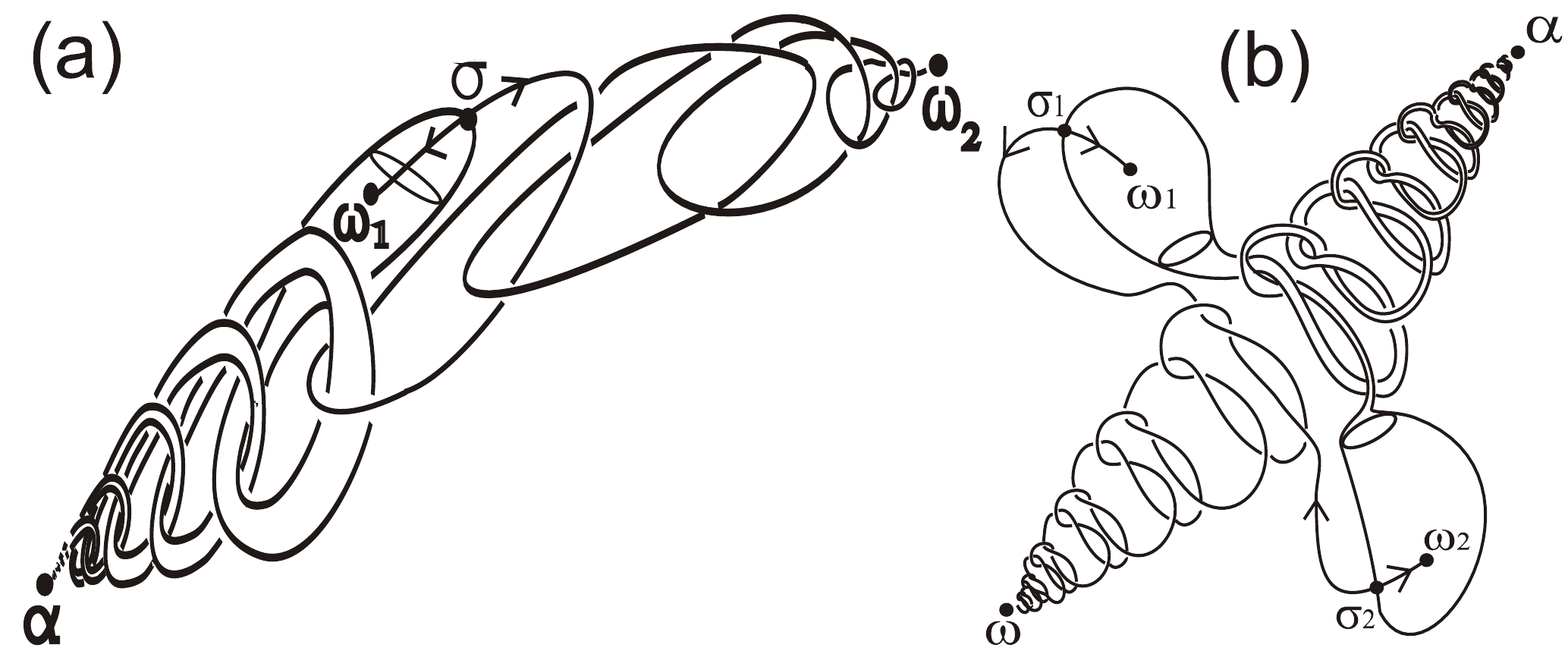, width=15. true cm, height=6. true cm}}
\caption{\footnotesize Morse-Smale diffeomorphisms with wild separatrices  on the $3$-sphere:
a) the separatrices of the saddle point $\sigma$ form an Artin-Fox arc with wild points in the source $\alpha$ and in the sink $\omega_2$; b) the separatrices of the saddle points $\sigma_1,\sigma_2$ form Debruner-Fox frame with wild points in the source $\alpha$ and in the sink $\omega$}  
\label{pr}
\end{figure}

Even Morse-Smale diffeomorphisms without heteroclinic intersections may not belong to a topological flow (see \cite{GrGuMePo12} for a characterization of this phenomenon,  see also \cite{BCVW} for an answer in the opposite direction). Thus the gluing of the dynamics from one periodic point to the next is not necessarily given by a flow, and can be very interesting from the topological viewpoint (see Figure~\ref{pr}). 

In view of the complexity of the classification of Morse-Smale diffeomorphisms on surface, and in view of the new topological behaviors allowed in dimensions $\geq 3$, a classification of Morse-Smale diffeomorphisms in dimension $3$ could appear hopeless. Nevertheless, below we describe a general approach to a classifiaction of  the dynamics of a Morse-Smale diffeomorphism in dimension $3$. 

The most simple among the Morse-Smale diffeomorphisms is the ``source-sink'' diffeomorphism whose non-wandering set consists of exactly two fixed points: a source and a sink. The ``source-sink'' diffeomorphisms have trivial dynamics: all the non-fixed points are wandering and under the action of the diffeomorphism they move from the source to the sink. Topological conjugacy of all these diffeomorphisms follows from the fact that the spaces of their wandering orbits are homeomorphic to $\mathbb{S}^{n-1}\times \mathbb{S}^1$ (in the orientation preserving case). 

The next simplest case in dimension $3$ is the orientation preserving Morse-Smale diffeomorphisms $f$ whose non-wandering sets consist of exactly four fixed points: two sinks $\omega_1,\omega_2$, one saddle point $\sigma$ with a $1$-dimensional unstable manifold, and a source $\alpha$ (see Figure \ref{pr} (a)).
The classification of such diffeomorphisms is the aim of \cite{BoGr} which shows that there is a one to one correspondence between the topological conjugacy classes of such Morse-Smale diffeomorphisms $f$ and 
the knots $\gamma$ in $\mathbb{S}^2\times \mathbb{S}^1$ in the homotopy class of the $\mathbb{S}^1$ factor, up to  homeomorphisms acting trivially on $\mathbb{Z}=H_1(\mathbb{S}^2\times \mathbb{S}^1)$. In that case, $\mathbb{S}^2\times \mathbb{S}^1$ is the space of orbits in 
the punctured basin $W^u_\alpha\setminus\{\alpha\}$, and a tubular neighborhood of the knot $\gamma$ is the projection on this orbit space of 
the punctured stable manifold $W^s_\sigma\setminus \{\sigma\}$. Note that the dynamics still looks line a ``source-sink'' dynamics, but the 
sink has been replaced by the invariant compact segment $\{\omega_1\}\cup W^u_\sigma\cup\{\omega_2\}$, which is an attracting set. 
This segment is wildly knotted if and only if the knot $\gamma$ is not trivial (where  \emph{trivial} means isotopic to the $\mathbb{S}^1$ factor). 

{Another  simple case consists in the Morse-Smale diffeomorphisms $f$ on a closed $3$-manifold $M^3$ which are the time one map of the gradient $X$ of a generic Morse function. 
In that case Smale's result \cite{Sm} implies the existence of a closed  (connected) surface $\Sigma_f$ transverse to $X$, disconnecting $M^3$ in two components: one repelling component 
and one attracting component; furthermore, the sources and the saddles with a $1$-dimensional stable manifolds  belong all to the repelling  component, 
and the sinks and the saddles with $1$-dimensional unstable manifolds belong all to the attracting component.
Here, proceeding our analogy with the source-sink dynamics, what plays the role of 
source (resp. sink) is the graph built by the union of the sources (sinks) and $1$-dimensional stable (unstable) manifolds of saddles.
Then, the surface $\Sigma_f$ intersects transversally every $2$-dimensional stable (resp. unstable) manifold of the saddles, 
and the intersection is exactly $1$ circle.  
Thus $\Sigma_f$ is equipped with a family $\mathcal{C}^s_f$ of disjoint stable circles, 
a family $\mathcal{C}^u_f$ of disjoint unstable circles, and these two families are transverse. 
Then $(\Sigma_f,\mathcal{C}^s_f,\mathcal{C}^u_f)$ is a complete invariant of topological conjugacy: 
the triple $(\Sigma_f,\mathcal{C}^s_f,\mathcal{C}^u_f)$ does not depend, up to homeomorphism, 
on the choice of the transverse surface $\Sigma_f$ and two such diffeomorphisms $f,f'$ are conjugated if 
and only if $(\Sigma_f,\mathcal{C}^s_f,\mathcal{C}^u_f)$ is homeomorphic to  
$(\Sigma_{f'},\mathcal{C}^s_{f'},\mathcal{C}^u_{f'})$ (see, for example, \cite{Pr}).}

When one studies a more complicated Morse-Smale diffeomorphism $f:M^3\to M^3$ on a closed 3-manifold $M^3$, we will see that the dynamics looks similar but ``the source'' and ``the sink'' then stand for the closed invariant sets of as simple topological structure as possible. One of them, denoted by $A_f$, 
is a compact attracting set and the other $R_f$ is a compact repelling set. Consider now the open set $V_f=M^3 \setminus(A_f\cup R_f)$. 
If the orbit space $\hat V_f=V_f/f$  can be canonically described, then it gives a good chance to get a complete topological invariant for these diffeomorphisms. 

More precisely, one considers the attractor $A_f$ (repeller $R_f$) as the closure of all one-dimensional unstable $W^u_p$ (stable $W^s_p$) manifolds of 
saddle points $p$. One can check that  $A_f$ ($R_f$) is a connected one-dimensional lamination in $M^3$ or a sink (source) in the exceptional case where there are no 
saddle points with the one-dimensional unstable (stable) manifold. These $1$-dimensional laminations do not disconnect the manifold $M^3$. Thus the set $V_f$ is a connected 
$3$-manifold, and one shows that the  orbit space  $\hat V_f$ is a closed connected manifold. Moreover, \cite{BoPa} proved that $\hat V_f$ is a {\it prime manifold}, 
that is a closed orientable 3-manifold which is either homeomorphic to $\mathbb S^2\times\mathbb S^1$ or irreducible (any smooth 2-sphere bounds a 3-ball there).

\vskip 5mm
\centerline{The closed manifold $\hat V_f$ is our \emph{first conjugacy invariant}\footnote{In this approach, the dimension $3$ play a specific role. In particular, 
in dimension $1$ and $2$, there is not such a canonical choice for the triple $A_f,R_f,V_f$, where all elements are connected; therefore  other approaches  have 
been used for the topological classification of Morse-Smale diffeomorphisms in dimension $1$ and $2$ (see section \ref{hi} below)}.}
\vskip 5mm
 
We denote by $$p_{_f}:V_f\to\hat V_f$$ the natural projection. Note that $p_{_f}$ is a cyclic cover whose deck transformation group is generated by $f$. Such a cyclic cover is 
associated to an epimorphism $\eta_{_f}:\pi_1(\hat{V}_f)\to \mathbb Z$ so that $[\hat c]\mapsto f^{\eta_{_f}([\hat c])}$ 
is the natural representation of $\pi_1(\hat{V}_f)$ in the deck transformation group. 

\vskip 5mm
\centerline{Then, $\eta_{_f}\in H^1(\hat{V}_f,\mathbb{Z})$ is our \emph{second conjugacy invariant.}}

\vskip 5mm
Our next invariants consist of  the projection in $\hat V_f$ of the $2$-dimensional invariant (stable or unstable) manifolds 
of the periodic saddle. The intersection with $V_f$ of the $2$-dimensional stable manifolds of the saddle points of $f$ is an invariant $2$-dimensional lamination $\Ga^s_{f}$, with finitely many leaves, and which is closed in $V_f$. Each leaf of this lamination is obtained by removing from a stable manifold its set of intersection points with the $1$-dimensional unstable manifold; this intersection is at most countable. As $\Ga^s_{f}$ is invariant under $f$, it passes to the quotient in a compact 
$2$-dimensional lamination $\hat\Ga^s_{f}$ on $\hat V_f$. Note that each $2$-dimensional stable manifold is a plane on which $f$ acts as a contraction, so that the quotient by $f$ of the punctured  stable manifold is either a torus or a Klein bottle.  Thus the leaves of $\hat\Ga^s_{f}$ are either tori or Klein bottles which are punctured along at most countable set.

One defines in the same way the unstable lamination $\hat \Ga^u_{f}$ as the quotient by $f$ of the intersection with $V_f$ of the $2$ dimensional unstable manifolds. The laminations $\hat\Ga^s_{f}$ and $\hat \Ga^u_{f}$ are transverse. 

\vskip 5mm

\centerline{The laminations $\hat \Ga^s_{_f}$ and $\hat \Ga^u_{_f}$ are our \emph{last conjugacy invariants}}

\vskip 5mm

We denote $S_f=(\hat V_f,\eta_{_f},\hat\Gamma^s_f,\hat\Gamma^u_f)$ and we call $S_f$ the \emph{scheme associated to $f$}.
More generally, in our setting, an \emph{abstract scheme} is a collection
$S= (\hat V,\eta,\hat\Ga^s,\hat\Ga^u)$ where $\hat V$ is a closed $3$-manifold, $\eta$ is an epimorphism from $\pi_1(\hat V)$ to $\mathbb{Z}$ and 
$\hat \Ga^s$, $\hat \Ga^u$ are two $2$-dimensional 
transverse compact laminations in $\hat V$. Two schemes $S=(\hat V,\eta,\hat\Ga^s,\hat\Ga^u)$ and $S'=(\hat V',\eta',\hat\Ga'^s, \hat\Ga'^u)$ are \emph{equivalent} 
if there is a homeomorphism $h\colon \hat V\to\hat V'$ so that $h_*(\eta)=\eta', h(\hat\Ga^s)=\hat \Ga'^s$ and $h(\hat\Ga^u)=\hat\Ga'^u$. 

{\bf Realization problem}\emph{ What abstract scheme $S=(\hat V,\eta,\hat\Ga^s,\hat\Ga^u)$ is equivalent to the scheme $S_f$ associated with a 
Morse-Smale diffeomorphism $f:M^3\to M^3$}? 

The realization problem has been solved in the general setting in  \cite{BoGrPo2017} generalizing the solution proposed in \cite{BoGrMePe3} 
for the gradient-like diffeomorphisms. 

The aim of this paper is to prove that $S_f$ is a complete conjugacy invariant for (orientation preserving) Morse-Smale diffeomorphisms: 

{\bf Theorem (Classification of Morse-Smale diffeomorphisms on 3-manifolds)}\emph{Let $M^3$ and $M'^3$  be  closed orientable $3$-manifolds 
and $f\colon M^3\to M^3$ and $f'\colon M'^3\to M'^3$ be two orientation preserving Morse-Smale diffeomorphisms. Then  $f,f'$ are topologically 
conjugate if and only if  the scheme $S_f$ is equivalent to the scheme  $S_{f'}$.}

The fact that, if $f$ is conjugated to $f'$ then the scheme $S_f$ and $S_{f'}$ are conjugated is almost ``by construction'': a conjugacy homeomorphism $h\colon M^3\to M'^3$ maps $V_f$ to $V_{f'}$ and the fact that it is a conjugacy implies that it passes to the quotient in an homeomorphism $\hat h\colon \hat V_f\to \hat V_{f'}$, which is an equivalence between the two schemes. The aim of the paper is to prove the reverse implication. 

\subsection{Intuitive idea for the proof}\label{intu}

Let $M^3$ and $M'^3$ be closed orientable $3$-manifolds and $f\colon M^3\to M^3$ and $f'\colon M'^3\to M'^3$ 
be two orientation preserving Morse-Smale diffeomorphisms. Assume that there is a homeomorphism $\hat\varphi\colon \hat V_f \to \hat V_{f'}$ which realises an equivalence between the scheme $S_{f}=(\hat V_f,\eta_f,\hat\Ga^s_{f},\hat\Ga^u_f)$  and $S_{f'}= (\hat V_{f'},\eta_{f'},\hat\Ga^s_{f'},\hat\Ga^u_{f'})$. 

As $\hat\varphi_*(\eta_f)=\eta_{f'}$ the homeomorphism $\hat\varphi$ admits a lift $\varphi\colon V_f\to V_{f'}$ and $h$ conjugates the generator of the deck transformation groups of the covers
$p_{_f}\colon V_f\to \hat V_f$ and $p_{_{f'}}\colon V_{f'}\to \hat V_{f'}$.  In other words, $\varphi$ conjugates the restrictions of $f$ and $f'$ to $V_f$ and $V_{f'}$, respectively. 
It is not hard to check that $\varphi$ can be extended on the periodic sources, sinks and saddles. 

If the homeomorphism $\varphi$ admits a continuous extension to the $1$-dimensional stable and unstable manifold of the saddle points, then it induces a conjugacy between $f$ and $f'$, ending the proof. However, in general, $\varphi$ does not admit such a continuous extension on the $1$-dimensional invariant manifolds:
we will need to modify $\varphi$ in order to get a new homeomorphism, extending on the $1$-dimensional invariant manifolds, and that is the aim of this paper. 

Let us first present a very tempting conceptual and global approach.  Each saddle point $\sigma$ is equipped with an invariant 
linearising neighborhood $N_\sigma$ with invariant stable $F^s_\sigma$ and unstable $F^s_\sigma$ foliations. Furthermore, these system of 
local foliations is \emph{compatible}: these invariant neigborhoods may intersect, and on the intersection, the stable leaves of dimension $1$ are contained 
in the $2$-dimensional stable leaves, the stable foliations of the same dimensions coincide, and so on. The existence of such compatible system of foliations comes 
back to the proof of the structural stability of the Morse-Smale diffeomorphisms. 

These locally invariant compatible foliations induce  on the quotient $\hat V_f$ and $\hat V_{f'}$ compatible foliations defined in the neighborhood of the stable and unstable laminations $\hat \Ga^s_f\cup \hat \Ga^u_f$ and $\hat \Ga^s_{f'}\cup \hat \Ga^u_{f'}$, respectively. If one can modify $\hat \varphi\colon \hat V_f\to\hat V_{f'}$ in order that it preserves these foliations, we are done:  each point in a one-dimensional  
stable or unstable manifold is the intersection point of this manifold with a $2$-dimensional leaf. 

This approach works well when the diffeomorphisms $f$ and $f'$ have no heteroclinic intersections. In that case the laminations $\hat \Ga_f^s$ and $\hat \Ga_f^u$ are disjoint and consist of finitely many tori or Klein bottles. Then we could modify $\hat\varphi$ in a neighborhood of the lamination $\hat \Ga_f^s$ and $\hat \Ga_f^u$ in order that it preserves the system of compatible foliations, and we could glue this 
local modification with the old homeomorphism far from the laminations.  This comes from the following general fact:
\begin{prop}\label{TT}For given a torus $T$ embedded to a $3$-manifold $M^3$ and a local orientation preserving  homeomorphism $\psi$ defined in a 
neighborhood $U(T)$ of $T$ and inducing the identity by restriction on $T$, there is homeomorphism $\Psi$ which is the identity map in a small neighborhood of 
$T$ and which coincides with $\psi$ in $\partial U(T)$. 
\end{prop}

We have not been able to get such a statement for the general and transversely intersecting laminations $\hat \Ga_f^s$ and $\hat \Ga_f^u$. 
For this reason we have given up on this global conceptual approach, and go back to a progressive, step by step, approach. What we do is to consider 
the saddle points one by one. First we consider the saddle $\sigma_0$ whose $1$-dimensional unstable manifold is not accumulated by other. In other words, 
we consider saddle whose $2$-dimensional stable manifold do no intersect the $1$-dimensional stable manifolds. For such a saddle, we modify $\varphi$ in a neighborhood of its 
unstable manifolds. Then we will consider the saddle $\sigma_1$ whose $1$-dimensional unstable manifold is only accumulated by the ones on which we have already done the 
modification.  Then we will perform a modification of $\varphi$ preserving the modifications which have been already done, and preparing the next modifications. 
In this way we will consider one by one the saddles $\sigma_0,\dots,\sigma_n$ with $1$-dimensional unstable manifold. When we performed all the modifications for getting 
an extension along every $1$-dimensional unstable manifold, we will consider the saddle points with the other index, with $1$-dimensional stable manifold.  We will check 
that we can perform the same kind of modifications without breaking the extensions which have been already done. 

Le us try to explain, very roughly, how we modify the conjugacy homeorphism $\varphi$ in the neighborhood of a periodic orbit. First notice that the homeomorphism $\varphi$, 
in restriction to a given
$2$-dimensional unstable manifold extends by continuity in a unique way on the heteroclinic points (which are at most countably many). Therefore $\varphi$ can be considered 
as defined on the whole $2$-dimensional unstable (or stable manifold). We also define a conjugacy homeomorphism $\psi^s$ on the union of the $1$-dimensional stable manifold, so that 
$\psi^s$ preserves the holonomies of the local $2$-dimensional (compatible) unstable foliations.  In any small linearizing neighborhood of a saddle point $\sigma$, one gets a local conjugacy 
homeomorphism $\phi_\sigma$ whose expression is the product of the restriction of $\varphi$ to $W^u_\sigma$ by $\psi^s$ on $W^s_\sigma$.  Thus $\xi_\sigma= \phi^{-1}_\sigma\varphi$ is a
homeomorphism defined in a neighborhood of $\sigma$ and commuting with the diffeomorphism $f$ and inducing the identity map on $W^u_\sigma$.

The goal is now to build a homeomorphism commuting with $f$, coinciding with the identity map in a neighborhood of $\sigma$ and coinciding with $\xi_\sigma$ out of 
a slightly larger neighborhood.  The difficulty is that this 
homeomorphism needs to preserve the intersection of this neighborhood with the invariant manifold crossing the neighborhood. We want also to preserve the extension we already have done in a neighborhood of the other saddles.  These neighborhoods are crossing the neighborhood of $\sigma$ as vertical tubes, that we want to preserve. As $W^u_\sigma$ locally disconnects the manifold, we work only on $1$ side of $W^u_\sigma$.  
One choose a  fundamental domain of the local dynamics which is the product of an annulus (fundamental domain in $W^u_\sigma$) by a vertical segment.  

The main topological difficulty we need to solve is to extend a homeomorphism defined in a neighborhood of this fundamental domain preserving  at the same time the $2$-dimensional horizontal unstable leaves, 
some $2$-dimensional stable leaves also crossing, and the homeomorphism in some crossing tube. This will be done by using variations of  Proposition \ref{TT}. 
For being able to solve this puzzle, we try first to simplify as most as possible the 
geometry of these intersections. For these modifications, we will use results from the topology of $3$-manifolds, allowing to extend in the interior of a region homeomorphism defined on the boundary. 
This topological tools are built from  classical results in Section~\ref{IIII}. 

\subsection{Illustrations to the scheme $S_f$}

Indeed one may see the factor space $\hat{V}_f=V_f/f$ as obtained by taking a fundamental domain in the basin of the  attractor $A_f$ and identifying its boundaries by $f$. 
It will be helpful  to define the epimorphism 
$$\eta_{_f}:\pi_1(\hat{V}_f)\to \mathbb Z$$ 
as follows. Take the homotopy class $[\hat c]\in\pi_1(\hat{V}_f)$ 
of a closed curve $\hat c\colon \mathbb R/\mathbb Z\to \hat{V}_f$. Then $\hat c\colon [0,1]\to \hat{V}_f$ 
lifts to a curve $c\colon [0,1] \to V_f$ connecting a point $x$ with a point $f^{n}(x)$ for some  $n\in \mathbb Z$, where $n$ is independent of the lift.
So define $\eta_{_f}[\hat c]=n$. Denote by $$p_{_f}:V_f\to\hat V_f$$ the natural projection. Due to the hyperbolicity of the periodic points, 
the restriction of $f^{m_p}$ to the two-dimensional stable  manifold of saddle point $p$ of the period $m_p$ is conjugated with the linear contraction. 
Then the space orbit of the two-dimensional separatrix is homeomorphic to the 2-torus or the Klein bottle depending on the map $f^{m_p}|_{W^s_p}$
preserves orientation or changes orientation, that we will indicate $\nu_p=+1$ or $\nu_p=-1$, accordingly. Thus, if $(W^s_p\setminus p)\subset V_f$ then the projection
$\hat W^s_p$ of $W^s_p$ to $\hat V_f$ is a smoothly embedded 2-torus or a Klein bottle. Moreover, $$\eta_{_f}(j_*(\pi_1(\hat W^s_p)))=m_p\mathbb Z$$ for the inclusion 
$j:\hat W^s_p\to \hat V_f$. If $W^s_p\setminus p$ has a non empty intersection with $A_f$ then $\hat W^s_p$ is a punctured torus or Klein bottle where the number of 
punctured points equals the number of heteroclinic orbits in $W^s_p$.  It will be similarly for two-dimensional unstable saddle separatrices.  Let us demonstrate it by examples.

\begin{figure}[h]
\centering\epsfig{file=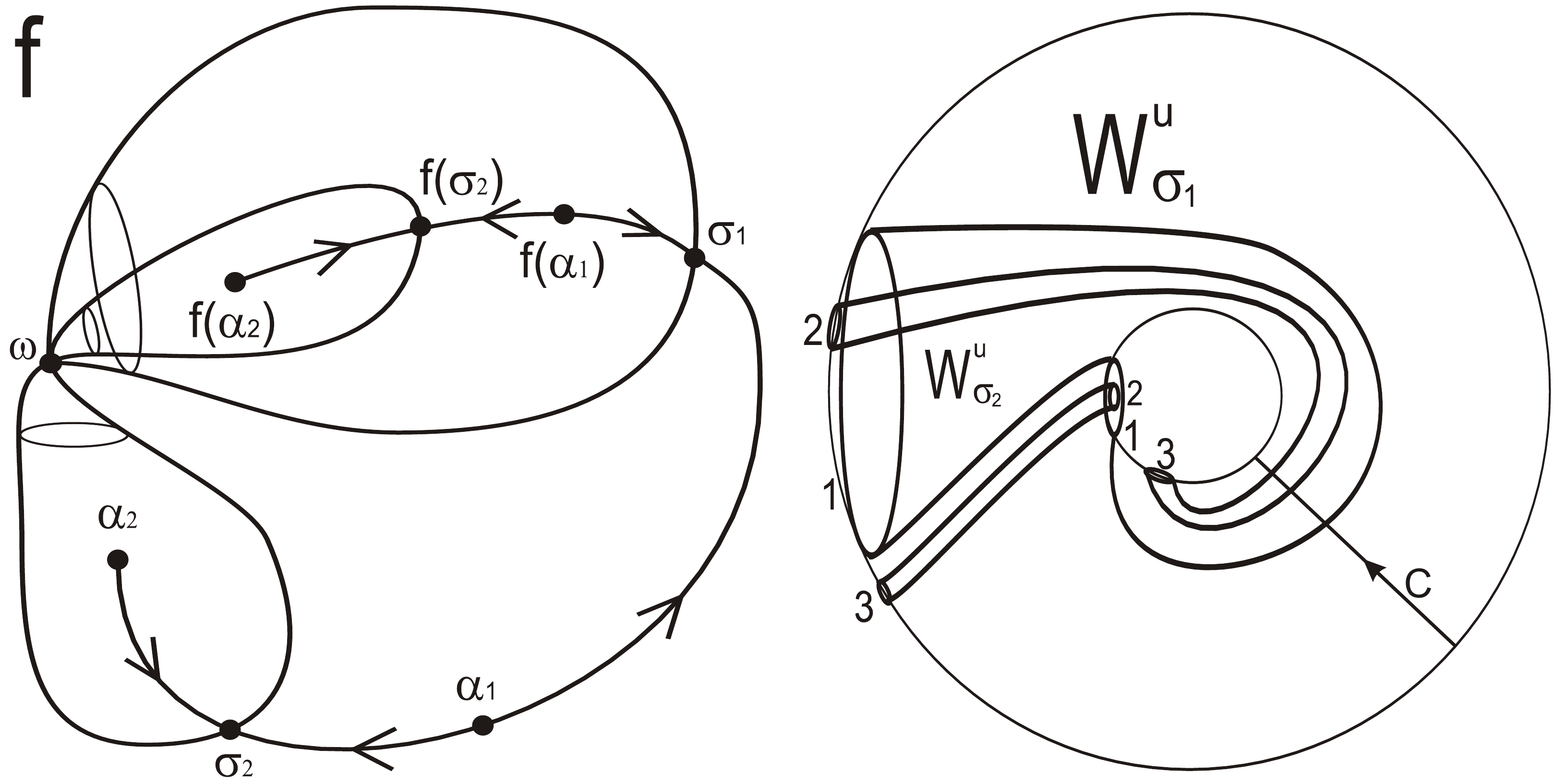, width=12. true cm,
height=6. true cm}\caption{\small The two-dimensional separatrices in a fundamental domain of $V_f$} \label{8888+}
\end{figure}

The left-side of Figure \ref{8888+} shows a Morse-Smale diffeomorphism $f:\mathbb S^3\to\mathbb S^3$ the non-wandering set of which consists of five periodic points $p$ with the following periodic data $p(m_p,dim\, W^u_p,\nu_p)$:  ${\omega}(1,0,+1)$, ${\sigma_1}(1,2,-1)$, ${\sigma_2}(2,2,+1)$, ${\alpha_1}(2,3,+1)$, ${\alpha_2}(2,3,+1)$. Then $A_f=\omega$ and the right-side of Figure \ref{8888+} shows the fundamental domains of the action of the diffeomorphism $f$ on ${W^s_{\omega}\setminus\omega}$. The domain is a 3-annulus from which the orbits spaces $\hat V_f$ is obtained by gluing the boundary spheres of the annulus by  $f$ and, hence, $\hat V_f$ is diffeomorphic $\mathbb S^2\times\mathbb S^1$ and $\eta_{_f}([p_{_f}(c)])=1$. The orbits spaces $\hat W^u_{{\sigma_i}},~i=1,2$ are obtained from the cylinders by a gluing the points with the same numbers so that $\hat W^u_{{\sigma_1}}$ is a Klein bottle and $\hat W^u_{{\sigma_2}}$ is a double-round torus. 
\begin{figure}[h]
\centering\epsfig{file=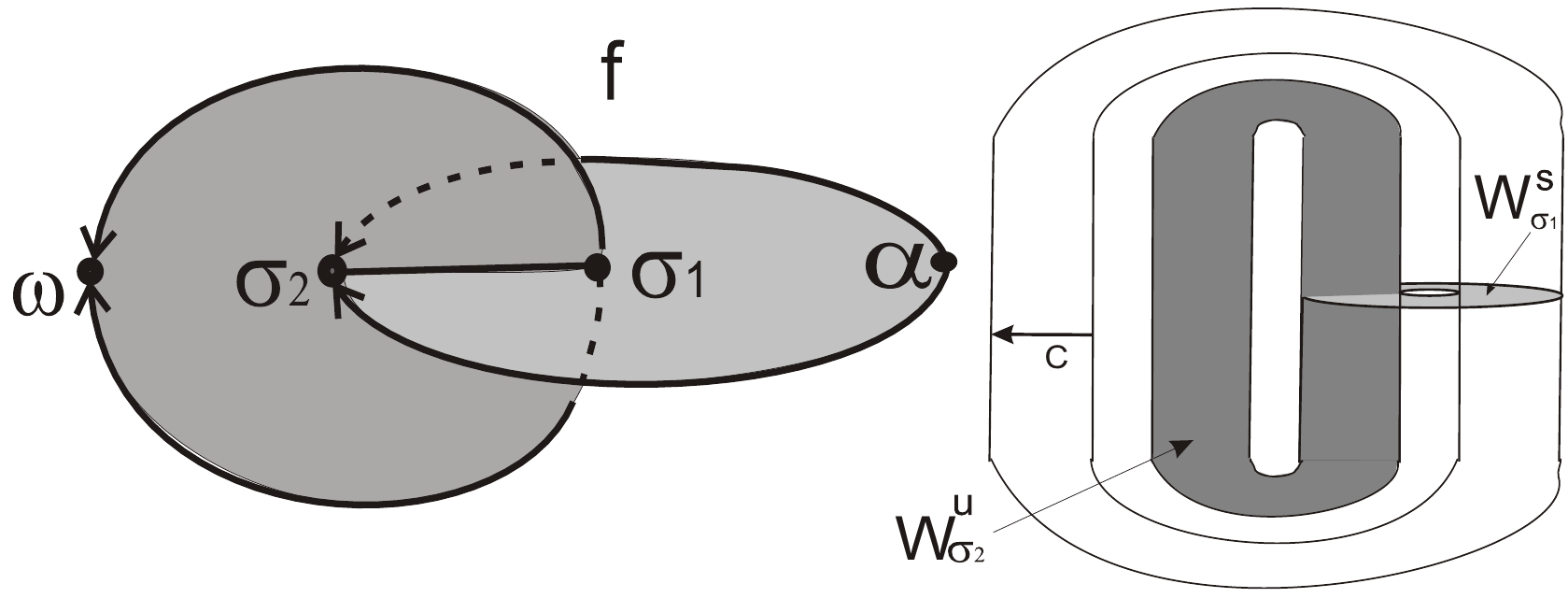, width=13. true cm,
height=6. true cm}\caption{\small The two-dimensional separatrices with a heterovlinic curve in a fundamental domain of $V_f$} \label{baterfgl+}
\end{figure}

Figure \ref{baterfgl+} shows the phase portrait of a Morse-Smale diffeomorphism $f$ of the 3-sphere which is the time-1
map of a flow on  $\mathbb S^3$. The non-wandering set of this diffeomorphism consists of the source $\alpha$, the sink $\omega$ and the two saddle points $\sigma_1$, $\sigma_2$ 
with one-dimensional and two-dimensional unstable manifolds, respectively. We also assume that the stable manifold of the point $\sigma_1$ intersects the unstable manifold of the point $\sigma_2$ by the unique non-closed non-compact heteroclinic curve. In this case the attractor $A_f=cl(W^u_\sigma)$ is a closed curve which is homeomorphic to the circle and it is tamely embedded into the ambient manifold $\mathbb S^3$. Moreover the attractor $A_f$ has a trapping neighborhood $M_f$ homeomorphic to the solid torus and such that $M_f\setminus int\,  f(M_f)$ is homeomorphic to $\mathbb T^2\times[0,1]$. Therefore the factor space $\hat V_f$ for this diffeomorphism is homeomorphic to the 3-torus $\mathbb T^3$ and $\eta_{_f}([p_{_f}(c)])=1$. The orbits spaces $\hat W^u_{{\sigma_i}},~i=1,2$ are obtained from the annuli by a gluing the boundary points so that $\hat W^u_{{\sigma_i}}$ are tori.
\begin{figure}[h]
\centerline{\includegraphics[width=12cm,height=8cm]{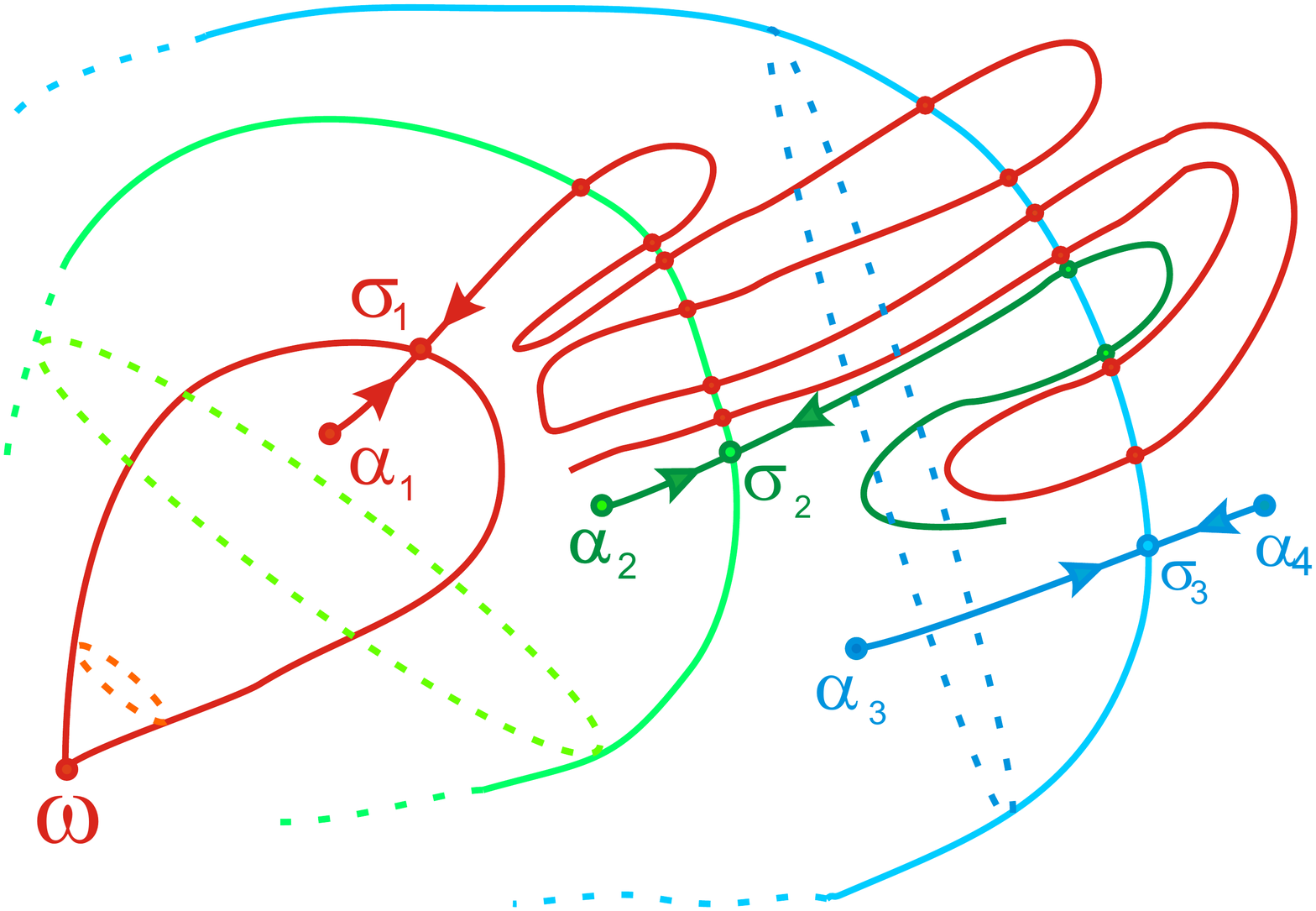}}\caption{\small Phase portrait of a Morse-Smale diffeomorphism with heteroclinic points on a 3-manifold} \label{faz}
\end{figure}
\begin{figure}[h]
\centerline{\includegraphics[width=16cm,height=6cm]{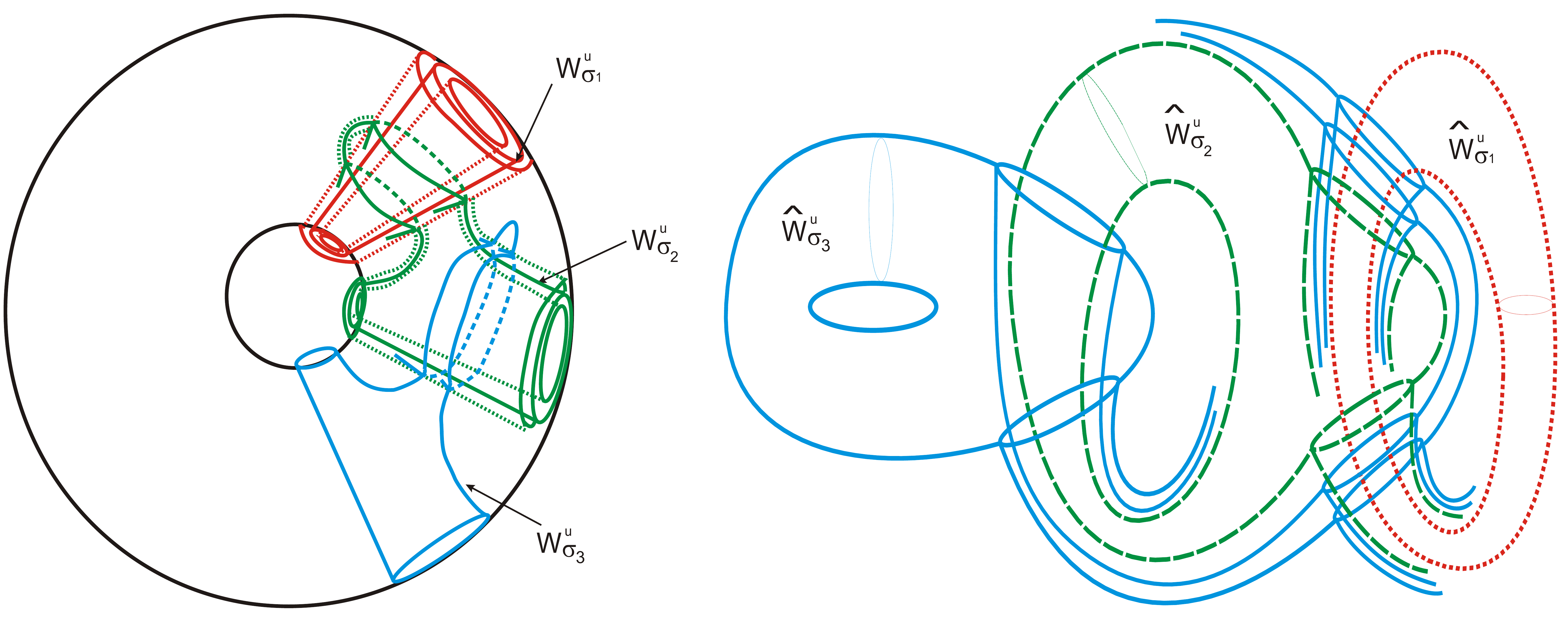}}\caption{\small An $u$-lamination of the  diffeomorphism $f$ from Figure \ref{faz}} \label{sch}
\end{figure}

To represent an embedding of a punctured surface let us consider a dynamical system whose phase portrait is on Figure \ref{faz}. For this diffeomorphism $A_f$ is a sink $\omega$ and fixed saddle points $\sigma_1,\sigma_2,\sigma_3$ satisfy the conditions: $W^u_{\sigma_1}\cap A_f=\emptyset$, $W^u_{\sigma_2}\cap A_f$ consists of a finite number heteroclinic orbits, $W^u_{\sigma_3}\cap A_f$ consists of a countable many heteroclinic orbits. 
The quotient $\hat V_f$ is diffeomorphic $\mathbb S^2\times\mathbb S^1$ and $\hat W^u_{\sigma_1}\cup\hat W^u_{\sigma_2}\cup \hat W^u_{\sigma_3}$ form an $u$-lamination 
(see Figure \ref{sch}).
   
For arbitrary Morse-Smale diffeomorphism $f:M^3\to M^3$ we get two transversally intersected $s$-lamination $\hat\Gamma^s_f$ and $u$-lamination 
$\hat\Gamma^u_f$ on $\hat V_f$, each leaf of which is either torus or Klein bottle with empty, finite or countable set of punctured points. 
     
\subsection{A historical background of the topological classification of Morse-Smale systems}\label{hi}

In 1937 A. Andronov and L. Pontryagin \cite{AP} introduced the concept of a {\it rough system} of   differential equations defined in a bounded part of the plane: a system that retains its 
qualitative properties under small changes in the right-hand side. They proved that the flow 
generated by such a system is characterized by the following properties:
\begin{itemize}
\item the set of fixed points and periodic orbits is finite and all its elements are 
hyperbolic;
\item there are no separatrices  from a saddle to a saddle;
\item all $\omega$- and $\alpha$-limit sets are 
contained in the union of fixed points and
periodic orbits (limit cycles).
\end{itemize}
The above properties are also known to characterize rough flows on a two-dimensional sphere. A similar result for flows with closed section and  without equilibrium states  on a two-dimension tori follows from the paper of  A. Mayer \cite{Mai} in 1939. The principal difficulty in passing from a two-dimensional sphere
to orientable surfaces of positive genus is the possibility that there may exist new
types of motion --- non-closed recurrent trajectories.  A. Andronov and L. Pontryagin  have shown also in \cite{AP} that the set of  the rough flows is  dense in the space of $C^{1}$-flows\footnote{This statement  was not explicitly formulated in \cite{AP} and  was  mentioned first time in papers of E. Leontovich \cite{Le} and M. Peixoto \cite{Pe0}. G. Baggis \cite{Ba} in 1955 made explicit  some details of the  proofs in \cite{AP}, which were not published.}. A similar criterion and density were proved by M. Peixoto \cite{Pe1}, \cite{Pe2} in 1962 for {\it structural stable flows}\footnote{This concept  is  a generalization of the notion of  roughness. The original definition of a rough flow involved the additional requirement that the conjugating 
homeomorphism  be $C^0$-close to the identity map. The concepts of ``roughness'' and ``structural stability'' are now known to be equivalent, though this fact is highly
non-trivial.} on orientable surfaces 
of genus greater than zero. An immediate generalization of properties of rough flows on orientable surfaces leads to Morse-Smale systems (continuous and discrete).

\begin{defi} \label{oM-S} A smooth dynamical system  given on an $n$-manifold  $(n\geq 1)$ $M^n$ is called Morse-Smale system if:

1) its non-wandering set consists of a finite number fixed points and periodic orbits each of them is hyperbolic;

2) the stable and unstable manifolds  $W^s_p$, $W^u_q$ intersect transversally
\footnote{Two smooth submanifolds  $X_1$, $X_2$ of an $n$-manifold $X$ 
{\it are transversaly intersected} if either $X_1\cap X_2=\emptyset$ or $T_{x}X_1+T_{x}X_2=T_xX$ 
for each point $x\in (X_1\cap X_2)$ (here $T_{x}A$ is denotation for tangent space to the manifold 
$A$ at the point $x$).} for any non-wandering points $p$, $q$.
\end{defi}

Morse-Smale systems are named after Smale's 1960 paper \cite{S2} in 1960, where he first 
studied  flows with the properties above and proved that they satisfy inequalities
similar to the Morse inequalities. That Morse-Smale systems are structurally stable
was later shown by S. Smale and J. Palis \cite{Pa}, \cite{PS}. However, already in 1961 S. Smale \cite{S1} that such systems do not exhaust the class of all rough systems, constructing
for this purpose a structurally stable diffeomorphism on the two-dimensional sphere
$\mathbb S^2$ with infinitely many periodic points. This diffeomorphism is now known as the
``Smale horseshoe''. Nevertheless, these systems have great value both in applications
(because they adequately describe any regular stable processes) and in studying the
topology of the phase space (because of the deep interrelation between the dynamics
of these systems and the ambient manifold). 

The crucial question in the study of dynamical systems is the determination of the set
of complete topological invariants: properties of a system that uniquely determine the decomposition of the phase space into trajectories up to topological equivalence
(conjugacy). We recall that two flows  $f^{t}$, $f^{\prime t}$ (two diffeomorphisms $f$, $f^{\prime}$) on an $n$-manifold $M^n$  are said to be {\it topologically equivalent} ({\it topologically conjugate}) if there exist a homeomorphism $h:M^n\to M^n$ that carries  trajectories of $f^{t}$ to trajectories of  $f^{\prime t}$ (which satisfies to the condition $f'h=hf$). Today, this problem has a rich history.

The equivalence class of Morse-Smale flows on a  circle is uniquely determined by the number of 
its fixed points. For  cascades on a circle, A. Mayer \cite{Mai} found in 1939
a complete topological invariant consisting of a triple of numbers: the number of
periodic orbits, their periods, and the so-called ordinal number. In 1955 E. Leontovich and A. Mayer \cite{LM2} introduced a compete topological invariant --- the scheme of a flow --- for flows with finitely many singular trajectories on a two-dimensional sphere. This scheme contained a description of singular trajectories (equilibrium states, periodic orbits, separatrices of saddle equilibrium states) and their relative
positions. In 1971, M. Peixoto \cite{Pe} formalized the notion of a Leontovich-Mayer's 
scheme and proved that for a Morse-Smale flow on an arbitrary surface a complete topological invariant is given by the isomorphism class of a directed graph associated with it whose vertices are in a one-to-one correspondence with the equilibrium states and closed trajectories and whose edges correspond to the connected components of the invariant manifolds of the equilibrium states and closed trajectories, where
the isomorphisms preserve specially chosen subgraphs\footnote{In \cite{OshSh}  
A. Oshemkov and V. Sharko pointed out a certain inaccuracy concerning the Peixoto invariant due to the fact that an isomorphism of graphs does not distinguish between types of decompositions into trajectories for a domain bounded by two periodic orbits.}.

Morse-Smale flows (cascades)  given on manifolds of dimension $n\geq 3~(n\geq 2)$ feature
a new type of motion compared with lower-dimensional systems, because of possible  intersections of the invariant manifolds of distinct saddle points --- {\it heteroclinic intersection}.  V. Afraimovich and L. Shilnikov \cite{ASh} proved that the restriction of
Morse-Smale flows to the closure of the set of heteroclinic trajectories is conjugate
to a suspension over a topological Markov chain. Nevertheless, an invariant similar to the Peixoto graph proved to be sufficient for describing a complete topological invariant for a broad subclass of such systems, and in particular for Morse-Smale diffeomorphisms on surfaces with finitely many heteroclinic orbits (A.  Bezdenezhnykh, V. Grines \cite{BG1}, \cite{BG3} in 1985, V. Grines \cite{G8} in 1993)\footnote{One should also point out that R. Langevin \cite{La} proposed a different approach to finding topological invariants for such diffeomorphisms. Though no classification results were given in \cite{La}, the
ideas there nevertheless turned out to be very fruitful and have been put to use in the classification of diffeomorphisms, as is demonstrated, in particular, in the present survey. A classification of
Morse-Smale diffeomorphisms on surfaces with infinitely many heteroclinic orbits which required
the machinery of topological Markov chains follows from the paper \cite{BoLa} by Ch. Bonatti and R. Langevin, where necessary and sufficient conditions for topological conjugacy of Smale diffeomorphisms ($C^1$-structurally stable diffeomorphisms) on surfaces were established.}; for flows with a finite number of singular trajectories  on 3-manifolds (Ya. Umansky \cite{Um90} in 1990);  for flows on the sphere $\mathbb S^{n},~n\geq 3$ without closed orbits (S. Pilyugin \cite{Pil78} in 1978); for diffeomorphisms on $M^n,~n>3$ 
with saddle points of the  Morse index  one (V. Grines, E. Gurevich, V. Medvedev \cite{GrGu07}, 
\cite{GrGuMe08} in 2007-2008).

Topological classification of even the simplest Morse-Smale diffeomorphisms on
3-manifolds does not fit into the concept of singling out a skeleton consisting of
stable and unstable manifolds of periodic orbits. The reason for this lies primarily
in the possible ``wild'' behaviour of separatrices of saddle points. More specifically,
even though the closure of a separatrix may differ from a separatrix by only one
point, it may fail to be even a topological submanifold. D. Pixton \cite{Pi} in 1977 was the
first to construct a diffeomorphism with wild separatrices --- for this he employed
the Artin-Fox curve \cite{ArFo} to realize the invariant manifolds of a saddle fixed point (see Figure \ref{ex4}). Ch. Bonatti and V. Grines  \cite{BoGr} in 2000 investigated the class of diffeomorphisms on a three-dimensional sphere (diffeomorphisms in the Pixton class)
that have non-wandering set consisting of four fixed points: two sinks, a source,
and a saddle. They showed that the Pixton class contains a countable set of pairwise
topologically non-conjugate diffeomorphisms. Furthermore, the topological conjugacy class of a diffeomorphism from the Pixton class is uniquely determined by the embedding type of a one-dimensional separatrix in the basin of a sink, which is described by a new topological invariant: a smooth embedding of the circle $\mathbb S^1$ (the orbit space
of a one-dimensional separatrix) in the manifold $\mathbb S^2\times\mathbb S^1$ (the space of wandering orbits in the basin of a sink).

Another difference between Morse-Smale diffeomorphisms in dimension 3 from
their surface analogues lies in the variety of heteroclinic intersections: a connected
component of such an intersection may be not only a point as in the two-dimensional
case, but also a curve, compact or non-compact. The problem of a topological classification of Morse-Smale cascades on 3-manifolds either without heteroclinic points (gradient-like cascades) or without heteroclinic curves was solved in a series of papers \cite{BoGr}, \cite{BoGrLaPo},  
\cite{BoGrMePe3}, \cite{BoGrPo2005}, \cite{BoGrPo2} from 2000  by Ch. Bonatti, V. Grines, F. Laudenbach, V. Medvedev, E. Pecou, O. Pochinka. In the present paper the topological classification of the set $MS(M^3)$ of  
preserving orientation Morse-Smale diffeomorphisms $f$ given on a smooth closed orientable 3-manifolds $M^3$ is obtained. 

Thus this paper is the last stone in a very long construction. But the result here is not a generalization of the previous results and there are a number of reasons for this. In the absence of heteroclinic points, each connected component of $\hat{\Gamma}^s_f,\hat{\Gamma}^u_f$ in the scheme $S_f$ is a compact surface, a torus or a Klein bottle (not a lamination, as in the general case) and a modification of the homeomorphism of schemes in the neighborhoods of these surfaces, preserving the two-dimensional foliation (see the idea of the proof in Section \ref{intu}), is reduced to the proof of topological facts of the type of Proposition \ref{TT}. In the absence of heteroclinic curves, the connected components in the scheme can already be by a lamination, as in the general case. But, the fact that laminations of different stability do not intersect greatly simplifies the dynamics of such diffeomorphisms; such diffeomorphisms exist only on manifolds that are connected sums of finitely many copies of $\mathbb S^2\times\mathbb S^1$. Accordingly, the construction of the conjugate homeomorphism is simpler than in the general case. As in the gradient-like case, we modify the homeomorphism of the schemes in the neighborhood of the lamination, preserving the two-dimensional foliation. In this case, the modification is done step by step from one component of the linear connection of the lamination to the other. We use author's topological facts that allow us to modify the homeomorphism on the current component without changing it on the previous one. However, all these achievements are still not enough for the general case, in which a transversal intersection of laminations of different stability is allowed. Modification in the neighborhood of each lamination will have to be done in such a way that changes in the neighborhood of stable lamination do not change the already adjusted homeomorphism in the neighborhood of unstable lamination. For this, the authors of the article invented thin exclusive topological constructions.

\subsection{The exact formulation of the results}
Let $f\in MS(M^3)$. For $q=0,1,2,3$ denote by $\Omega_q$ the set of all periodic 
points of $f$ with the $q$-dimensional  unstable manifold. Let us represent the dynamics 
of $f$ in a form ``source-sink'' by a following way. Set  $A_f=W^u_{\Omega_0\cup\Omega_1}$, 
$R_f=W^s_{\Omega_2\cup\Omega_3}$ and $V_f=M^3\setminus(A_f\cup R_f)$. Then the set $A_f~(R_f)$  is a  connected attractor (repeller)\footnote{A compact set
$A\subset M^n$ is {\it an attractor of a diffeomorphism $f:M^n\to M^n $} if
there is a neighborhood $U$ of the set $A$ such that $f (U) \subset int~U$ and
$A=\bigcap\limits_{n\in\mathbb{N}}f^n(U)$. A set $R\subset M^n $ is called {\it a
repeller} of  $f$ if it is an attractor of  $f ^ {-1} $.} of $f$ with the topological dimension less or equal than 1,  the set $V_f$ is a connected 3-manifold and $V_f=W^s_{A_f\cap \Omega_f}\setminus A_f=W^u_{R_f\cap \Omega_f}\setminus R_f$. Moreover, a quotient $\hat V_f=V_f/f$ is a closed connected orientable 3-manifold, on which the natural projection $p_{_f}:V_f\to\hat V_f$ induces an epimorphism 
$\eta_{_f}:\pi_1(\hat V_f)\to\mathbb Z$, assigning a homotopy class 
$[c]\in\pi_1(\hat V_f)$ of a closed curve $c\subset\hat V_f$ an integer 
$n$ such that its lift on  $V_f$ joints a point $x$ with point  $f^n(x)$. 
Set $\hat{\Gamma}^s_f=p_{_f}(W^s_{\Omega_1}\setminus A_f)$ and  
$\hat{\Gamma}^u_f=p_{_f}(W^u_{\Omega_2}\setminus R_f)$. 

\begin{defi} \label{s-s-s} The collection $S_{f}=(\hat V_{f},\eta_{_{f}},\hat{\Gamma}^s_{f},\hat{\Gamma}^u_{f})$ is called by a scheme of the  diffeomorphism $f\in MS(M^3)$.
\end{defi}

\begin{defi} \label{eqg03} Shemes $S_f$ and  $S_{f'}$ of diffeomorphisms $f,f'\in MS(M^3)$ are  
called equivalent if there is a homeomorphism  $\hat\varphi:\hat V_f\to\hat V_{f'}$ with the 
following properties:

1) $\eta_{_f}=\eta_{_{f'}}\hat\varphi_*$;

2) $\hat\varphi(\hat{\Gamma}^s_{f})=\hat{\Gamma}^s_{f'}$ and  $\hat\varphi(\hat{\Gamma}^u_{f})=\hat{\Gamma}^u_{f'}$.
\end{defi}  

\begin{theo} \label{t.invariant} Morse-Smale diffeomorphisms $f,f'\in MS(M^3)$ are 
topologically conjugate if and only if their schemes are equivalent.
\end{theo}

The structure of the paper is the following:
\begin{itemize}
\item In Section \ref{I} we give informal and formal formulations of the classification results for Morse-Smale 3-diffeomorphisms and a historical background of this classification problem.
\item In Section \ref{AR} we represent the general properties of Morse-Smale diffeomorphisms and their space of wandering orbits, which  are necessary for the topological classification. 
\item In Section \ref{II} we construct a compatible system of neighbourhoods, which is a key point for the  construction of a conjugating homeomorphism. 
\item In Section \ref{III} we construct a conjugating homeomorphism.     
\item In Section \ref{IIII} we prove some topological lemmas which we used in the classification theorem.
\end{itemize}

{\it Acknowledgements.} This work was supported by the Russian Science Foundation (project 17-11-01041) ,
Russian Foundation for Basic Research (project nos. 15-01-03687-a, 16-51-10005-Ko\_a), the  Basic Research Program at the HSE (project 90) in 2017. We also thank the Institut de Math\'ematiques de Bourgogne for its warm hospitality. We are very appreciate the referee for advises which allowed us to strongly improve the introduction.

\section{General properties of Morse-Smale diffeomorphisms}\label{AR}

In this part we represent the general properties of diffeomorphisms from the class 
$MS(M^n)$ of orientation preserving  Morse-Smale diffeomorphisms $f$ given on a  smooth closed orientable n-manifolds $M^n$, which  are necessary for the topological classification. 
Proofs of facts below can be found in following sources: \cite{BGMP}, 
\cite{BoGrPo2005}, \cite{BoPa}, \cite{grin}, \cite{GrMePoZh}, \cite{GrPo2011}, \cite{Ko}, 
\cite{Pa}, \cite{PaMe1998}, \cite{Robinson-book99}, \cite{S3}, \cite{Th}. 

\subsection{Dynamics} 

Let $f\in MS(M^n)$. According to Definition \ref{oM-S}, the non-wandering set $\Omega_f$ 
of the diffeomorphism $f$ consists of finite number of the periodic points ($\Omega_f=Per_f$). The hyperbolic structure of the set $\Omega_f$ implies the existence of the invariant manifolds 
for each periodic point $p\in\Omega_f$ of period $m_p$: {\it stable} $W^s_p$ and {\it unstable} 
$W^u_p$ which are defined in topological terminus as the  following: 

$$W^s_p=\{x\in M^n~:~\lim\limits_{n\to+\infty}d(f^{nm_p}(x),p)=0\},$$ 
$$W^u_p=\{x\in M^n~:~\lim\limits_{n\to+\infty}d(f^{-nm_p}(x),p)=0\},$$ 
where $d$ is a metric on $M^n$. Moreover  $\dim~W^s_p=n-q_p$ ($\dim~W^u_p=q_p$), where 
 $q_p$ is the number of the  eigenvalues of Jacobian  $\left(\frac{\partial f^{m_p}}{\partial x}\right)\vert_{p}$
 with the absolute value greater then 1 ({\it Morse index}). In further for any subset $P\subset\Omega_f$
 will denote by $W^{u}_{P}$ ($W^s_{P}$) a union of the  unstable (stable) manifolds of all points 
from $P$. A connected component $\ell^s_p~(\ell^u_p)$ of the set $W^s_p\setminus p$ 
($\dim~W^u_p\setminus p$) is called a {\it stable (unstable) separatrix of the point $p$}. 
A number $\nu_p$, which equals $+1$ if the map $f^{m_p}|_{W^u_p}$ preserves orientation and equals $-1$ if $f^{m_p}|_{W^u_p}$ changes orientation, is called a {\it type of orientation} 
of the point $p$. 
The triple $(m_p,q_p,\nu_p)=(m_{_{\mathcal O_p}},q_{_{\mathcal O_p}},\nu_{_{\mathcal O_p}})$ 
is called a {\it periodic date} of the point $p$ or the orbit $\mathcal O_p$.  

A periodic point  $p$ is called by {\it saddle} if $0<q_p<n$ and called by {\it node} in opposite 
case, moreover $p$ is a {\it sink (source)} if $q_p=0~(q_p=n)$. As a  diffeomorphism $f$ 
preserves orientation then type orientation of any node point equals $+1$ but for a 
saddle point both values $+1,-1$ are possible.  For $q\in\{0,\dots,n\}$ denote by 
$\Omega_q$ the set of all periodic points with the Morse index $q$ and by $k_f$ 
the number of periodic orbits of  $f\in MS(M^n)$.

Dynamical properties and topological type of a Morse-Smale diffeomorphism
are largely determined by the properties of the embedding and by the mutual disposition of the invariant manifolds of the periodic points. The key role here belongs to
the study of asymptotic properties of the invariant manifolds of the saddle periodic
points.

\begin{stat} \label{M-Sm-bas} Let $f\in MS(M^n)$. Then

(1)  $M^n=\bigcup\limits_{p\in\Omega_f}W^u_p$;

(2)  $W^u_p$ is a smooth submanifold\footnote{A subset $A$ of an $C^r$-manifold $X$ ($r\geq 0$) 
is called by {\it $C^r$-submanifold} if for some integer $0\leq k\leq n$ each point of $A$ belongs 
to a chat $(U,\psi)$ of $X$ such that 
$\psi(U\cap A)=\mathbb R^k$ or $\psi(U\cap A)=\mathbb R^k_+$, where $\mathbb
R^k=\{(x_1,\dots,x_n)\in\mathbb R^n:x_{k+1}=\dots=x_n=0\}$ and $\mathbb
R^k_+=\{(x_1,\dots,x_n)\in\mathbb R^k:x_{k}\geq 0\}$. Herewith $A$ becomes $C^r$-manifolds 
with chats $\{(U\cap A,\psi\vert_{U\cap
A})\}$. $C^0$-submanifold is called also {\it topological submanifold}.

A classical topological fact says that a subset  $A$ of a $C^r$-manifold $X$ with $r\geq 1$ 
is  $C^r$-submanifold if and only if it is image of  a $C^r$-embedding. That is there is 
a $C^r$-manifold $B$ and a regular $C^r$-map $g:B\to X$ (the rank of Jacoby matrix of $g$ 
at each point equals to the dimensional of the manifold $B$) and which homeomorphically 
sends $B$ to the subspace $A=g(B)$ with induced from $X$ topology. The map $g$ is called 
{\it $C^r$-embedding}.} of the manifold $M^n$ which is diffeomorphic to $\mathbb R^{q_p}$ 
for every periodic point $p\in\Omega_f$;

(3) $cl(\ell^u_p)\setminus
(\ell^u_p\cup p)=\bigcup\limits_{r\in\Omega_f:\ell^u_p\cap W^s_r\neq\emptyset}W^u_r$ for 
every unstable separatrix $\ell^u_p~$ of a periodic point  $p\in\Omega_f$, where $cl(\cdot)$ stands for the closure of $(\cdot)$.
\end{stat}

According to item (2) of Statement  \ref{M-Sm-bas},  the map  
$f|_{W^u_{\mathcal O_p}}:W^u_{\mathcal O_p}\to W^u_{\mathcal O_p}$ is a 
diffeomorphism. Furthermore, the class of topological conjugacy of the diffeomorphism
$f^{m_{p}}|_{W^u_p}$ is completely determined by the Morse index $q_p$ and the orientation type 
$\nu_p$ of the point $p$. Namely, according to theorem on the local topological classification of the hyperbolic 
fixed point of a diffeomorphism (see theorem 5.5 in  \cite{PaMe1998}), the map $f^{m_p}$ is 
locally conjugated at $p$ to a  linear diffeomorphism   $a_{q_p,\nu_p}:\mathbb{R}^n\to\mathbb{R}^n$ 
given by the formula 
$$a_{q_p,\nu_p}(x_1,\dots,x_n)=(\nu_p\cdot 2x_1,2x_2,\dots,2x_{q_p},
\nu_p\cdot\frac{x_{q_p+1}}{2},\frac{x_{q_p+2}}{2},\dots,\frac{x_{n}}{2}).$$ 

Let us call  $a_{q,\nu}:\mathbb{R}^n\to\mathbb{R}^n$ by a  {\it canonical diffeomorphism}. 
Furthermore, denote by  $a^u_{q,\nu}$, $a^s_{q,\nu}$ the restrictions of the canonical 
diffeomorphism to $Ox_1\dots x_q$, $Ox_{q+1}\dots x_n$ and called the diffeomorphism 
$a^u_{q,\nu},a^s_{q,\nu}$ by a {\it canonical expansion}, {\it canonical contraction}, 
accordingly. According to item (2) of Statement \ref{M-Sm-bas}, $W^u_{\mathcal O_p}$ is  a 
smooth submanifold of $M^n$ and, hence, the map 
$f|_{W^u_{\mathcal O_p}}:W^u_{\mathcal O_p}\to W^u_{\mathcal O_p}$ is a diffeomorphism. 
Thus we have the following global topological classification of the maps $f|_{W^u_{\mathcal O_p}}$.
  
\begin{stat}  \label{conjtogom} Let $f\in MS(M^n)$. Then for every  periodic point 
$p\in\Omega_f$  the diffeomorphism  $f^{m_{p}}|_{W^u_p}:W^u_p\to W^u_p$ is topologically 
conjugate to the canonical expansion  $a^u_{q_p,\nu_p}:\mathbb R^{q_p}\to\mathbb R^{q_p}$.
\end{stat}

If a periodic point is saddle then the embedding of its $f$-invariant neighborhood  is also of important.  We begin with the linear
case.

For $q\in\{1,\dots,n-1\}$, $t\in(0,1]$ let   
$$\mathcal N_{q}^t=\{(x_1,\dots,x_n)\in\mathbb{R}^n~:~ (x_1^2+\dots+x_q^2)(x_{q+1}^2+\dots+x_{n}^2)< t\}$$
 and $\mathcal N^1_{q}=\mathcal N_{q}$. Notice that the set $\mathcal N_{q}^t$ is invariant with respect to the canonical diffeomorphism $a_{q,\nu}$ which has the only fixed saddle point
at the coordinate origin $O$, its unstable manifold being $W^u_O=Ox_1\dots x_q$ and its stable manifold $W^s_O=Ox_{q+1}\dots x_n$.  

\begin{defi} \label{adop} Let $f\in MS(M^n)$. We call a neighborhood $N_\sigma$ of a saddle point  
$\sigma\in\Omega_f$ linearizing if there is a homeomorphism  
${\mu}_\sigma:N_\sigma\to {\mathcal N}_{q_\sigma}$ which conjugates the diffeomorphism  
$f^{m_{\sigma}}\vert_{{N}_\sigma}$ to the canonical diffeomorphism  
$a_{q_\sigma,\nu_\sigma}|_{{\mathcal N}_{q_\sigma}}$.

The neighborhood $N_{\mathcal O_\sigma}=\bigcup\limits_{k=0}^{m_{\sigma}-1}f^k(N_\sigma)$ 
equipped with the map $\mu_{\mathcal O_\sigma}$ made up of the homeomorphisms 
$\mu_\sigma f^{-k}:f^k(N_\sigma)\to \mathcal N_{n,q_\sigma},~k=0,\dots,m_{\sigma}-1$ is 
called the  linearizing neighborhood of the orbit $\mathcal O_\sigma$ and the map $\mu_{\mathcal O_\sigma}$
 is called  linearizing.
\end{defi}

\begin{stat} \label{addo} Every saddle point (orbit) of the  diffeomorphism $f\in MS(M^n)$ has  
a linerizing neighborhood.
\end{stat}

Due to the linear dynamics near saddle points we have the following fact.

\begin{stat} Let $\sigma$ be a saddle  point of a diffeomorphism  
$f\in MS(M^n)$, let $T_\sigma\subset W^s_\sigma$ be a compact neighborhood of the point   $\sigma$ and $\xi\in T_\sigma$. Then for every sequence of points $\{\xi_m\}\subset (M^n\setminus T_\sigma)$ converging to the point $\xi$ there are a subsequence $\{\xi_{m_j}\}$, a sequence of natural numbers  
$k_{m_{j}}\to+\infty$ and a point $\eta\in (W^u_\sigma\setminus \sigma)$ such that the sequence of points  
$\{f^{k_{m_{j}}}(\xi_{m_j})\}$ converges to the point $\eta$. 
\label{sequen}
\end{stat}

\begin{figure}[h]
\centerline{\includegraphics[width=6cm,height=5.5cm]{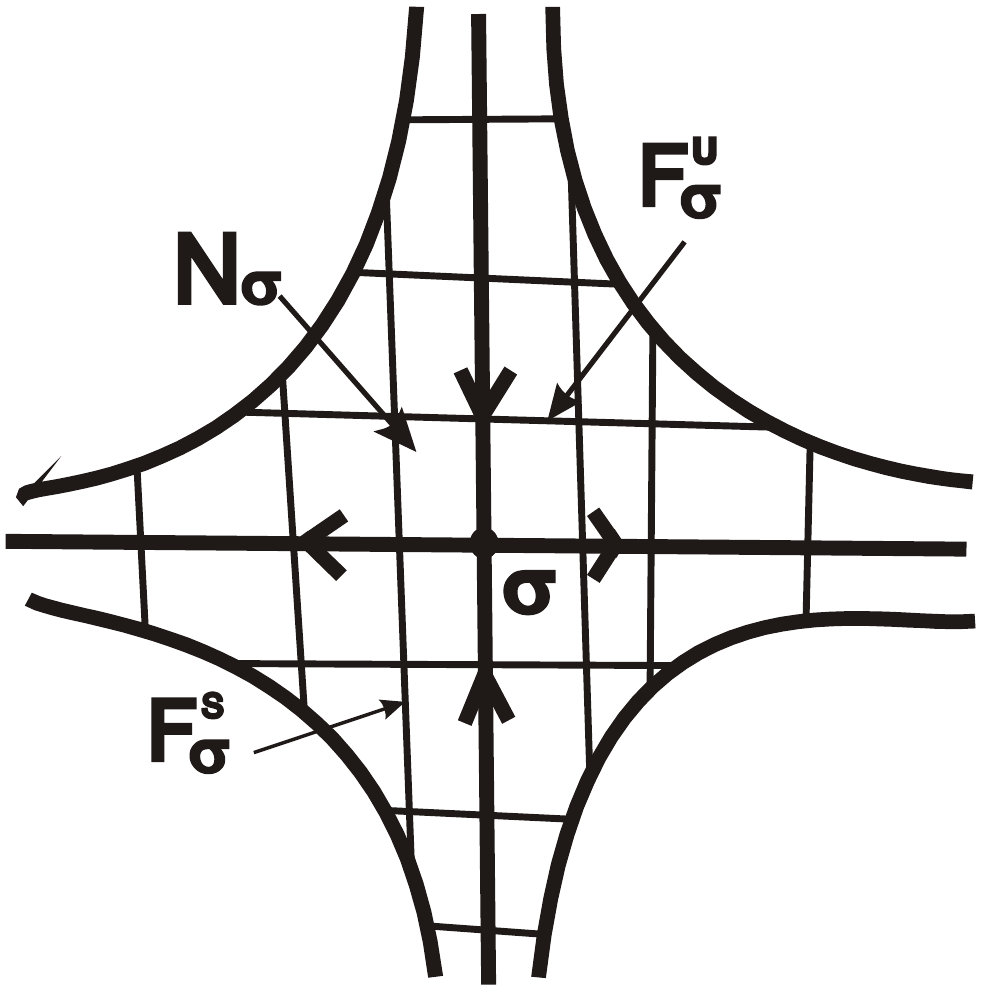}} 
\caption{\small Linearizing foliations in linearizing neighborhood}\label{5}
\end{figure}

Define in the neighborhood $\mathcal N_{q}$ a pair of transversal foliations 
$\mathcal{F}^u_q,~\mathcal{F}^s_{q}$ by next way:   
$$\mathcal{F}^u_q=\bigcup\limits_{(c_{q+1},\dots,c_n)\in Ox_{q+1}\dots x_n}
\{(x_1,\dots,x_n)\in \mathcal N_{q}~:~(x_{q+1},\dots,x_n)=(c_{q+1},\dots,c_n)\},$$ 

$$\mathcal{F}^s_{q}=\bigcup\limits_{(c_1,\dots,c_q)\in Ox_1\dots x_q}
\{(x_1,\dots,x_n)\in \mathcal N_{q}~:~(x_1,\dots,x_q)=(c_1,\dots,c_q)\}.$$
Notice that the canonical diffeomorphism $a_{q,\nu}$ sends the leaves of the 
foliation $\mathcal{F}^u_q$ ($\mathcal{F}^s_{q}$) to the leaves of the same foliation. 
According to Statement  \ref{addo}, for any saddle point $\sigma$ of $f\in MS(M^3)$, 
the foliations $\mathcal{F}^u_{q_\sigma},~\mathcal{F}^s_{q_\sigma}$ induce, by means 
linearizing map,  $f$-invariant foliations  ${F}^u_{\mathcal O_\sigma},~{F}^s_{\mathcal O_\sigma}$ 
on the linearizing neighborhood $N_{\mathcal O_\sigma}$, which are called {\it linearizing} 
(see Figure \ref{5}). 

According to item (1) of Statement  \ref{M-Sm-bas}, the invariant manifold of a periodic point of a diffeomorphism $f\in MS(M^n)$ is a  smooth submanifold of $M^n$. Nevertheless, its closure can have a complicated topological structure. The nature of this phenomenon is either dynamical or topological. The first 
case corresponds to a situation when a separatrix of a saddle point takes part in the 
heteroclinic intersections. 

\begin{defi} \label{grara} If ${\sigma_1},{\sigma_2}$ are distinct periodic saddle 
points of a diffeomorphism 
$f\in MS(M^n)$ for which  $W^s_{\sigma_1}\cap W^u_{\sigma_2}\neq\emptyset$ then the intersection 
$W^s_{\sigma_1}\cap W^u_{\sigma_2}$ is called heteroclinic. 

\begin{itemize}
\item If  $\dim(W^s_{\sigma_1}\cap W^u_{\sigma_2})>0$ then a connected component of the intersection $W^s_{\sigma_1}\cap W^u_{\sigma_2}$ is called a  heteroclinic 
manifold and if $\dim(W^s_{\sigma_1}\cap W^u_{\sigma_2})=1$ then it is called a heteroclinic curve;

\item If $\dim(W^s_{\sigma_1}\cap W^u_{\sigma_2})=0$ then the intersection  
$W^s_{\sigma_1}\cap W^u_{\sigma_2}$ is countable, each point of this set is called heteroclinic point and the orbit of a heteroclinic point  is called the heteroclinic orbit. 
\end{itemize}
\end{defi}

\begin{defi} \label{grad-like} A diffeomorphism $f\in MS(M^n)$ is said to be  gradient-like if from $W^s_{\sigma_1}\cap
W^u_{\sigma_2}\neq\emptyset$ for different points  ${\sigma_1},{\sigma_2}\in\Omega_f$ it follows that
$\dim~W^u_{\sigma_1}<\dim~W^u_{\sigma_2}$. 
\end{defi}

It follows from the transversality of intersection of invariant manifolds of the  periodic point that  a diffeomorphism $f\in MS(M^n)$ is gradient-like if and only if it has no heteroclinic points. 

According to item (3) of Statement \ref{M-Sm-bas} the closure of a separatrix of a saddle 
point which has heteroclinic intersections is not a topological manifold in general,
but the closure of a separatrix of a saddle with no heteroclinic intersections is a
topologically embedded manifold\footnote{$C^0$-map $g:B\to X$ is called a {\it topological embedding} of a topological manifold $B$ to a topological 
manifold $X$ if it homeomorphically  sends $B$ to a subspace $g(B)$ with the induced from $X$ topology. The image $A=g(B)$ is called a {\it topologically embedded manifold}. 
Notice that a topologically embedded manifold is not a topological submanifold in general. If $A$ is a  submanifold then it is called by {\it tame or tamely embedded}, 
in the opposite case by {\it wild or wildly embedded} and the points of violation of submanifold's condition are called {\it points of wildness}.}. The following 
statement holds. 

\begin{stat} \label{bezget} Let $f\in MS(M^n)$ and $\sigma$ be a saddle point 
of $f$ such that its unstable separatrix $\ell^u_\sigma$ does not take part 
in the heteroclinic intersections. Then $$cl(\ell^u_\sigma)\setminus
(\ell^u_\sigma\cup \sigma)=\{\omega\},$$ where  $\omega$ is a sink periodic 
point of $f$. Herewith, if $q_\sigma= 1$ then $cl(\ell^u_\sigma)$ is a 
topologically embedded arc in  $M^n$, if $q_\sigma\ge 2$ then   $cl(\ell^u_\sigma)$ 
is a topologically embedded $q_\sigma$-sphere in $M^n$.
\end{stat}

According to item  (2) of Statement  \ref{M-Sm-bas}, $\ell^u_\sigma\cup \sigma$ is 
a smooth submanifold of the manifold $M^n$. However the manifold $cl(\ell^u_\sigma)$ 
may be wild at the point $\omega$, in this case the separatrix $\ell^u_\sigma$ is 
called {\it wild} and it is called {\it tame} in the opposite case.

\begin{figure}[h]\centerline{\includegraphics[width=12cm,height=7cm]{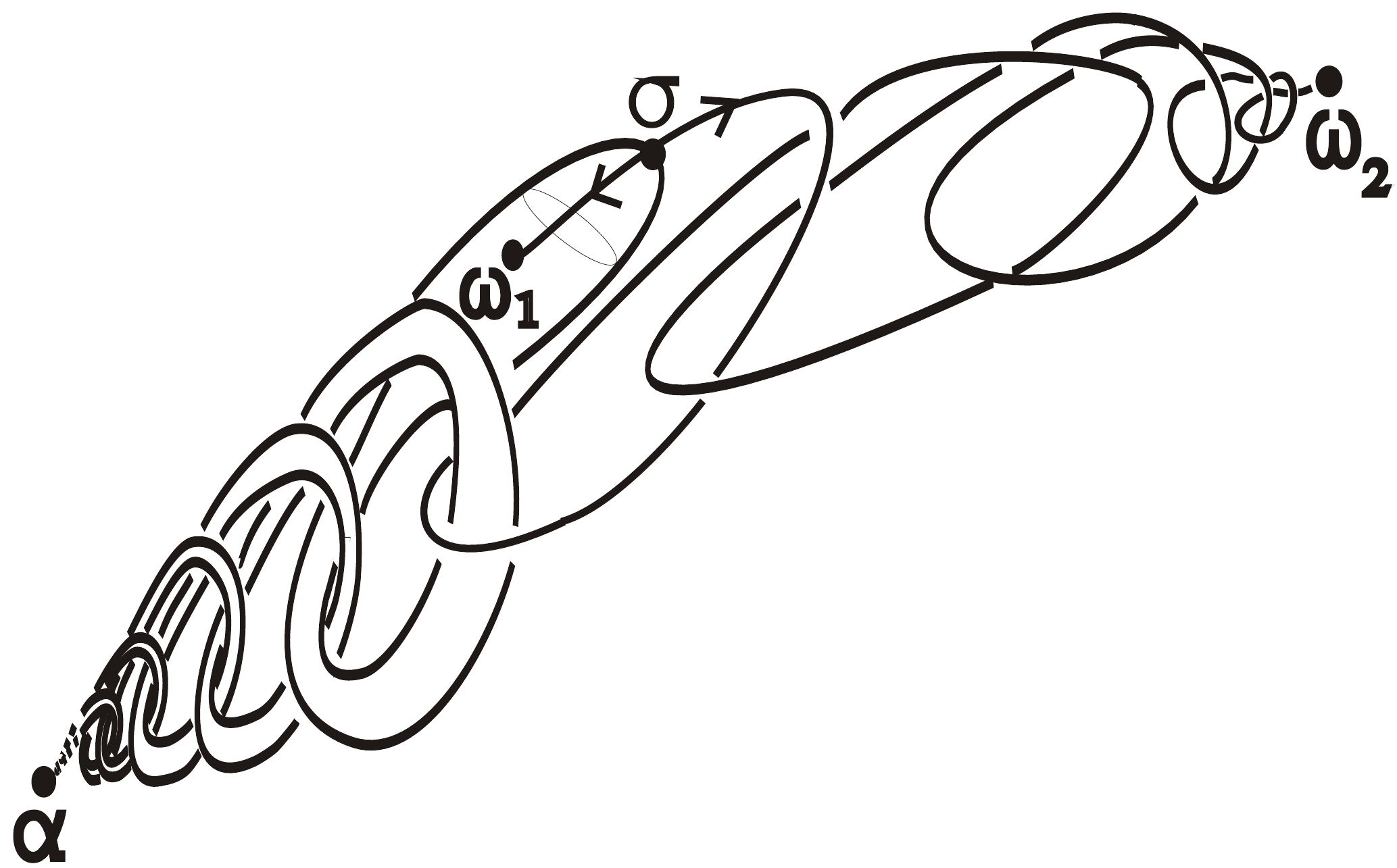}}
\caption{\small Pixton's example}\label{ex4}\end{figure}

For $n=2$ according to E. Moise's results \cite{Moise}, every compact arc and, hence, any separatrix with no heteroclinic points
is tamely embedded into  $M^2$. 
An example of a wild compact arc in $\mathbb{S}^3$ (that has nothing to do with the dynamics) was constructed by E. Artin and R. Fox  in 1948  \cite{ArFo}. The first example of a Morse-Smale diffeomorphism with wildly embedded separatrices belongs to D. Pixton 
\cite{Pi} (see Figure \ref{ex4}) and it is based on Artin-Fox arc. The following statement proved in \cite{Har} and  \cite{Bur} contains a criteria of the tame embedding of separatrices of the saddle points for a diffeomorphism $f\in MS(M^3)$.   

\begin{stat} \label{kr1} Let $f\in MS(M^3)$,  $\omega$ be a sink and $\ell^u_\sigma$ be a one-dimensional  (two-dimensional) separatrix of a saddle $\sigma$ such that  
$\ell^u_\sigma\subset W^s_\omega$. The separatrix $\ell^u_\sigma$ is tamely embedded into $M^3$ iff there is a smooth 3-ball $B_\omega\subset W^s_\omega$ containing $\omega$ 
in its interior and such that $\ell^u_\sigma$ intersects $\partial B_\omega$ at exactly 
one point (one circle).
\end{stat}

According to S. Smale \cite{S3},  it is possible to define a partial order in the set of saddle points for a Morse-Smale diffeomorphism $f$ as follows: for different periodic orbits $\mathcal O_p\neq\mathcal O_q$, one sets $$\mathcal O_p\prec\mathcal O_q~~if~~and~~only~~if~~W^u_{\mathcal O_q} \cap W^s_{\mathcal O_p} \neq\emptyset.$$ 
In that case, it follows from lemma 1.5 of  \cite{Pa} that there is a sequence of distinct periodic orbits $\mathcal O_{p_0}, \dots,\mathcal O_{p_n}$  satisfying the following conditions: $\mathcal O_{p_0}=\mathcal O_p $, $\mathcal O_{p_n}=\mathcal O_q $ and $\mathcal O_{p_i}\prec\mathcal O_{p_{i+1}}$. In that case the sequence $\mathcal O_{p_0}, \dots,\mathcal O_{p_n}$ is said to be an {\it $n$-chain  connecting $\mathcal O_p$ with $\mathcal O_q$}. The length of the longest chain connecting $\mathcal O_p,\mathcal O_q$ is denoted by $$beh(\mathcal O_q|\mathcal O_p).$$ Suppose that $beh(\mathcal O_q|\mathcal O_p)=0$ if $W^u_{\mathcal O_q} \cap W^s_{\mathcal O_p} =\emptyset$. For a subset $P$ of the periodic orbits let us set $beh(\mathcal O_q|P)=\max\limits_{\mathcal O_p\in P}\{beh(\mathcal O_q|\mathcal O_p)\}$. 

Let us divide the set of the  saddle orbits of $f$ by two parts $\Sigma_A,\Sigma_R$ such that $W^s_{\Sigma_R}\cap W^u_{\Sigma_A}=\emptyset$. Let us set 

$$A=W^u_{\Sigma_A}\cup\Omega_0,~~R=
W^s_{\Sigma_R}\cup W^s_{\Omega_3},~~V=M^n\setminus(A\cup R)~~~~~~~(*)$$

\begin{stat} \label{attris} Let $f\in MS(M^n)$ then:

1) the set $A~(R)$ is an attractor (repeller) of $f$, moreover, if 
$\dim~A\leq (n-2)$ $(\dim~R\leq (n-2))$ then the repeller $R$ 
$($attractor  $A)$ is connected and if $\dim~(A\cup R)\leq (n-2)$ 
then the manifold $V$ is connected;

2) $V=W^s_{A\cap\Omega_f}\setminus A=W^u_{R\cap\Omega_f}\setminus R$.
\end{stat}

We called $V$ by {\it characteristic manifold}. 
Below we study orbit space of some wandering sets and, in particular, 
{\it characteristic space} $\hat V=V/f$. 

\subsection{Orbit spaces}
\label{sec1}

In this section we interest the topology of an orbit space for some diffeomorphism 
$g:X\to X$ on a manifold $X$. We use denotation $X/g$ for {\it $g$-orbits on $X$} and $p_{_{X/g}}:X\to X/g$ for the natural projection. Let us recall that 
{\it a fundamental domain of action $g$ on $X$} is a closed set $D_g\subset X$ 
such that there is a set $\tilde D_{g}$ with the following properties:

1) $cl(\tilde D_{G})=D_{G}$;

2) $g^k(\tilde D_{G})\cap \tilde D_{G}=\emptyset$ for all 
$k\in(\mathbb Z\setminus\{0\})$;

3) $\bigcup\limits_{k\in\mathbb Z}g^k(\tilde D_{g})=X$. 

We say that $g$ acts {\it discontinuously} on $X$ if for each compact set  
$K\subset X$ the set of elements $k\in\mathbb Z$ such that 
$g^k(K)\cap K\neq\emptyset$ is finite.  
In the case of such action the projection  $p_{_{X/g}}$ is a cover (see Statement \ref{disc} below) 
and then we can make the following construction. Suppose that the space $X/g$ is connected and denote by $n_{_X}$ the number of connected components of $X$ and  
by $p^{-1}_{_{X/g}}(\hat x)$ the preimage of a point $\hat x\in X/g$ with respect to  
the cover $p_{_{X/g}}:X\to X/g$ (it is an orbit of some point 
${x}\in p ^ {-1} _ {_ {X/g}}(\hat x)$). Let $\hat c$ be a loop in  
$X/g$ such that $\hat c(0)=\hat c(1)=\hat x$. Due to monodromy theorem 
(see, for example, corollary 16.6 in \cite{Ko}) there is a unique path 
${c}$ in $X$ with the beginning at the point $x$ (${c} (0)=x$), 
which is the lift of $c$. Therefore, there is an element $k\in n_{_X}\mathbb Z$ 
\footnote{Here  $n_{_X}\mathbb Z$ denotes the set of integers  multiples by  $n_{_X}$.}  such that  ${c}(1)=g^k(x)$. Let $\eta_{_ {X/g}}:\pi_1(X/g)\to n_{_X}\mathbb Z$ be a 
map sending $[\hat c]$ to $k$.  
 
\begin{stat} \label{disc} Let a diffeomorphism $g$ acts discontinuously on 
$n$-manifold $X$. Then:

1) the natural projection $p_{_{X/g}}: X\to X/g$ is a cover;

2) the quotient $X/g$ is an $n$-manifold;

3) for a fundamental domain $D_g$ of action $g$ on $X$ the orbit spaces $D_g/g$ and $X/g$ are homeomorphic;

4) the map $\eta _ {_ {X/g}}:\pi_1(X/g)\to n_{_X}\mathbb Z$ is an epimorphism.
\end{stat}

\begin{stat} Let diffeomorphisms $g,~g'$ act discontinuously on manifolds $X,~X'$, 
accordingly, and $X/g,~X'/g'$ are connected. Then: 

1) if $h:X\to X$ is a homeomorphism such that  $hg=g'h$ then a map  
$\hat h: X/g\to X'/g'$ given by the formula $\hat h=p_{_{X'/g'}} hp^{-1}_{_{X/g}}$ is 
a homeomorphism and  $\eta_{_{X/g}}=\eta_{_{X'/g'}}\hat h_*$;

2) if $\hat h:X/g\to X'/g'$ is a homeomorphism such that 
$\eta_{_{X/g}}=\eta_{_{X'/g'}}\hat h_*$ then  for some $ x\in X$ and 
$x'\in p^{-1}_{_{X'/g'}}(\hat h(p_{_{X/g}}(x)))$ there is a unique homeomorphism 
 $h:X\to X'$ being a lift of $\hat h$ and such that $hg=g'h$, $h(x)=x'$. \label{conj}
\end{stat}

\begin{figure}[h]\centerline{\includegraphics[width=8.0cm,height=12.0cm]{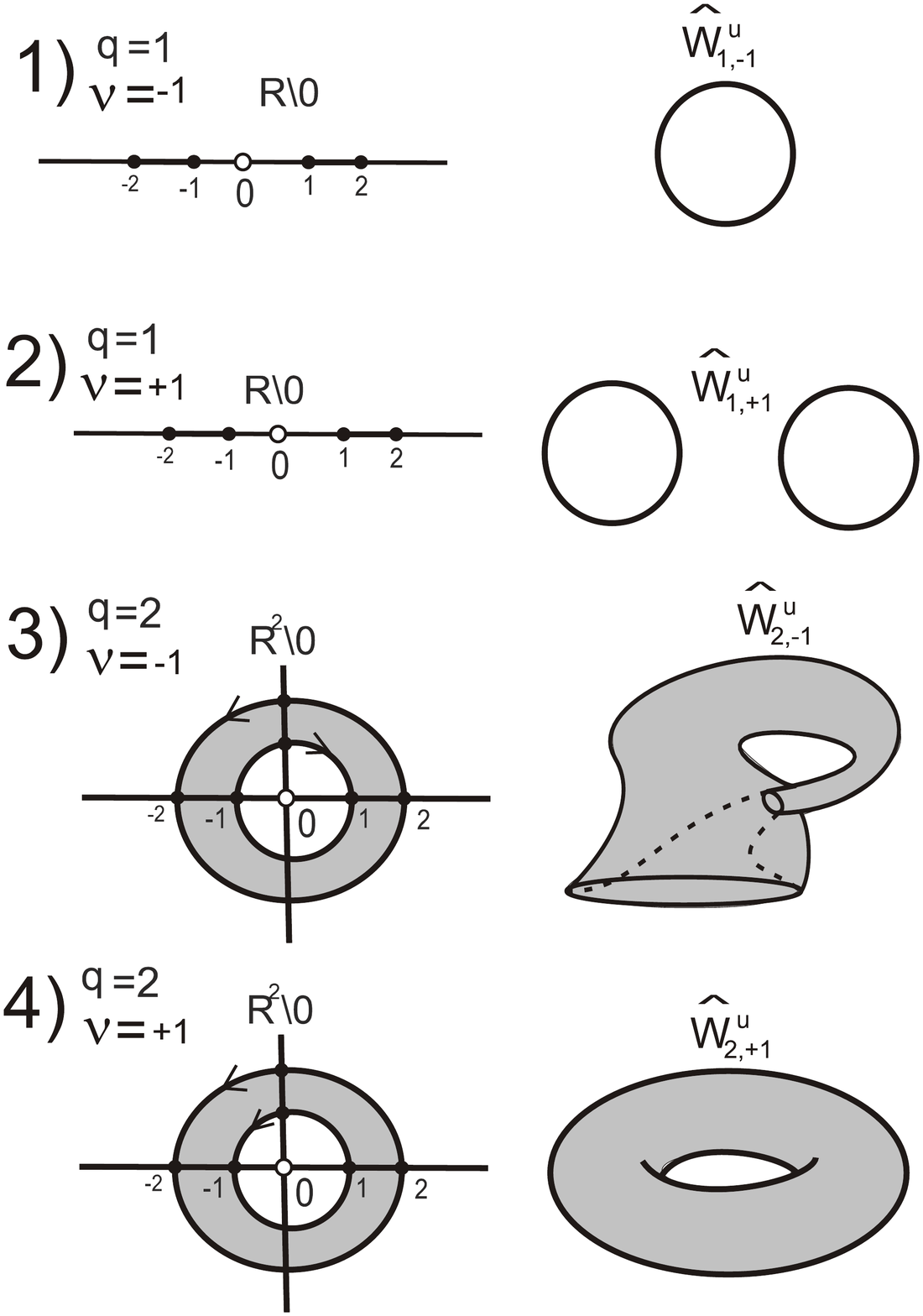}}
\caption{\small Orbit spaces of the canonical expansion}\label{factor}\end{figure}

Let us illustrate the facts above on the orbit space
 $\hat{\mathcal W}^u_{q,\nu}=(\mathbb R^{q}\setminus O)/a^u_{q,\nu}$ of 
action of the canonical expansion $a^u_{q,\nu}$ on $\mathbb R^q\setminus O$ for  
$q\in\{1,\dots,n\},~\nu\in\{+1,-1\}$. It is obvious that this action is discontinuous 
and its fundamental domain is the annulus 
$\{(x_1,\dots,x_q)\in\mathbb R^q:1\leq x_1^2+\dots+x_q^2\leq 4\}$ 
(see Figure \ref{factor}) what implies following list of possibilities.

\begin{stat}\label{++1} $ $ 

1) the space $\hat{\mathcal W}^u_{1,-1}$ is homeomorphic to the circle;

2) the space $\hat{\mathcal W}^u_{1,+1}$ is homeomorphic to the pair of circles;

3) the space $\hat{\mathcal W}^u_{2,-1}$ is homeomorphic to the Klein bottle;

4) the space $\hat{\mathcal W}^u_{2,+1}$ is homeomorphic to the torus $\mathbb T^2$;

5) the space $\hat{\mathcal W}^u_{q,-1},~q\geq 3$ is homeomorphic to a generalized  Klein bottle\footnote{{\it A generalized Klein bottle} is a topological space which is obtained from $\mathbb S^{q-1}\times[0,1]$ by identification of its boundary with 
respect to map $g:\mathbb S^{q-1}\times\{0\}\to\mathbb S^{q-1}\times\{1\}$ given by the 
formula $g(x_1,x_2,\dots,x_q,0)=(-x_1,x_2,\dots,x_q,1)$.};

6) the space $\hat{\mathcal W}^u_{q,+1},~q\geq 3$ is homeomorphic to  
$\mathbb S^{q-1}\times\mathbb S^1$.
\end{stat}

Now let $r$ be a periodic point of $f\in MS(M^n)$ with the  Morse index $q_r\geq 1$. Consider the orbit space $\hat W^u_{\mathcal O_{r}}=(W^u_{\mathcal O_{r}}\setminus\mathcal O_{r})/f$.
Next statement illustrates interrelation between $\hat W^u_{\mathcal O_{r}}$ and the  linear model. 

\begin{stat} \label{space-orbit} Let $r$ be a periodic point of a diffeomorphism $f\in MS(M^n)$ 
with the period $m_r$, the orientation type $\nu_r$ and the Morse index $q_r\geq 1$. Then the natural 
projection $p_{_{\hat W^u_{\mathcal O_{r}}}}$ is a cover which induces a structure of a smooth 
orientable $q_r$-manifold on the space $\hat W^u_{\mathcal O_{r}}$ and there is a homeomorphism 
$\hat h^u_{\mathcal O_r}:\hat W^u_{\mathcal O_r}\to\hat{\mathcal W}^u_{q_r,\nu_r}$ such that 
$\eta_{_{\hat W^u_{\mathcal O_{r}}}}([\hat c])=m_r\eta_{_{\hat{\mathcal W}^u_{q_r,\nu_r}}}([\hat h^u_{\mathcal O_r}(\hat c)])$ for every closed curve $\hat c\subset\hat W^u_{\mathcal O_{r}}$. 
\end{stat}

It is similar for a definition of the space orbit  
$\hat{\mathcal W}^s_{q,\nu}=(\mathbb R^{n-q}\setminus O)/a^{s}_{q,\nu}$ of 
{\it canonical contraction} for $q\in\{0,\dots,n-1\},~\nu\in\{+1,-1\}$ and 
the space $\hat W^s_{\mathcal O_{r}}=(W^s_{\mathcal O_{r}}\setminus\mathcal O_{r})/f$ 
for periodic point $r$ with the Morse index $q_r\leq (n-1)$.

\begin{figure}[h]\centerline{\includegraphics[width=9.5cm,height=13cm]{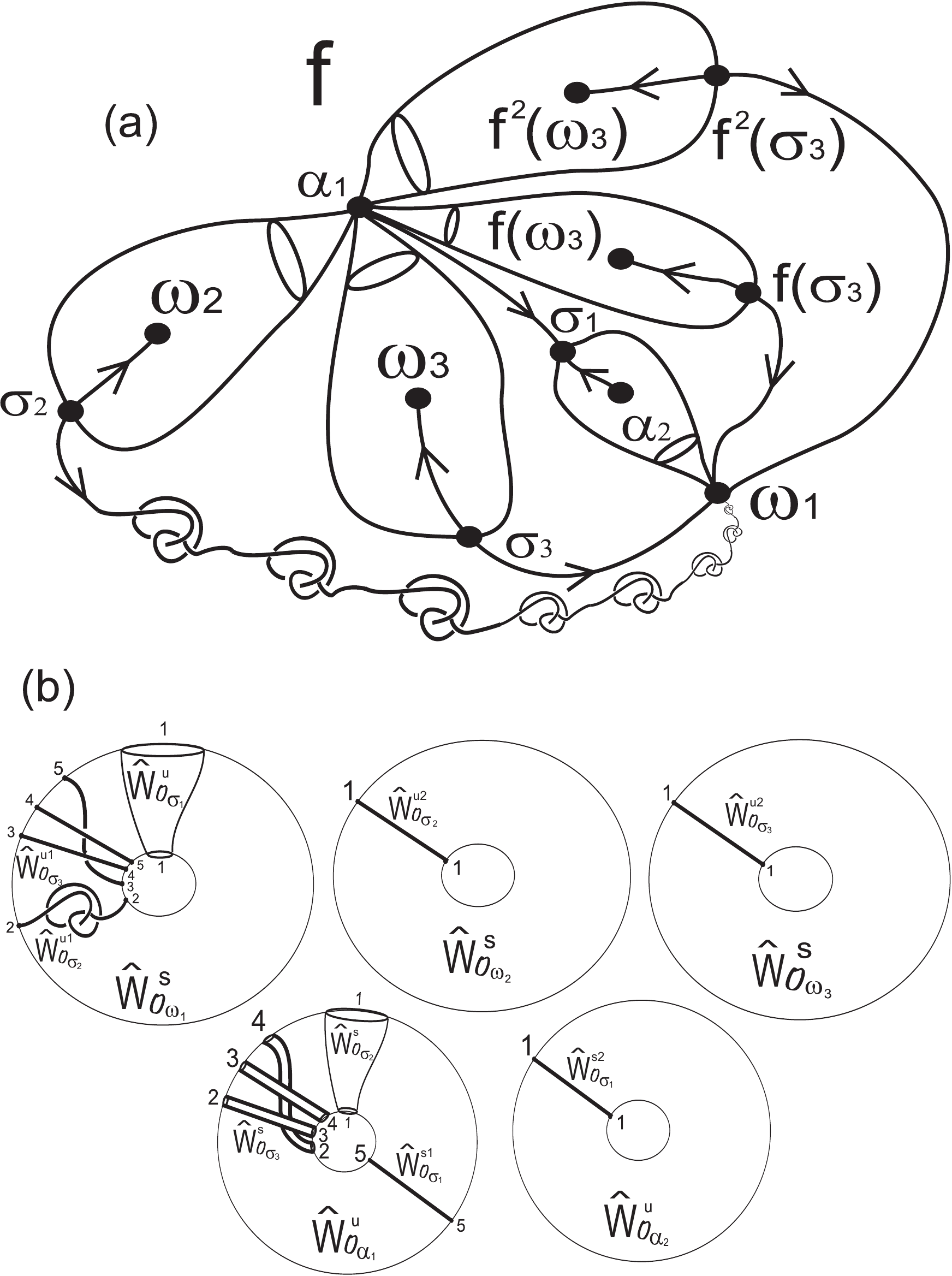}}
\caption{\small Orbit spaces of the separatrices of the periodic points}
\label{p8}\end{figure}

Figure \ref{p8} (a) shows a Morse-Smale diffeomorphism $f:\mathbb S^3\to\mathbb S^3$ the non-wandering set of which consists of eight periodic points with the following periodic data:  $\mathcal O_{\omega_1}(1,0,+1)$, $\mathcal O_{\omega_2}(1,0,+1)$, $\mathcal O_{\omega_3}(3,0,+1)$, $\mathcal O_{\sigma_1}(1,2,+1)$, $\mathcal O_{\sigma_2}(1,1,+1)$, $\mathcal O_{\sigma_3}(3,1,+1)$, $\mathcal O_{\alpha_1}(1,3,+1)$, $\mathcal O_{\alpha_2}(1,3,+1)$. Figure \ref{p8} (b) shows the fundamental domains of the action of the diffeomorphism $f$ on ${W^s_{\mathcal O_{\omega_i}}\setminus\mathcal O_{\omega_i}},~i=1,2,3$, ${W^u_{\mathcal O_{\alpha_i}}\setminus\mathcal O_{\alpha_i}},~i=1,2$. Each fundamental domain is the 3-annulus from which the orbits spaces  $\hat W^s_{\mathcal O_{\omega_i}},~i=1,2,3$, $\hat W^u_{\mathcal O_{\alpha_i}},~i=1,2$ are obtained by gluing the boundary spheres by the diffeomorphism  $f^{m_{\omega_i}},~i=1,2,3$, $f^{m_{\alpha_i}},~i=1,2$ respectively. The orbits spaces $\hat W^s_{\mathcal O_{\sigma_i}},~\hat W^u_{\mathcal O_{\sigma_i}},~i=1,2,3$ are obtained from the arcs and the cylinders by gluing the points with the same numbers and of the circles with the same numbers.

\begin{figure}[h]\centerline{\includegraphics[width=9cm,height=10cm]{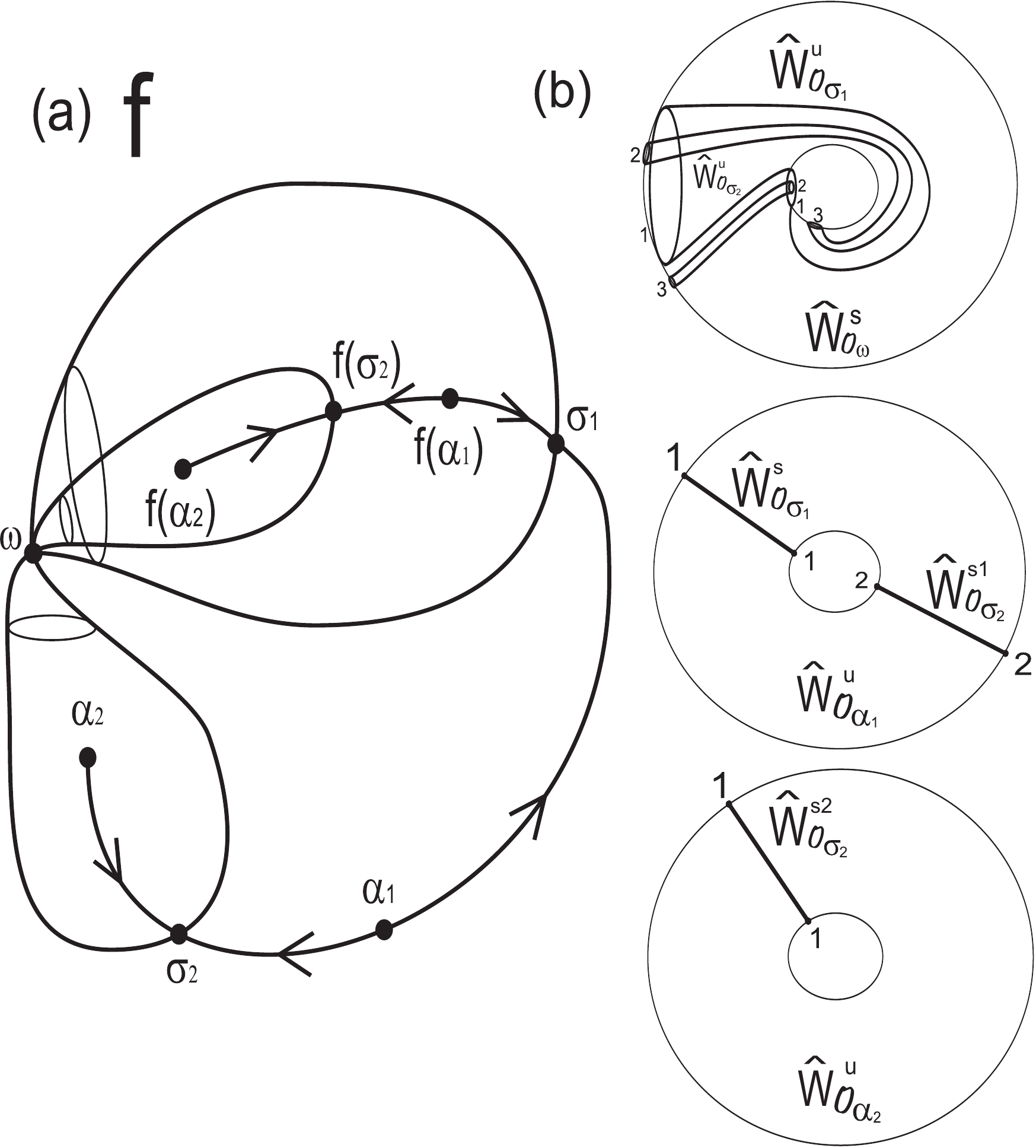}}
\caption{\small Orbit spaces of the separatrices of the periodic points}
\label{8888}\end{figure}

Figure \ref{8888} (a) shows a Morse-Smale diffeomorphism $f:\mathbb S^3\to\mathbb S^3$ the non-wandering set of which consists of five periodic points with the following periodic data:  $\mathcal O_{\omega}(1,0,+1)$, $\mathcal O_{\sigma_1}(1,2,-1)$, $\mathcal O_{\sigma_2}(2,2,+1)$, $\mathcal O_{\alpha_1}(2,3,+1)$, $\mathcal O_{\alpha_2}(2,3,+1)$. Figure \ref{8888} (b) shows the fundamental domains of the action of the diffeomorphism $f$ on ${W^s_{\mathcal O_\omega}\setminus\mathcal O_\omega}$ and ${W^u_{\mathcal O_{\alpha_i}}\setminus\mathcal O_{\alpha_i}},~i=1,2$. Each fundamental domain is the 3-annulus from which the orbits spaces $\hat W^s_{\mathcal O_\omega}$,  $\hat W^u_{\mathcal O_{\alpha_i}},~i=1,2$ are obtained by gluing the boundary spheres of the annulus by the the diffeomorphism $f^{m_{\omega}}$, $f^{m_{\alpha_i}},~i=1,2$ respectively. The orbits spaces $\hat W^s_{\mathcal O_{\sigma_i}},~\hat W^u_{\mathcal O_{\sigma_i}},~i=1,2$ are obtained from the arcs and the cylinders by gluing the points with the same numbers and the circles with the same numbers.

On the set ${\mathcal{N}}^{u}_{q}=\mathcal{N}_{q}\setminus W^s_O$ the action of the group $A_{q,\nu}=\{a_{q,\nu}^k,\,k\in\mathbb Z\}$ is discontinuous again. Due to Statement \ref{disc} the space orbit $\hat{\mathcal{N}}^u_{q,\nu}=
(\mathcal{N}^{u}_{q})/a_{q,\nu}$ is a smooth  $n$-manifold. As  
$a_{q,\nu}|_{W^u_O\setminus O}=a^u_{q,\nu}|_{W^u_O\setminus O}$ then   
$\hat{\mathcal{N}}^u_{q,\nu}$ is tubular neighborhood of the space  
$\hat{\mathcal{W}}^u_{q,\nu}$. Furthermore  $\hat{\mathcal{W}}^u_{q,+1}$ is homeomorphic to $\mathbb S^{q-1}\times\mathbb S^1\times\{0\}$ and its tubular neighborhood  
$\hat{\mathcal{N}}^u_{q,+1}$ is homeomorphic to  
$\mathbb S^{q-1}\times\mathbb S^1\times\mathbb D^{n-q}$. As $a^2_{q,-1}=a^2_{q,+1}$ 
and the diffeomorphisms $a^2_{q,+1}$ and $a_{q,+1}$ are topologically conjugated due to  
Statement  \ref{conjtogom}, hence, the manifold $\hat{\mathcal{W}}^u_{q,+1}$ is the 2-fold cover for the manifold $\hat{\mathcal{W}}^u_{q,-1}$  and the manifold $\hat{\mathcal{N}}^u_{q,+1}$ is the 2-fold cover for the neighborhood  $\hat{\mathcal{N}}^u_{q,-1}$.  

Similarly one defines the orbits space $\hat{\mathcal N}^s_{q,\nu}=\mathcal{N}^{s}_{q}/a^{s}_{q,\nu}$ (where ${\mathcal{N}}^{s}_{q}=\mathcal{N}_{q}\setminus W^u_O$), the covering map $p_{_{\hat{\mathcal{N}}^s_{q,\nu}}}:\mathcal{N}^{s}_{q}\to\hat{\mathcal{N}}^s_{q,\nu}$ and the map $\eta_{_{\hat{\mathcal{N}}^s_{q,\nu}}}$ from the union of the fundamental groups of the connected components of the manifold $\hat{\mathcal{N}}^s_{q,\nu}$ into the group $\mathbb Z$.

\begin{figure}[h]\centerline{\includegraphics[width=10cm,height=11cm]{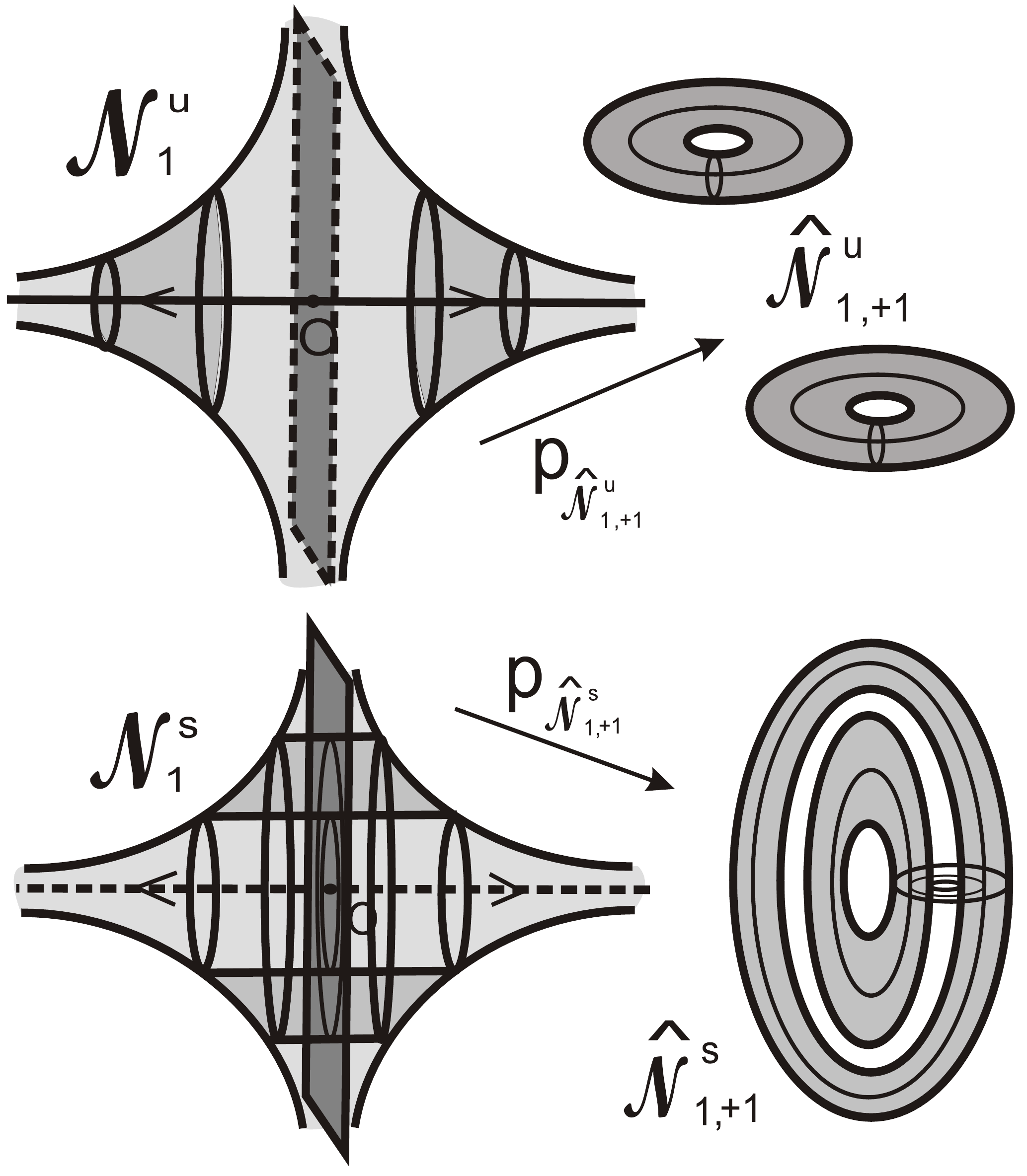}}
\caption{\small Neighborhoods of the orbit spaces of the canonical contraction and the expansion for 
$n=3$}\label{perestt3}\end{figure}

Figure \ref{perestt3} shows these objects for $n=3$; $q=1$; $\nu=+1$. To make the structure of the orbits space $\hat{\mathcal{N}}^s_{q,\nu},\hat{\mathcal{N}}^u_{q,\nu}$ more clear we mark out the fundamental domain of the action of the canonical diffeomorphism $a_{q,\nu}$ on the sets ${\mathcal{N}}^s_{q},{\mathcal{N}}^u_{q}$.

Now let $\sigma$ be a saddle periodic point with Morse index $q_\sigma$ of a diffeomorphism $f\in MS(M^n)$ and let $N_{\mathcal O_\sigma}$ be a linearizing neighborhood of the orbit $\mathcal O_\sigma$. Denote $N^u_{\mathcal O_\sigma}=N_{\mathcal O_\sigma}\setminus W^s_{\mathcal O_\sigma}$. Consider the orbits space $\hat N^u_{\mathcal O_{\sigma}}=N^u_{\mathcal O_{\sigma}}/f$ of the action of the diffeomorphism $f$ on $N^u_{\mathcal O_{\sigma}}$.  Denote by $p_{_{\hat N^u_{\mathcal O_{\sigma}}}}:N^u_{\mathcal O_{\sigma}}\to\hat N^u_{\mathcal O_{\sigma}}$ the natural projection. The following statement shows the connection between the orbits space $\hat N^u_{\mathcal O_{\sigma}}$ and the linear model.
\begin{stat} \label{space-orbit-n} Let $\sigma$ be a saddle periodic point of period $m_\sigma$ with orientation type $\nu_\sigma$ and Morse index $q_\sigma$ for a diffeomorphism $f\in MS(M^n)$. Then the projection $p_{_{\hat N^u_{\mathcal O_{\sigma}}}}$  is the covering map; it induces a structure of a smooth $n$-manifold on the orbits space $\hat N^u_{\mathcal O_{\sigma}}$ and it induces a map $\eta_{_{\hat N^u_{\mathcal O_{\sigma}}}}$  from the union of the fundamental groups of the connected components of the manifold $\hat N^u_{\mathcal O_{\sigma}}$ into the group $\mathbb Z$ such that there is a homeomorphism  $\hat\mu^u_{\mathcal O_\sigma}:\hat N^u_{\mathcal O_\sigma}\to\hat{\mathcal N}^u_{q_\sigma,\nu_\sigma}$ which satisfies $\eta_{_{\hat N^u_{\mathcal O_{\sigma}}}}([\hat c])=m_\sigma\eta_{_{\hat{\mathcal N}^u_{q_\sigma,\nu_\sigma}}}([\hat\mu^u_{\mathcal O_\sigma}(\hat c)])$ for any closed curve $\hat c\subset\hat N^u_{\mathcal O_{\sigma}}$. 
\end{stat}

Similarly one defines the orbits space $\hat N^s_{\mathcal O_{\sigma}}=N^s_{\mathcal O_{\sigma}}/f$ of the action of the group $F$ on $N^s_{\mathcal O_{\sigma}}=N_{\mathcal O_{\sigma}}\setminus W^u_{\mathcal O_{\sigma}}$,  the covering map $p_{_{\hat{N}^s_{\mathcal O_{\sigma}}}}:{N}^{s}_{\mathcal O_{\sigma}}\to\hat{{N}}^s_{\mathcal O_{\sigma}}$ and the map $\eta_{_{\hat{{N}}^s_{\mathcal O_{\sigma}}}}$ consisting of nontrivial homomorphisms into the group $\mathbb Z$ on the fundamental group of each connected component of the manifold  $\hat{{N}}^s_{\mathcal O_{\sigma}}$. 

Below for any $t\in(0,1]$ we denote  $N^t_{{\sigma}}=(\mu_{\sigma})^{-1}({\mathcal N}^t_{q_\sigma})$, $N^t_{\mathcal O_{\sigma}}=\bigcup\limits_{k=0}^{m_{\sigma}-1}f^k(N^t_\sigma)$, ${\mathcal N}^{ut}_{q_\sigma}={\mathcal N}^t_{q_\sigma,\nu_\sigma}\setminus W^s_O$, $N^{ut}_{{\sigma}}=(\mu_{\sigma})^{-1}({\mathcal N}^{ut}_{q_\sigma})$, $N^{ut}_{\mathcal O_{\sigma}}=\bigcup\limits_{k=0}^{m_{\sigma}-1}f^k(N^{ut}_\sigma)$, ${\mathcal N}^{st}_{q_\sigma}={\mathcal N}^t_{q_\sigma}\setminus W^u_O$, $N^{st}_{{\sigma}}=(\mu_{\sigma})^{-1}({\mathcal N}^{st}_{q_\sigma})$, $N^{st}_{\mathcal O_{\sigma}}=\bigcup\limits_{k=0}^{m_{\sigma}-1}f^k(N^{st}_\sigma)$, $\hat{\mathcal N}^{ut}_{q_\sigma,\nu_\sigma}=p_{_{\hat{\mathcal N}^{ut}_{q_\sigma,\nu_\sigma}}}({\mathcal N}^{ut}_{q_\sigma})$ and $\hat{\mathcal N}^{st}_{q_\sigma,\nu_\sigma}=p_{_{\hat{\mathcal N}^{st}_{q_\sigma,\nu_\sigma}}}({\mathcal N}^{st}_{q_\sigma})$.

\begin{stat} \label{Ge} For every $t\in (0,1)$, the  neighborhood $N^t_\sigma$ is linearizing.
Generically, the boundary of $N^t_\sigma$ does not contain any heteroclinic point.
\end{stat}

Let us recall that we divided the set of the saddle orbits of $f$ by two parts $\Sigma_A,\Sigma_R$ such that $W^s_{\Sigma_R}\cap W^u_{\Sigma_A}=\emptyset$ and set $A=W^u_{\Sigma_A}\cup\Omega_0,~~R=
W^s_{\Sigma_R}\cup \Omega_3,~~V=M^n\setminus(A\cup R)$.

\begin{stat} The orbit space $\hat V=V/f$ is a closed orientable n-manifold. 
\end{stat}

\section{Compatible foliations}\label{II}

In this section for any diffeomorphism $f\in MS(M^3)$ the existence of a compatible system of neighbourhoods is proved. This system is a key point for the  construction of conjugating homeomorphism. 
 
\begin{defi}\label{dopsystem} Let $f\in MS(M^3)$. A collection $N_f$ of linearizing 
neighbourhoods ${N}_{\mathcal O_1},\dots,{N}_{\mathcal O_{k_f}}$ of all saddle orbits ${\mathcal O_1},\dots,{\mathcal O_{k_f}}$ of  $f$ is called compatible and the corresponding foliations are called compatible if for every saddle orbits $\mathcal O_i,~\mathcal O_j$ following properties are hold:

1) if $W^{s}_{\mathcal O_i}\cap W^{u}_{\mathcal O_j}=\emptyset$ and $W^{u}_{\mathcal O_i}\cap W^{s}_{\mathcal O_j}=\emptyset$  then ${N}_{\mathcal O_i}\cap{N}_{\mathcal O_j}=\emptyset$;

2) if $W^{s}_{\mathcal O_i}\cap W^{u}_{\mathcal O_j}\neq\emptyset$ then $({F}^s_{\mathcal O_j,x}
\cap{N}_{\mathcal O_i})\subset{F}^s_{\mathcal O_i,x}$ and $({F}^u_{\mathcal O_i,x}\cap
{N}_{\mathcal O_j})\subset{F}^u_{\mathcal O_j,x}$ for  $x\in(N_{\mathcal O_i}\cap N_{\mathcal O_j})$.
\end{defi}

\begin{rema} The compatible system of neighborhoods is a modification of the  admissible system of the tubular families introduced by J. Palis and S. Smale in papers \cite{Pa} and \cite{PS}. But we construct a compatible system of neighborhoods for arbitrary diffeomorphism $f\in MS(M^3)$ independently (see Theorem \ref{use2-1/2}).
\end{rema}

\begin{figure}[h]
\centerline{\includegraphics[width=9cm,height=13cm]{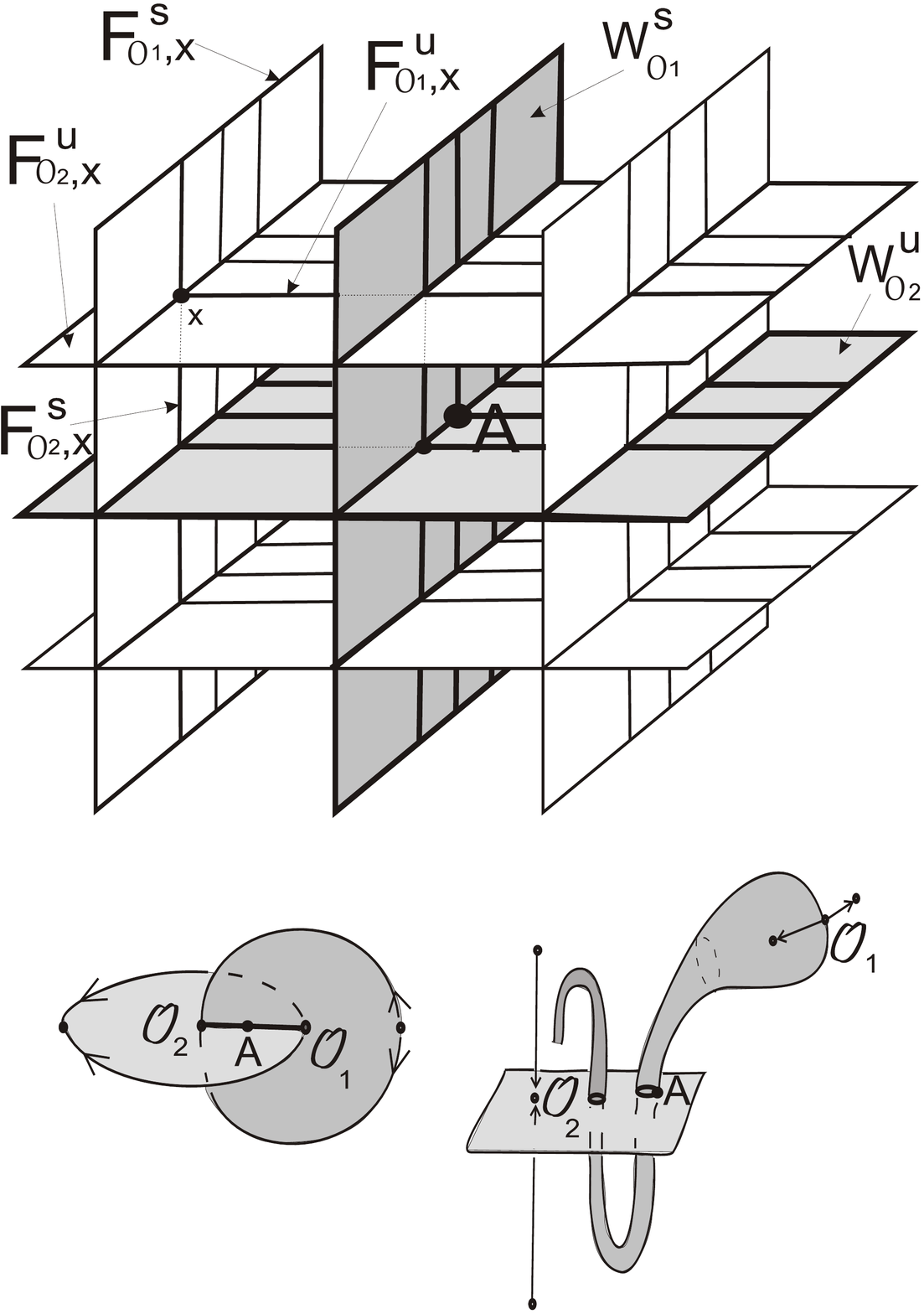}} \caption{\small Foliated 
neighborhood of a point on a heteroclinic curve}\label{dopusti}
\end{figure}

Figure \ref{dopusti} shows a foliated neighborhood of a point $A$ belonging to a heteroclinic curve. Below some Morse-Smale diffeomorphisms with heteroclinic curves on $\mathbb S^3$ are represented.

Let $f\in MS(M^3)$. It follows from Statement \ref{attris} that if the set $\Omega_2$ is empty then $R_f$ consists of unique source. If $\Omega_2\neq\emptyset$ then denote by $n$ the length of the longest chain connecting some $p,q$ from $\Omega_2$. Divide the  set $\Omega_2$ by $f$-invariant disjoint parts $\Sigma_0,\Sigma_1,\dots,\Sigma_n$ using the rule: $beh(\mathcal O|(\Omega_2\setminus\mathcal O))=0$ for each orbit $\mathcal O\in \Sigma_{0}$ and $beh(\mathcal O|\Sigma_i)=1$ for each orbit $\mathcal O\in \Sigma_{i+1},~i\in\{0,\dots,n-1\}$.     

Let us introduce the following notations.
\begin{den}\label{ate} $ $

\begin{itemize}
\item[-] $W^u_i:=W^u_{\Sigma_i}$, $W^s_i:=W^s_{\Sigma_i}$; 
\item[-] $N^t_i:=\bigcup\limits_{\mathcal O\in\Sigma_i}N^t_\mathcal O$ and $\mu_i$ composed from $\mu_{\mathcal O},~\mathcal O\in\Sigma_i$;
\item[-] for every point $x\in N_i$, denote  ${F}^u_{i,x}$ (resp. ${F}^s_{i,x}$) the  leave of the foliation 
 ${F}^u_{i}$ (resp. ${F}^s_{i}$) passing through $x$;
\item[-] for every point $x\in N_i$, set $x^u_i=W^u_i\cap F^s_{i,x}$ and $x^s_i=W^s_i\cap F^u_{i,x}$. Thus, we
have $x=(x^u_i,x^s_i)$ in the coordinates defined by $\mu_i$;
\item[-]  for $i\in\{0,\dots,n\}$, set
$A_{i}:=A_f\cup\bigcup\limits_{j=0}^i W_j^u,\hfill\break
V_i:= W^s_{A_i\cap\Om_f}\smallsetminus A_i,~~\hat V_i:= V_i/f$. Observe that $f$ acts freely 
on $V_i $ and denote the natural projection by $p_{_{i}}:V_{i}\to\hat V_{i}$ and $\eta_i$ corresponding epimorphism;
\item[-]  for complicity we set $A_{-1}=A_f,~R_{-1}=R_f$, $V_{-1}=V_f,~\hat V_{-1}=\hat V_f,~p_{_{-1}}=p_{_f},~\eta_{-1}=\eta_f$; 
\item[-] for $j,k\in\{0,\dots,n\}$ and $t\in(0,1)$, set $\hat W^s_{j,k}=p_k(W^{s}_{j}\cap V_k)$, $\hat W^u_{j,k}=p_k(W^{u}_{j}\cap V_k)$,  $\hat{N}^{t}_{j,k}=p_{k}(N^t_j\cap V_{k})$; 
\item[-]  $L^u:=\bigcup\limits_{i=0}^nW^u_i,~L^s:=\bigcup\limits_{i=0}^nW^s_i$,
$L^u_i:=L^u\cap V_i$\,,
~$L^s_i:= L^s\cap V_i$\,,
$\hat{L}^u_i:=p_i(L^u_i)$, $\hat{L}^s_i:=p_i(L^s_i)$. 
\end{itemize}
\end{den}
As $\Omega_1$ for $f$ is $\Omega_2$ for $f^{-1}$ then it is possible to divide the periodic orbits of the set $\Omega_1$ by parts by a similar way.

\begin{theo} \label{use2-1/2} For each diffeomorphism $f\in MS(M^3)$ there is a 
compatible system of neighbourhoods. 
\end{theo}
\begin{demo} The proof consists of four steps. 

\nd {\bf Step 1.} Due to Lemma 3.5 in \cite{BoGrLaPo} there exist $f$-invariant neighbourhoods $U^{s}_{0},\dots,U^{s}_{n}$ of 
the sets $\Sigma_{0},\dots,\Sigma_{n}$ respectively,  equipped with two-dimensional $f$-invariant  foliations   ${F}^{u}_{0},\dots,{F}^{u}_{n}$ whose leaves are smooth
such that  the following properties hold for each $i\in\{0,\dots,n\}$:

{\rm (i)} the unstable manifolds $W^u_i$ are leaves of the foliation ${F}^{u}_{i}$ and 
each leaf of the foliation ${F}^{u}_{i}$ is transverse to  $L^s_i$;
 
{\rm (ii)} for any  $0\leq i <k\leq n$ and $x\in U^s_{i}\cap U^s_{k}$, we have the inclusion 
$F^u_{k,x}\cap {U}^s_{i}\subset F^u_{i,x}$. 

Denote by $F^u_{\Omega_2}$ constructed two-dimensional foliation. We have similarly two-dimensional foliation $F^s_{\Omega_1}$.     
Set $$\hat F^u_{\Omega_2}=p_{_f}(F^u_{\Omega_2}),~~~\hat F^s_{\Omega_1}=p_{_f}(F^s_{\Omega_1}).$$

{\bf Step 2.} Let us construct an  $f$-invariant neighbourhood $U(H)$ of the set of the heteroclinic curves $H=W^s_{\Omega_1}\cap W^u_{\Omega_2}$ equipped by an $f$-invariant  
foliation $G$ consisting of two-dimensional discs which are transversal to $H$ and to both foliations $F^u_{\Omega_2}$ and $F^s_{\Omega_1}$. 

Set  $\hat H=\hat {L}^{s}_{f}\cap\hat {L}^{u}_{f}$. By the construction the set $\hat H$ is compact and  consists of at most countable set of curves being the projection of the heteroclinic curves. Let us divide the set $\hat H$ by parts $\hat H_0,\dots,\hat H_m$ next way: $\hat H_0$ consists of all compact curves and $\hat H_{i+1}$ consists of curves $\hat\gamma$ such that $cl~(\hat\gamma)\setminus\hat\gamma\subset cl~(\hat H_i)$ for $i\in\{0,\dots,m-1\}$.   

Any curve $\hat\gamma\subset\hat H$ belongs to the intersection 
$p_{_f}({W}^{s}_{\sigma_1})\cap p_{_f}(W^{u}_{\sigma_2})$ for some $\sigma_1\in\Omega_1,~\sigma_2\in\Omega_2$  
(depending on $\hat\gamma$). Then there is a tubular neighbourhood 
$U(\hat \gamma)$ of the curve $\hat\gamma$ foliated by two-dimensional discs 
$\hat G_{\hat\gamma}=\{\hat d_x,~x\in\hat \gamma\}$ which are 
transversal to the leaves of the foliations $\hat{ F}^{s}_{\Omega_1}$, $\hat{ F}^{u}_{\Omega_2}$. Denote by $\hat G_{\hat\gamma,\hat x}$ a leaf of the foliation $\hat G_{\hat\gamma}$ passing through a point $\hat x$. Due to the  compatibility of two-dimensional foliations  ${F}^{s}_{\Omega_1}$, ${F}^{u}_{\Omega_2}$ we can construct these tubular neighbourhoods with discs satisfying the properties: if $U(\hat \gamma_i)\cap U(\hat \gamma_j)\neq\emptyset$ for $\hat\gamma_i\subset\hat H_i,~\hat\gamma_i\subset\hat H_i,~0\leq i<j\leq m$ then $$(\hat G_{\hat\gamma_j,\hat x}
\cap{U}_{\hat\gamma_i})\subset\hat G_{\hat\gamma_i,\hat x}~~for~~\hat x\in({U}_{\hat\gamma_i}\cap{U}_{\hat\gamma_j}).$$

Denote by $\hat G$ a two-dimensional foliation  which is formed by the two-dimensional 
discs of the foliations $\hat G_{\hat\gamma},\hat\gamma\subset\hat H$ 
and by  $$\hat{F}^s_{\hat H}~(\hat{F}^{u}_{\hat H})$$ a one-dimensional foliation which is formed
by the intersection of the leaves of the foliation $\hat G$ with the leaves  of the foliations 
$\hat F^s_{\Omega_1}~(\hat F^u_{\Omega_2})$. Set $U(\hat H)=\bigcup\limits_{\hat\gamma\subset\hat H}U(\hat\gamma)$, $U(H)=p_{_{k_1}}^{-1}(U(\hat H))$ and denote by  ${F}^s_{H}~({F}^u_{H})$ a foliation on $U(H)$ consisting from the preimages with respect to the projection $p_{_{f}}$ of leaves of the foliation 
$\hat{F}^s_{\hat H}~(\hat{F}^u_{\hat H})$. Without loss of generality we can suppose that the 
following projection $\hat\pi^{s}_{\hat H}:
{U(\hat H)}\to\hat {L}^{u}_{f}~(\hat\pi^{u}_{\hat H}:
{U(\hat H)}\to\hat {L}^{s}_{f})$ along the leaves $\hat{F}^s_{\hat H}~(\hat{ F}^{u}_{\hat H})$
is well-defined.

We also have the following statement.

\begin{lemm} \label{la} There exist  $f$-invariant neighborhoods $U^{u}_{0},\dots,U^{u}_{n}$ of  the sets $\Sigma_{0},\dots,\Sigma_{n}$  respectively, equipped with one-dimensional $f$-invariant foliations ${F}^{s}_{0},\dots,{F}^{s}_{n}$ with smooth leaves such that the  following properties  hold for each $i\in\{0,\dots,n\}$:

{\rm (iii)} the stable manifold $W^s_i$ is a leaf of the foliation ${F}^{s}_{i}$ and each leaf of the foliation ${F}^{s}_{i}$ is transverse to $L^u_i$;
 
{\rm (iv)} for any $0\leq j<i$ and $x\in U^u_{i}\cap U^u_{j}$, we have the inclusion 
$({F}^s_{j,x} \cap{U}^u_{i})\subset
{F}^s_{i,x}$;

{\rm (v)} the intersection of a leave of the foliation ${F}^{s}_{i}$ with the set 
$U(H)$ is a leave of the foliation ${F}^s_{H}$.
\end{lemm}
\begin{demo} The proof is done  by an increasing induction  from $i=0$; it is skipped due to similarity to the Step 1.
\end{demo}

{\bf Step 3.} We prove the following statement for each $i=0,\dots,n$.
\begin{lemm} \label{lal} $ $

{\rm (vi)} There exists  an $f$-invariant neighborhood $\tilde N_{i}$ of the set $\Sigma_i$ contained in $U^s_{i}\cap U^u_{i}$ and such that the restrictions of the foliations ${F}^u_i$ and ${F}^s_i$ to $\tilde N_i$ are transverse.
\end{lemm}
\begin{demo} For this aim, let us choose a fundamental domain   $K^s_i$ of the restriction of $f$ to $W^s_i\setminus\Sigma_i$ and take  a tubular neighborhood $N(K^s_i)$ of  $K^s_i$ whose disc fibres are contained in leaves of $F^u_i$. Due to property (i), $F^u_i$ is transverse to $W^s_i$. Since $F^s_i$ is a $C^1$-foliation, if $N(K^s_i)$ is small enough, the foliations $F^s_i$ and $F^u_i$ have transverse intersection in $N(K^s_i)$. Set  $$\tilde N_i:=W^u_i\,\bigcup_{k\in\mathbb Z}f^k\left(N(K^s_i)\right). $$ 
This is a neighbourhood  of $\si_i$; it satisfies condition (vi) and the previous properties (i)--(v) still hold.  A priori the boundary of $\tilde N_i$ is only piecewise smooth; but, by choosing  $N(K^s_i)$ correctly at its corners we may arrange that $\partial\tilde N_i$ be smooth.
\end{demo}

{\bf Step 4.} For proving the theorem it remains to show the existence of linearizing neighbourhoods 
$N_i\subset\tilde N_i,~i=0,\dots,n$, for which the required foliations are the restriction 
to $N_i$ of the foliations 
${F}^u_i,~{F}^s_i$.  

For each orbit of $f$ in $\Si_i$, choose one $p$.  Let $\tilde N_p$ be a connected component of  $\tilde N_i$ containing $p$. There is a homeomorphism $\varphi^u_p:W^u_p\to W^u_{O}$ (resp. $\varphi^s_p:W^s_p\to W^s_{O}$) conjugating the diffeomorphisms $f^{per~p}|_{W^u_p}$ and $a_{q_p,\nu_p}|_{W^u_{O}}$ (resp. $f^{per~p}|_{W^s_p}$ and $a_{q_p,\nu_p}|_{W^s_{O}}$). In addition, for any point $z\in\tilde N_p$ there is unique pair of points $z_s\in W^s_p,~z_u\in W^u_p$ such that  $z=F^s_{i,z_u}\cap F^u_{i,z_s}$. We define a topological embedding $\tilde{\mu}_{p}:\tilde N_{p}\to\mathbb R^3$ by the formula $\tilde{\mu}_{p}(z)=(x_1,x_2,x_3)$ where  
$(x_1,x_2)={\varphi}^u_{p}(z_u)$ and $x_3={\varphi}^s_{p}(z_s)$. Choose $t_0\in(0,1]$ such that $\mathcal N^{t_0}_{q_p}\subset\tilde\mu_p(\tilde{N}_{p})$. Observe that $a_{q_p,\nu_p}\vert_{\mathcal{N}^{t_0}_{q_p}}$ is conjugate to   $a_{q_p,\nu_p}\vert_{\mathcal{N}_{q_p}}$ by $$h(x_1,x_2,x_3)= (\frac{x_1}{\sqrt{t_0}},\frac{x_2}{\sqrt{t_0}},\frac{x_3}{\sqrt{t_0}}).$$ 

Set $N_{p}=\tilde\mu^{-1}_p(\mathcal{N}^{t_0}_{q_p})$ and 
$\mu_{p}= h \tilde\mu_p:N_{p}\to \mathcal{N}_{q_p}$. Then, $N_p$ is the wanted neighbourhood with its
linearizing homeomorphism $\mu_{p}$. Set $ N_{f^k(p)}=f^k(N_p)$ and denote by $\tilde\mu_i$ a map composed from $\tilde\mu_p$ for all $p\in\Sigma_i$ such that $\tilde\mu_{f^k(p)}(x)=\tilde\mu_{p}(f^{-k}(x))$ for $x\in N_{f^k(p)}$ and $k\in\{1,\dots,per~p\}$. 

For saddle points with index Morse $1$ it is possible to prove lemmas similar to \ref{la}, \ref{lal} and, hence, construct compatible neighbourhoods.
\end{demo}

\section{The proof of the classification Theorem \ref{t.invariant}}\label{III}

Let us prove that diffeomorphisms $f,f'\in MS(M^3)$ are topologically conjugated 
if and only if their schemes are equivalent.

{\bf Necessity.} Let diffeomorphisms $f,
f'\in{MS}(M^3)$ are topologically conjugated by a homeomorphism $h:M^3\to M^3$. Set  
$\varphi=h|_{V_f}$. Then the homeomorphism  $\varphi:V_f\to V_{f'}$ conjugates the 
diffeomorphisms $f|_{V_f}$ and $f'|_{V_{f'}}$. As the  natural projects  
$p_{_f}$, $p_{_{f'}}$ are covers and $\varphi$ sends the invariant manifolds 
of the periodic points of $f$ to the invariant manifolds of the periodic points 
of $f'$ preserving dimension and stability then, due to Statement  \ref{conj}, a map  
$\hat\varphi=p_{_{f'}}\varphi p_f^{-1}:\hat V_f\to\hat V_{f'}$ is the required homeomorphism 
doing the schemes  $S_f$, $S_{f'}$ equivalent.

{\bf Sufficiency.} For proving the sufficiency of the conditions in Theorem \ref{t.invariant}, let us consider  a homeomorphism $\hat\varphi:\hat V_{f}\to\hat V_{f'}$ such that:

(1) $\eta_{_f}=\eta_{_{f'}}\hat\varphi_*$;

(2)  $\hat\varphi(\hat{\Gamma}^s_{f})=\hat{\Gamma}^s_{f'}$ and  $\hat\varphi(\hat{\Gamma}^u_{f})=\hat{\Gamma}^u_{f'}$.

From now on, the dynamical objects attached to $f'$ will be denoted by {$L'^u, L'^s, \Si'_i, \ldots$ with 
the same meaning as  $L^u, L^s, \Si_i, \ldots$} have with respect to  $f$. Due to property (1), $\hat\vp$ lifts
 to an {\it equivariant}\,\footnote{For brevity, equivariance stands for $(f,f')$-equivariance.}
homeomorphism $\varphi: V_{f}\to V_{f'}$,  that is: $f'|_{V_{f'}}= \vp f \vp^{-1}|_{V_{f'}}$\,. Due to property (2), $\varphi$ maps  $ \Gamma_f^u$ to $\Gamma_{f'}^u$ and  $ \Gamma_f^s$ to $\Gamma_{f'}^s$.
Thanks to Theorem \ref{use2-1/2} we may use compatible linearizable neighbourhoods of the saddle points of $f$ (resp. $f'$). 

An idea of the proof is the following: we  modify the homeomorphism $\varphi$ in a neighborhood of $\Ga_f^u$ such that {the} final homeomorphism preserves {the} compatible foliations, then we 
do  similar modification near $\Ga^s_f$. So we get a homeomorphism $h:M\setminus(\Om_0\cup\Om_3)\to M\setminus(\Om'_0\cup\Om'_3)$ conjugating $f|_{M\setminus(\Om_0\cup\Om_3)}$ with $f'|_{M\setminus(\Om'_0\cup\Om'_3)}$. Notice that $M\setminus(W^s_{\Om_1}\cup W^s_{\Om_2}\cup\Om_3)=W^s_{\Om_0}$ and $M\setminus(W^s_{\Om'_1}\cup W^s_{\Om'_2}\cup\Om'_3)=W^s_{\Om'_0}$. Since $h(W^s_{\Om_1})=W^s_{\Om'_1}$ and 
$h(W^s_{\Om_2})=W^s_{\Om'_2}$ then $h(W^s_{\Om_0}\setminus\Om_0)=W^s_{\Om'_0}\setminus\Om'_0$. Thus for each connected component $A$ of $W^s_{\Om_0}\setminus\Om_0$ there is a sink $\omega\in\Om_0$ such that $A=W^s_\omega\setminus\omega$. Similarly $h(A)$ is a connected component of $W^s_{\Om'_0}\setminus\Om'_0$ such that  $h(A)=W^s_{\omega'}\setminus\omega'$ for a sink $\omega'\in\Om'_0$. Then we can continuously extend $h$ to $\Om_0$ assuming $h(\omega)=\omega'$ for every $\omega\in\Om_0$. A similar extension of $h$ to $\Om_3$ finishes the proof. 

Indeed all constructions we reduce to $\mathbb R^3$ using linearizing maps. So in subsection below  we prove the crucial Proposition for linear model and in subsection \ref{red} in a sequence of lemmas we explain only how to reduce a non linear situation to the linear. 

\subsection{Linear model}

Let us recall that we denoted by $a:\mathbb R^3\to\mathbb R^3$ {the canonical  linear diffeomorphism}
with the unique fixed point {$O=(0,0,0)$} which is a saddle point with the plane $Ox_1x_2$ as the unstable manifold and the axis $Ox_3$ as the stable manifold; for simplicity, we assume  that 
 $a$ has a sign $\nu=+$. Let $$N=\{(x_1,x_2,x_3)\in\mathbb R^3:0\leq(x_1^2+x_2^2)x_3\leq 1\}.$$  Let $\rho>0,\,\delta\in(0,\frac{\rho}{4})$ and $$d=\{(x_1,x_2,x_3)\in\mathbb R^3:x_1^2+x_2^2\leq\rho^2,x_3=0\},$$ $$U=\{(x_1,x_2,x_3)\in\mathbb R^3:(\rho-\delta)^2\leq x_1^2+x_2^2\leq\rho^2,x_3=0\},$$  $$c=\{(x_1,x_2,x_3)\in\mathbb R^3:x_1^2+x_2^2=\rho^2,x_3=0\},$$ $$c^{0}=\{(x_1,x_2,x_3)\in\mathbb R^3:x_1^2+x_2^2=(\rho-\frac{\delta}{2})^2,x_3=0\},$$ $$c^{1}=\{(x_1,x_2,x_3)\in\mathbb R^3:x_1^2+x_2^2=(\rho-\delta)^2,x_3=0\}.$$
Let $K=d\setminus int\,a^{-1}(d)$, $V=(K\cup a(K))\cap\{(x_1,x_2,x_3)\in\mathbb R^3:x_1\geq 0,x_3=0\}$ and $\beta=U\cap Ox^+_1$, where $Ox_1^+=\{(x_1,x_2,x_3)\in\mathbb R^3:x_3^2+x_2^2=0,x_1>0\}$. Let $T\subset Ox_1x_2$ be an $a$-invariant closed one-dimensional lamination such that either $c=U\cap T$ or $T$ (maybe empty) is transversal to $\partial U$ and every connected component of $T\cap U$ is segment which has unique intersection point with each connected component of $\partial U$. 

\begin{figure}[h]\centerline{\includegraphics[width=12cm,height=14cm]{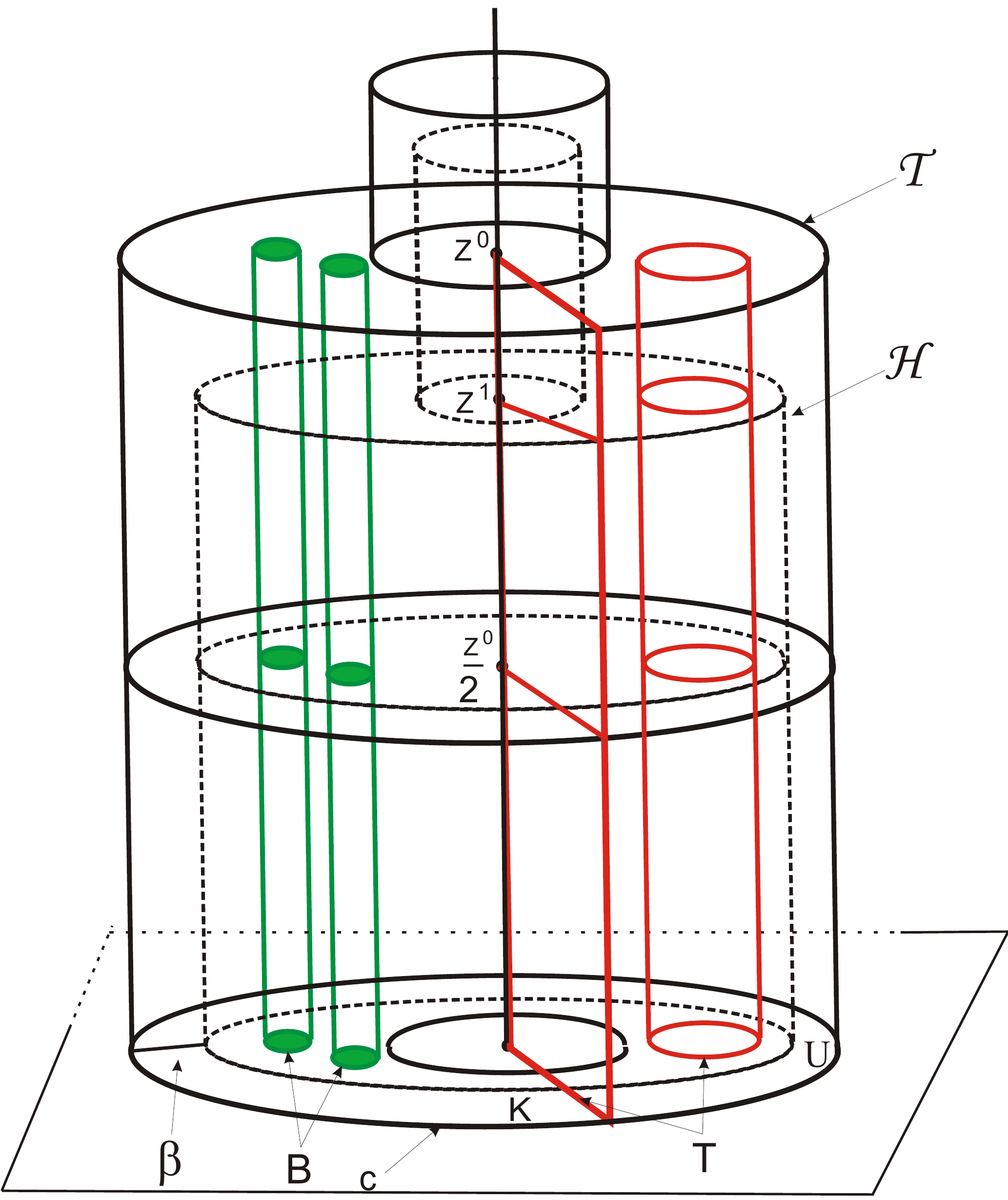}}
\caption{\small A linear model} \label{fin}
\end{figure}

Choose a point $Z^0=(0,0,z^0)\in Ox_3^+$ such that $\rho^2\cdot z^0<\frac{1}{4}$. Then,
choose a point $Z^1=(0,0,z^1)$ in $Ox_3^+$ so that $z^0>z^1>\frac{z^0}{2}$. Let  $\Pi(z)=\{(x_1,x_2,x_3)\in\mathbb R^3:x_3=z\}$. For every set $A\subset Ox_1x_2$, let $\tilde A=A\times[0,z^0]$. Denote by {$\mathcal W$} the 3-ball bounded by the surface $\tilde c$ and the two planes $\Pi(z^0)$ and $\Pi(\frac{z^{0}}{2})$. Let $\Delta$ be a closed 3-ball bounded by the surface $\tilde c^{1}$ and the two planes $Ox_1x_2$ and $\Pi(z^1)$. Let 
$$\mathcal T=\bigcup\limits_{k\in\mathbb Z}a^{k}(\tilde d){\quad{\rm and}\quad}
\mathcal H=\bigcup\limits_{k\in\mathbb Z}a^{k}(\Delta).$$
Notice that {the construction yields} $\mathcal H\subset \mathcal T$ and makes {$\mathcal W$}
a fundamental domain for the action of $a$ on $\mathcal T$ (see Figure \ref{fin}). 

\begin{prop} \label{F1} Let $z^0>z^1>\dots>z^m>\frac{z^0}{4}>0$ and $\xi:\mathcal T\setminus Ox_3\to N$ be a topological embedding with the following properties:

{\rm (i)} $\xi a=a\xi$;

{\rm (ii)} $\xi$ is the identity on $Ox_1x_2$;

{\rm (iii)} $\xi(\Pi(z^0\cap \mathcal T)) =\Pi(z^0)$ and $\xi(\Pi(z^i)\cap\partial\mathcal T)\subset \Pi(z^i),\,i\in\{2,\dots,m\}$;

{\rm (iv)} $\xi(\tilde c)\cap\tilde c^0=\emptyset$, $\xi(\tilde c^1)\cap\tilde c^0=\emptyset$ and $\xi(\tilde\beta)\subset\tilde V$;

{\rm (v)} $\xi(\tilde T\cap\mathcal T)\subset\tilde T$ and $\xi^{-1}(\tilde T\cap\mathcal T)\subset\tilde T$. 

\nd Then there is a topological embedding $\zeta:\mathcal T\to N$ such that  

{\rm 1)} $\zeta a=a\zeta$;

{\rm 2)} $\zeta$ is the identity on $\mathcal H$;

{\rm 3)} $\zeta(\Pi(z^i)\cap\mathcal T)\subset \Pi(z^i),\,i\in\{0,2,\dots,m\}$

{\rm 4)} $\zeta$ is $\xi$ on $\partial{\mathcal T}$;

{\rm 5)} $\zeta(\tilde T\cap\mathcal T)\subset\tilde T$ and $\zeta^{-1}(\tilde T\cap\mathcal T)\subset\tilde T$. 

Moreover, if $\xi$ is identity on $\tilde B$ for a set 
$B\subset (K\setminus U)$ then $\zeta$ is also identity on $\tilde B$.
\end{prop} 
\begin{demo} By the construction $\mathcal W$ is a fundamental domain of the diffeomorphism $a$ restricted to $\mathcal T\setminus Ox_1x_2$. Therefore to proof the proposition it is enough to construct the topological embedding  $\zeta_{\mathcal W}=\zeta|_{\mathcal W}$ with the properties 1)-5) on the set  $\mathcal{W}$. Because we can extend $\zeta_{W}$ to $\mathcal T$ by the formula $\zeta(x)=a^{-k}(\zeta_{W}(a^{k}(x)))$, where  $a^{k}(x)\in W$ on $\mathcal T\setminus Ox_1x_2$ and define $\zeta$ to be identity on $Ox_1x_2$.   

Firstly we define $\zeta_{\mathcal W}$ to be identity on $\mathcal W\cap\mathcal H$. Divide the remain part $\mathcal Q=\mathcal W\setminus\mathcal H$ on the following pieces: $Q_1=\tilde U\cap\mathcal Q$, $Q_2=a^{-1}(\tilde U)\cap \mathcal Q$, $Q_3=\mathcal Q\setminus(Q_1\cup Q_2)$ (see Figure \ref{fin}).  Define a topological embeddings $h_{Q_3}:Q_3\to\mathcal Q$ by the following way. Let $\kappa:[z^{1},z^0]\to[0,z^0]$ be a homeomorphism given by the formula $\kappa(z)=\frac{z^0(z-z^1)}{z^0-z^1}$ and  $h_{Q_3}=\kappa^{-1}\xi\kappa$ on $Q_3=(K\setminus int~U)\times[z^1,z^0]$. Thus we get the desired  embedding $\zeta_{\mathcal W}$ on $\mathcal T\cup Q_3$.  

By property 4) of the map $\xi$ we have $\xi(\tilde c^1)\cap \tilde c=\emptyset$ $(\xi(\tilde c)\cap \tilde c^1=\emptyset)$ then the surfaces $\tilde c$, $\zeta_{\mathcal W}(\tilde c^1)$, $\Pi(z^0)$, $\Pi(z^{m+1})$  $(\zeta_{\mathcal W}(a^{-1}(\tilde c)\cap\mathcal W), a^{-1}(\tilde c^1)\cap\mathcal W$, $\Pi(z^0)$, $\Pi(z^{1}))$ bound a closed 3-dimensional set, denote it by $\check Q_1,\,(\check Q_2)$. Due to Proposition \ref{solid} the sets $\check Q_1,\,\check Q_2$ are solid tori. Let $$A^i=Q_1\cap \Pi(z^i),\quad\check A^i=\check Q_1\cap \Pi(z^i),\quad i\in\{0,2,\dots,m+1\}.$$ Let $S^i\,(\check S^i),\,i\in\{2,\dots,m+1\}$ be the closure of a connected component of the set $Q_1\setminus\bigcup\limits_{j=2}^m\Pi(z^i)\,(\check Q_1\setminus\bigcup\limits_{j=2}^m\Pi(z^i))$ bounded by $\Pi(z^i)$ below. By Proposition \ref{solid}  $S^i\,(\check S^i)$ is a solid torus. Further let us consider two possibilities distinctly: case 1) $c=U\cap T$ or $T$ is empty; case 2) $T$ is transversal to $\partial U$ and every connected component of $T\cap U$ is segment which has unique intersection point with each connected component of $\partial U$.  

{\bf Case 1.} Let us define a topological embedding $h_{A^i}:A^i\to \check A^i$ such that:

- $h_{A^i}$ is $h_{\mathcal W}$ on $\tilde c^1\cap A^i$, is $\xi$ on $\tilde c\cap A^i$, $h_{A^i}(\tilde\beta\cap A^i)\subset \tilde V$ and $h_{A^i}(\tilde T\cap A^i)=\tilde T\cap\check A^i$ for $i\in\{2,\dots,m+1\}$;

- $h_{A^0}=\xi|_{A^0}$.

 Moreover we have a homeomorphism $h_{\partial S^i}:\partial S^i\to\partial{\check S}^i$ which coincides with $h_{A^j}$ on $S^i\cap A^j$, with $\xi$ on $S^i\cap\tilde c$ and with $\zeta_{\mathcal W}$ on $S^i\cap\tilde c^1$. By the construction a curve $\mu^i=\partial (S^i\cap \tilde \beta)$ is a meridian of $S^i$ and $\check\mu^i=h_{\partial S^i}(\mu^i)$ is a meridian of $\check S^i$. Then (see \cite{Ro} for example) there is a homeomorphism $h_{S^i}:S^i\to\check S^i$ such that $h_{S^i}|_{\partial S^i}=h_{\partial S^i}$.

Similarly we have a homeomorphism $h_{\partial Q_2}:\partial Q_2\to\partial{\check Q_2}$ such that  $h_{\partial Q_2}|_{Q_2\cap\Pi(z^1)}=id|_{Q_2\cap\Pi(z^1)}$,  and $h_{Q_2\cap\Pi(z^0)}=a^{-1}h_{A_1^{m+1}}a|_{Q_2\cap\Pi(z^0)}$, $h_{\partial Q_2}|_{Q_2\cap a^{-1}(\tilde c)}=\zeta_{\mathcal W}|_{Q_2\cap a^{-1}(\tilde c)}$ and $h_{\partial Q_2}|_{Q_2\cap a^{-1}(\tilde c^1)}=id|_{Q_2\cap a^{-1}(\tilde c^1)}$. Hence there is a homeomorphism $h_{Q_2}:Q_2\to\check Q_2$ such that $h_{Q_2}|_{\partial Q_2}=h_{\partial Q_2}$.
  
Thus the required homeomorphism is defined by the formula
$$\zeta|_{W}(x)=\cases{x, x\in\mathcal T;\cr
h_{Q_3}(x), x\in Q_3;\cr
h_{S^i}(x), x\in S^i, i\in\{2,\dots,m+1\};\cr
h_{Q_2}(x),x\in Q_2.\cr}.$$

{\bf Case 2.} According to Proposition \ref{B2} and Remarks \ref{B1} there is a topological embedding $h_{A^i}:A^i\to \check A^i$ such that:

- $h_{A^i}$ is $h_{\mathcal W}$ on $\tilde c^1\cap A^i$, is $\xi$ on $\tilde c\cap A^i$ and $h_{A^i}(\tilde T\cap A^i)=\tilde T\cap\check A^i$ for $i\in\{2,\dots,m+1\}$. 

Let $h_{A^0}=\xi|_{A^0}$. Thus we have a homeomorphism $h_{\partial S^i}:\partial S^i\to\partial{\check S}^i$ which coincides with $h_{A^j}$ on $S^i\cap A^j$, with $\xi$ on $S^i\cap\tilde c$ and with $\zeta_{\mathcal W}$ on $S^i\cap\tilde c^1$, also $h_{\partial S^i}(\partial S^i\cap\tilde T)=\partial{\check S}^i\cap\tilde T$. According to Proposition \ref{B2} there is a homeomorphism $h_{S^i}:S^i\to S^i$ such that $h_{S^i}$ is $h_{\partial S^i}$ on $\partial S^i$ and $h_{S^i}(\tilde T\cap S^i)=\tilde T\cap\check S^i$.  

Similarly we can construct a homeomorphism $h_{\partial Q_2}:\partial Q_2\to\partial{\check Q_2}$ such that  $h_{\partial Q_2}|_{Q_2\cap\Pi(z^1)}=id|_{Q_2\cap\Pi(z^1)}$,  and $h_{Q_2\cap\Pi(z^0)}=a^{-1}h_{A_1^{m+1}}a|_{Q_2\cap\Pi(z^0)}$, $h_{\partial Q_2}|_{Q_2\cap a^{-1}(\tilde c)}=\zeta_{\mathcal W}|_{Q_2\cap a^{-1}(\tilde c)}$ and $h_{\partial Q_2}|_{Q_2\cap a^{-1}(\tilde c^1)}=id|_{Q_2\cap a^{-1}(\tilde c^1)}$. Hence there is a homeomorphism $h_{Q_2}:Q_2\to\check Q_2$ such that $h_{Q_2}|_{\partial Q_2}=h_{\partial Q_2}$ and $h_{Q_2}(Q_2\cap \tilde T)=\check Q_2\cap \tilde T$.
  
Thus the required homeomorphism is defined by the formula
$$\zeta|_{W}(x)=\cases{x, x\in\mathcal T;\cr
h_{Q_3}(x), x\in Q_3;\cr
h_{S^i}(x), x\in S^i, i\in\{2,\dots,m+1\};\cr
h_{Q_2}(x),x\in Q_2.\cr}.$$
\end{demo}

\subsection{Reduction to the linear model}\label{red}

Recall the partition $\Si_0\sqcup\cdots\sqcup \Si_n$ associated with the Smale order on the periodic points {of index 2}. Firstly we formulated needed results from \cite{BoGrLaPo}.  

\begin{lemm}[\cite{BoGrLaPo}, Lemma 4.1]\label{leaf_number} {There is a unique continuous extension $\vp:\Si_i\to\Si'_i$ for every 
$i= 0, \dots, n$. This extension is equivariant and bijective from $\Om_2$ to $\Om'_2$ which preserves the type of the orientation and the period of points.}
\end{lemm}

Let us introduce the {\it radial functions} $r^u_i, r^s_i: N_i\to [0,+\infty)$ defined by: 
$$r^u_i(x)=\Vert\mu_i(x^u_i)\Vert^2\quad{\rm and} \quad r^s_i(x)=  (\mu_i(x^s_i))^2.$$
With this definition at hand, the neighbourhood $N_i^t$ of $\Si_i$ is defined by the inequation
$$r^u_i(x)\,.\,r^s_i(x)<t\,.$$ Observe that the radial function $r^s_i$ endows each stable separatrix $\gamma_p$ of $p\in \Si_p$ with a natural order (and similarly  with $^\prime$). 

\begin{lemm}[\cite{BoGrLaPo}, Lemma 4.2] \label{pusk} There is a unique continuous  extension  {of $\vp|_{\Ga^u_f}$}
$$\vp^{us}:L^u\to L'^u$$ such that {the following holds}: 

\nd {\rm 1)} If  $x\in W^u_j\cap W^s_i$, $j>i$, then $\vp^{us}(x)\in  W'^u_j\cap W'^s_i$.
  
\nd {\rm 2)} If $x$ and $y$ lie in $\gamma_p\cap L^u$ with $r^s_p(x)<r^s_p(y)$,  then $\vp^{us}(x)$ and 
$\vp^{us}(y)$ lie in $ \gamma'_{\varphi(p)}\cap L'^u$ with
 $r'_{\varphi(p)}(\vp^{us}(x))<r'_{\varphi(p)}(\vp^{us}(y))$.
\end{lemm}

\begin{rema} \label{mai} Due to Lemma \ref{Ge} we may assume  that in all lemmas below we choose values  $t=\beta_i,a_i,...$  such that the boundary of the linearizing neighbourhood  $N^t_i$ does not contain any heteroclinic point. 
\end{rema} 
 
\begin{lemm}[\cite{BoGrLaPo}, Lemma 4.4] \label{psibe} There are numbers $\beta_0,\dots,\beta_n\in (0,1]$
such that, for  every {$i\in\{0,\dots,n\}$, for every point $p\in\Si_i$ and $x\in N^{\beta_i}_p\cap L^u$, 
 the next inequality holds:
  $$r'^u_{i}(\vp^{us}(x^u_i))\,.\, r'^s_{i}(\vp^{us}(x^s_i))<\frac12\,.$$} 
\end{lemm}

Let us set $a_0=\beta_0$ and $a_1=\beta_1$.

\begin{lemm} \label{psu} If $n\geq 2$ then there are numbers $a_2\in (0,\beta_2],\dots,a_n\in (0,\beta_n]$ with following property for each $i\in\{2,\dots,n\}$: for $0\leq k\leq i-2$ the intersection $$W^s_{k}\cap (N^{a_i}_{i}\smallsetminus(\bigcup\limits_{\mu=k+1}^{i-1} N^{a_\mu}_{\mu}))$$ is either empty or consists of open arcs each of which is a leaf of the foliation  $F^s_{i}\cap{N}^{a_{i}}_{i}$.
\end{lemm}
\begin{demo} We will construct the sequence by induction on $i=2,\dots,n$. 

For $i=2$ the unstable manifolds of points from $\Si_{2}$ have only heteroclinic intersections with the stable manifolds of  saddles from $\Si_{k}$ with $k<2$. Thus, the projection $\hat W^{u}_{2,1}$ is a union of finite number pairwise disjoint smooth torus and Klein bottles in the manifold $\hat V_{1}$. The set $\hat{ L}^{s}_{1}$ is a compact one-dimensional lamination, the intersection  $\hat W^{u}_{2,1}\cap\hat{ L}^{s}_{1}$ is transversal and consists of at most countable set of points which are the projection with respect to $p_{1}$ of the heteroclinic orbits from the unstable manifolds $W^{u}_{2}$. The set $\hat W^{u}_{2,1}\cap(\hat{ W}^{s}_{0,1}\smallsetminus\hat N^{a_1}_{1,1})$ is either empty or  consists of finite number $q_0\in\mathbb N$ points. Due to remark \ref{mai} there is a number $a_2\in(0,\beta_2]$ such that the intersection $\hat N^{a_2}_{2,1}\cap(\hat{W}^{s}_{0,1}\smallsetminus\hat N^{a_1}_{1,1})$ is either empty or consists of exactly $q_0$ intervals each of which is a leaf of the foliation  $\hat F^s_{2}\cap\hat{N}^{a_{2}}_{2}$ (see Figure \ref{arc} where $q_0=1$). Hence, the intersection $W^s_{0}\cap (N^{a_2}_{2}\smallsetminus N^{a_1}_{1}))$ is either empty or consists of open arcs each of which is a leaf of the foliation  $F^s_{2}\cap{N}^{a_{2}}_{2}$.

\begin{figure}[ht]
\centerline{\includegraphics[width=13cm,height=9cm]{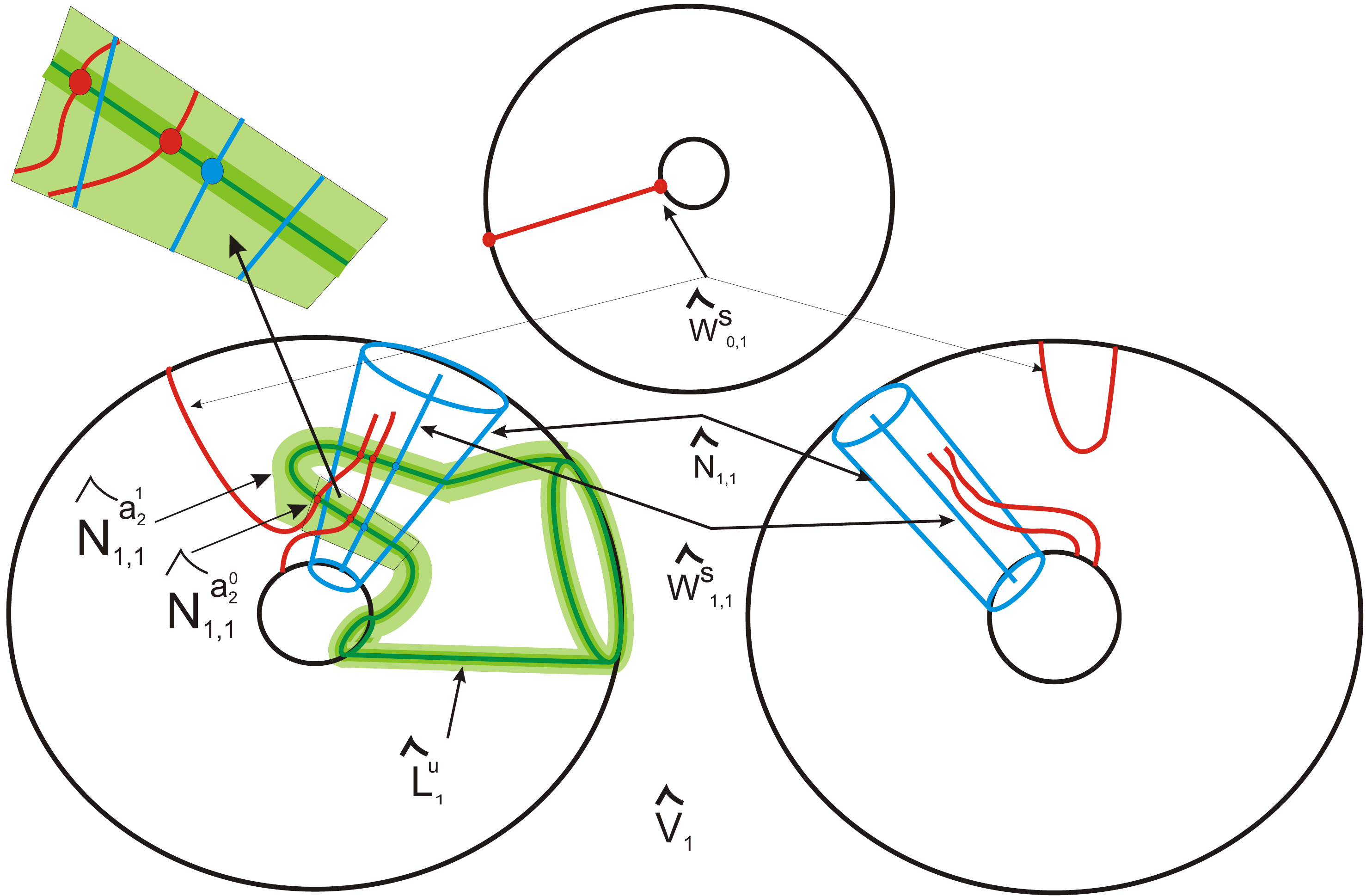}}\caption{\small Illustration {of} the proof of 
Lemma \ref{psu} for $i=2$.} \label{arc}
\end{figure}

Let us describe a finding of the number  $a_{i},~i>0$ supposing that the numbers  $a_{j}$ with desired  properties are already constructed for all $j={0},\dots,i-1$. 

The unstable manifold of points from $\Si_{i}$ have only heteroclinic intersections with the stable manifolds of  saddles from $\Si_{k}$ with $k<i$. Thus, the projection $\hat W^{u}_{i,i-1}$ is a smooth torus in the manifold $\hat V_{i-1}$. The set $\hat{ L}^{s}_{i-1}$ is a compact one-dimensional lamination, the intersection  $\hat W^{u}_{i,i-1}\cap\hat{ L}^{s}_{i-1}$ is transversal and consists of at most countable set of points which are the projection with respect to $p_{i-1}$ of the heteroclinic orbits from the unstable manifolds $W^{u}_{i}$.  For each $k=0,\dots,i-2$ the intersection $\hat W^{u}_{i,i-1}\cap(\hat{ W}^{s}_{k,i-1}\smallsetminus\bigcup\limits_{\mu=k+1}^{i-1}\hat N^{a_\mu}_{\mu,i-1})$ is either empty or consists of $q_k\geq 0$ points. Due to remark \ref{mai} there is a number $a_i^{k}\in(0,\beta_i]$ such that  the intersection $\hat N^{a^k_i}_{i,i-1}\cap(\hat{ W}^{s}_{k,i-1}\smallsetminus\bigcup\limits_{\mu=0}^{k-1}\hat N^{a_\mu}_{\mu,i-1})$ is either empty or
  consists of exactly $q_k$ intervals each of which is a leaf of the foliation  $\hat F^s_{i}\cap\hat{N}^{a_i^{k}}_{i}$.

Thus $a_i=min\{a^0_i,\dots,a_i^{i-2}\}$ is the required.
\end{demo}

The corollary  below immediately follows  from lemma \ref{psu}.

\begin{coro} \label{con} For each $k\in\{0,\dots,n-1\}$ the intersection $\hat W^{s}_{k,k}\cap (\bigcup\limits_{i=k+1}^n\hat{N}^{a_{i}}_{i,k})$ consists of finite number of open arcs $\hat I^{k}_1,\dots,\hat I^{k}_{r_{k}}$ such that $\hat I^{k}_l$ for each $l=1,\dots,r_{k}$ is a connected component of intersection $\hat W^{s}_{k,k}\cap\hat{N}^{a_{i}}_{i,k}$ for some $i>k$.
\end{coro}

For brevity, for $i=0,\ldots, n$, we denote by $\vp^{u}_i$ the restriction $\vp^{us}|_{W^u_i}$ 
in the rest of the proof of Theorem \ref{t.invariant}. Let {$\psi_i^s: W^s_i\to W'^s_i$} be any equivariant homeomorphism which extends {$\vp^{us}\vert_{W_i^s\cap L^u}$} and let   $t_i\in (0,1)$ be a small enough number so that, for every $x\in N_i^{t_i}$, the next inequality holds:
$$(*)_i\quad\quad r'^s(\vp_i^u(x^u_i))\, . \,r'^u(\psi^s_i(x^s_i))<1.$$
In this setting, one derives an equivariant embedding $\phi_{\vp^u_i ,\psi^s_i}: N_i^{t_i} \to N'_i$ which is defined by sending $x\in N_i^{t_i}$ to $\left(\vp_i^u(x^u_i),\psi^s_i(x^s_i)\right)$.

\begin{lemm} \label{1dim} There is a homeomorphism $\psi^s:L^s\to L^{\prime s}$ consisting of conjugating  homeomorphisms $\psi^s_{0}:W^s_{0}\to W^{\prime s}_{0},\dots,\psi^s_{n}:W^s_{n}\to W^{\prime s}_{n}$ such that for each $i\in\{0,\dots,n\}$:

1) $\psi^s_{i}|_{W^s_{i}\cap L^u}=\vp^{us}|_{W^s_i\cap L^u}$;

2) the topological embedding ${\phi}_{{\vp}^{u}_i,\psi^s_i}$ is well-defined on $N^{a_i}_i$;

3) if $x\in (W^s_i\cap N_j^{a_j}),j>i$ then $\psi^s_i(x)={\phi}_{{\vp}^{u}_j,\psi^s_j}(x) $.
\end{lemm}
\begin{demo} We are going to construct $\psi^s_i$, by a decreasing induction on $i$ from $i=n$ up to $i=0$. 

The stable manifolds of the saddles from $\Si_n$ does not have any heteroclinic intersection. The projection $\hat W^{s}_{n,n}$ is a smooth submanifold of the manifold $\hat V_{n}$ consisting of finite number connected component homeomorphic to circle. The same holds for $\hat W^{\prime s}_{n,n}$. Denote by  $\hat\psi^s_{n}:\hat W^{s}_{n,n}\to\hat W^{\prime s}_{n,n}$ a homeomorphism with following property: let $\hat\gamma_n$ be a connected component of $\hat W^{s}_{n,n}$ and $U(\hat \gamma_n)$ be its neighbourhood, $\hat\gamma'_n=\hat\psi^s_{n}(\hat\gamma_{n})$ and $U(\hat\gamma'_n)$ be  a  neighbourhood of $\hat\gamma'_n$ then $\hat\varphi_n(U(\hat \gamma_n))\cap U(\hat \gamma'_n)\neq\emptyset$. Let $\tilde\psi^{s}_{n}:W^s_{n}\smallsetminus\Si_{n}\to {W}^{\prime s}_{n}\smallsetminus\Si'_{n}$ be a lift with respect to $p_{n}$ of $\hat\psi^s_{n}$ on $W^s_{n}\smallsetminus\Si_{n}$.  

Let $p\in\Si_n$ and $\gamma_p$ is a connected component of $W^s_p\smallsetminus p$. Let us choose a fundamental domain $I_{\gamma_p}$ of $f^{per~\gamma_p}|_{\gamma_p}$. Set $N_{I_{\gamma_p}}=\{x\in N^{a_n}_p~|~x^s_n\in I_{\gamma_p}\}$  and $\la'^u_{\gamma_p}=sup\{r'^u_p(\vp^{us}(x^u_n))~|~x\in N_{I_{\gamma_p}}\}$. For $k\in\mathbb Z$ we set $\la'^s_{\gamma_p}(k)=\frac{1}{2^k}\cdot sup\{r'^s_{p'}(\tilde\psi^s_n(x))~|~x\in I_{\gamma_p}\}$. As $\la'^s_{\gamma_p}(k)$ tends to 0 as $k$ tends to $\infty$ then there is $k_*\in\mathbb N$ such that $\la'^s_{\gamma_p}(k_*)\cdot\la'^u_{\gamma_p}<1$. Set $\psi^s_n|_{\gamma_p}=f'^{k_*}\tilde \psi^s_n|_{\gamma_p}$. We define similarly $\psi^s_n$ on other connected component of $W^s_p\smallsetminus p$ (different from $\gamma_p$) and set  $\psi^{s}_{n}(p)=p'$. Thus $r'^u_{p'}(\vp^{us}(x^u_n))\cdot r'^s_{p'}(\psi^{s}_n(x^s_n))\leq \la'^s_{\gamma_p}(k_*)\cdot\la'^u_{\gamma_p}<1$ and, hence, the topological embedding ${\phi}_{{\vp}^{u}_n,\psi
 ^s_n}$ is well-defined on $N^{a_n}_p$. Then we do the same for each point $p\in\Si_i$. 

Let us describe a construction of the homeomorphism  $\psi^{s}_{i},~i<n$ supposing that the homeomorphisms  $\psi^{s}_{n},\dots,\psi^s_{i+1}$ are already constructed. 

The stable manifolds of points from $\Si_{i}$ have heteroclinic intersections only with unstable manifolds of the  saddles $\Si_{j}$ with $j>i$. Thus, the projection $\hat W^{s}_{i,i}$ is a smooth submanifold of the manifold $\hat V_{i}$ consisting of finite number connected components homeomorphic to circle. Due to corollary \ref{con}, the intersection $\hat W^{s}_{i,i}\cap (\bigcup\limits_{j=i+1}^n\hat{N}^{a_{j}}_{j})$ consists of finite number of open arcs $\hat I^{i}_1,\dots,\hat I^{i}_{r_{i}}$ such that $\hat I^{i}_l$ for each $l=1,\dots,r_{i}$ is a connected component of intersection $\hat W^{s}_{i,i}\cap\hat{N}^{a_{j}}_{j,i}$ for some $j>i$.  

Denote by $I^{i}_l$ a connected component of the set $p_{{i}}^{-1}(\hat I^{i}_l)$. The arc $I^{i}_l$ is an arc in $N^{a_{j}}_{j}$ intersecting $W^u_{j}$ at unique point $x^{i}_l$. Set $x^{\prime {i}}_l=\vp^{us}(x^{i}_l)$ and denote by $I^{\prime {i}}_l$ the connected component of $W^{\prime s}_{i}\cap N^{\prime}_{j}$ passing through the point $x^{\prime {i}}_l$. By assumption of the induction on the set $I^{i}_l$ the homeomorphism $\phi_{\vp^{u}_j,\psi^s_j}$ is well-defined, let us set $\psi_{I^{i}_l}=\phi_{\vp^{us}_j,\psi^s_j}|_{I^i_l}:I^{i}_l\to I^{\prime {i}}_l$.

\begin{figure}[h]
\centerline{\includegraphics[width=13cm,height=9cm]{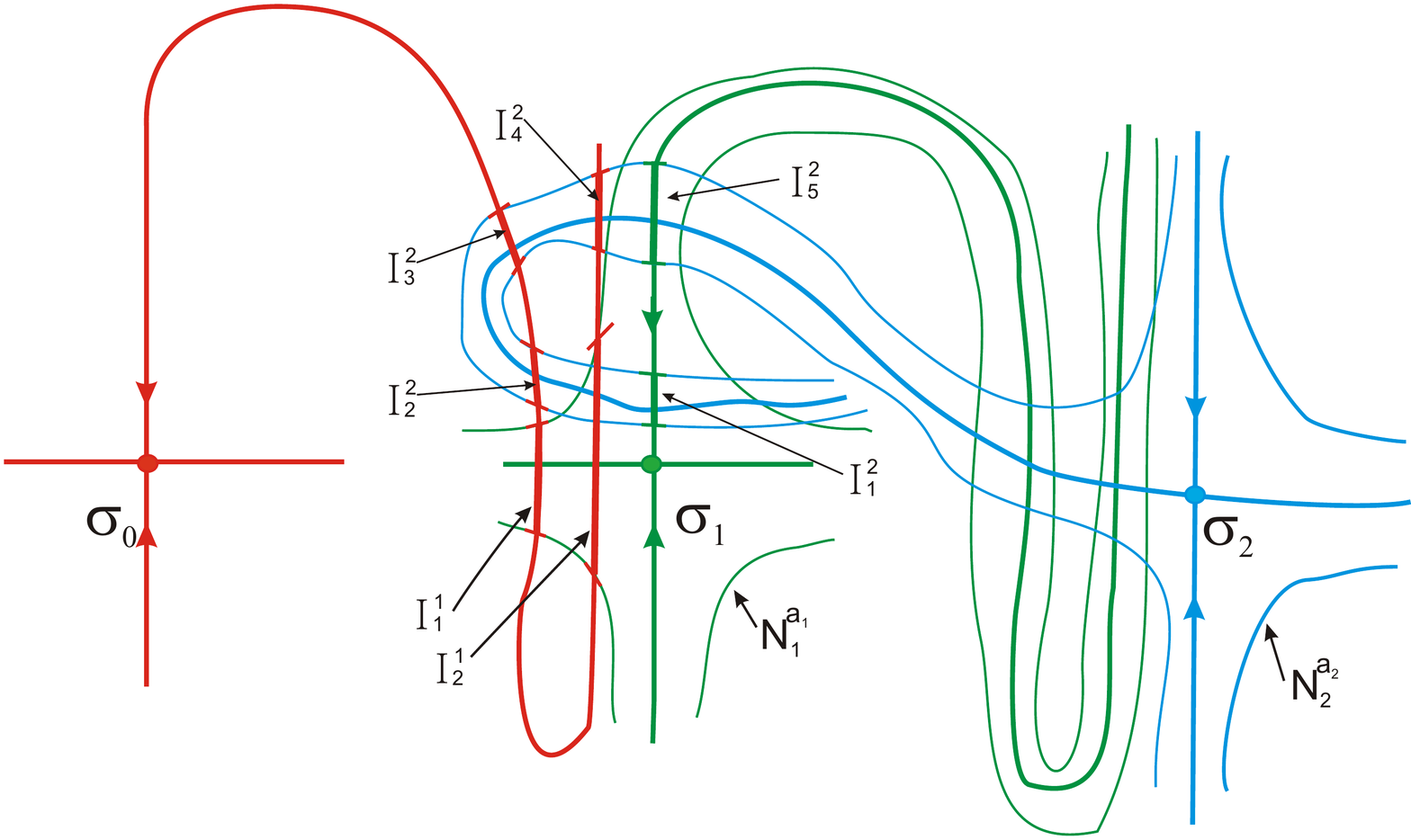}}\caption{\small Illustration {of} the proof of Lemma \ref{1dim}.} \label{segm}
\end{figure}

Set 
$$\hat\psi^s_{{i},l}=p'_{{i}} \psi_{I^{i}_l}(p_{{i}}|_{I^{i}_{l}})^{-1}\colon \hat I^{i}_l\to \hat{W}^{\prime s}_{{i},{i}}.$$ 

Notice that $\hat\psi^s_{{i},l}$  does not depend on the choice of the lift $I^{i}_l$ of $\hat I^{i}_l$. Indeed, if $\tilde I^{i}_l$ be a connected component of the set $p_{{i}}^{-1}(\hat I^{i}_l)$ different from $I^{i}_l$. Then there is unique $z\in\mathbb Z$ such that $\tilde I^{i}_l=f^z(I^{i}_l)$. As $\psi^s_{j}$ conjugates $f|_{W^s_{j}}$ and $f'|_{W^{\prime s}_{j}}$ then $\psi_{\tilde I^{i}_l}(f^z(x))={f'}^z(\psi_{I^{i}_l}(x))$ for any $x\in I^{i}_l$. It  means that $p'_{{i}}  \psi_{I^{i}_l}  (p_{{i}}|_{I^{i}_{l}})^{-1}=p'_{{i}}  \psi_{\tilde I^{i}_l}  (p_{{i}}|_{\tilde I^{i}_{l}})^{-1}$. 

By the construction the map $\hat\psi^s_{{i},l},~l=1,\dots,{r_{i}}$ coincides on $\hat I^{i}_{l}\cap\hat{ L}^{u}_{{i}}$ with $p'_{{i}} \vp^{us}(p_{{i}}|_{\hat I^{i}_{l}\cap\hat{ L}^{u}_{{i}}})^{-1}$. Due to lemma \ref{pusk}, the map $\vp^{us}$ sends ${W}^{s}_{{i}}\cap L^u$ to ${W}^{\prime s}_{{i}}\cap {L'}^u$ with preserving the order on each connected component $W^s_i\smallsetminus\Si_i$ and $W'^s_i\smallsetminus\Si'_i$. Then there is a homeomorphism $\hat\psi^s_{i}:\hat W^{s}_{{i},{i}}\to\hat W^{\prime s}_{{i},{i}}$ coinciding with $\hat\psi^s_{{i},l}$ on $\hat I^{i}_{l},~l=1,\dots,r_{{i}}$. Moreover $\hat\psi^s_{i}$ possesses  following property: let $\hat\gamma$ be a connected component of $\hat W^{s}_{i,i}$ and $U(\hat \gamma)$ be its neighbourhood, $\hat\gamma'=\hat\psi^s_{i}(\hat\gamma)$ and $U(\hat\gamma')$ be  a  neighbourhood of $\hat\gamma'$ then $\hat\varphi_i(U(\hat \gamma))\cap U(\hat \gamma')\neq\emptyset$. Denote by $\tilde\psi^{s}_{i}:W^s_{i}\to {W}^{\prime s}_
 {i}$ a homeomorphism which is a lift with respect to $p_{{i}}$ of $\hat\psi^s_{{i}}$ on $W^s_{{i}}\smallsetminus\Si_{i}$  such that it coincides with $\vp^{us}$ on $W^s_{i}\cap L^u$ and $\psi^{s}_{i}|_{\Si_{i}}=\vp|_{\Si_{i}}$. 

Let $p\in\Si_i$ and $\gamma_p$ be a connected component of $W^s_p\smallsetminus p$ containing heteroclinic points. Let $N^{a_i}_{\gamma_p}$ be a connected component of $N^{a_i}_{p}\smallsetminus W^u_{p}$ containing $\gamma_p$. Set $\gamma_{p'}=\tilde\psi^{s}_{i}(\gamma_p)$ and $\psi^{s}_i|_{\gamma_p}=\tilde\psi^{s}|_{\gamma_p}$. Let us show that $r'^s_{p'}(\psi^{u}_i(x^u_i))\cdot r'^u_{p'}(\psi^{s}_i(x^s_i))<1$ for each $ x\in N_{\gamma_p}^{a_i}$. Indeed, let us choose a heteroclinic point $y\in \gamma_p$. Set  $\la'^u_{\gamma_p}=sup\{r'^u_{p}(\vp^{us}(x^u_i))~|~x\in (N^{a_i}_{\gamma_p}\cap F^u_{i,y})\}$ and $\la'^s_{\gamma_p}=r'^s_{p'}(\vp^{us}(y))$. By lemma \ref{psibe} $\la'^u_{\gamma_p}\cdot\la'^s_{\gamma_p}<\frac12$. Denote by $Q_{\gamma_p}$ a subset of $N^{a_i}_{\gamma_p}$ which bounded by $\partial N^{a_i}_{p}, F^u_{i,y}$ and $f^{per~\gamma_p}(F^u_{i,y})$. Notice that $Q_{\gamma_p}$ is a fundamental domain of  $f^{per~{\gamma_p}}|_{N_{\gamma_p}^{a_i}}$. By the construc
 tion $r'^u_{p'}(\psi^{u}_i(x^u_i))\leq 2\la'^u_{\gamma_p}$ and $r'^s_{p'}(\psi^{s}_i(x^s_i))\leq \la'^s_{\gamma_p}$ for any $x\in Q_{\gamma_p}$ and, hence, $r'^s_{p'}(\psi^{u}_i(x^u_i))\cdot r'^u_{p'}(\psi^{s}_i(x^s_i))\leq 2\la'^u_{\gamma_p}\cdot\la'^s_{\gamma_p}<1$ for each $x\in Q_{\gamma_p}$. Thus ${\phi}_{{\vp}^{u}_i,\psi^s_i}(N_{p}^{a_i})\subset N'_{p'}$.

Let $p\in\Si_i$ and $\ell_p$ be a connected component of $W^s_p\smallsetminus p$ which does not contain heteroclinic points. Let $N^{a_i}_{\ell_p}$ be a connected component of $N^{a_i}_{p}\smallsetminus W^u_p$ containing $\ell_p$. Set $\ell'_p=\tilde\psi^{s}_{i}(\ell_p)$. Let us choose a fundamental domain $I_{\ell_p}$ of $f^{per~\ell_p}|_{\ell_p}$. Set $N_{I_{\ell_p}}=\{x\in N^{a_i}_{p}~|~x^s\in I_{\ell_p}\}$  and $\la'^u_{\ell_p}=sup\{r'^u_{p'_-}(\psi^{u}_i(x^u_i))~|~x\in N(I_{\ell_p})\}$. For $k\in\mathbb Z$ we set $\la'^s_{\ell_p}(k)=\frac{1}{2^k}\cdot sup\{r'^s_{p'}(\tilde\psi^s_i(x))~|~x\in I_{\ell_p}\}$. As $\la'^s_{\ell_p}(k)$ tends to 0 as $k$ tends to $\infty$ then there is $k_*\in\mathbb N$ such that $\la'^s_{\ell_p}(k_*)\cdot\la'^u_{\ell_p}<1$. Set $\psi^s_i|_{\ell_p}=f'^{k_*}\tilde \psi^s_i|_{\ell_p}$. Thus ${\phi}_{{\vp}^{u}_i,\psi^s_i}(N_{\ell_p}^{a_i})\subset N'_{p}$.

Finally, a map $\psi^s:L^s\to L^{\prime s}$ consisting of the homeomorphisms $\psi^s_{0}:W^s_{0}\to W^{\prime s}_{0},\dots,\psi^s_{n}:W^s_{n}\to W^{\prime s}_{n}$ is a homeomorphism due to following property: if $x\in (W^s_i\cap N_j^{a_j}),j>i$ then $\psi^s_i(x)={\phi}_{{\vp}^{u}_j,\psi^s_j}(x) $.  
\end{demo}

Let $n\geq 1$, $i\in\{0,\dots,n-1\}$ and let $G_i$ be the union of all 
stable one-dimensional separatrices of saddle points {in} $\Si_i$ {which contains} heteroclinic points. Let $\check G_i\subset G_i$ be {the} union of separatrices {in} $G_i$ such that $G_i=\bigcup\limits_{\ga\in\check G_i}orb(\ga)$ and, for every pair $(\ga_1, \ga_2)$ of distinct separatrices in $\check G_i$ {and  every $k\in \mathbb Z$,} one has $\ga_{2}\neq f^k(\ga_1)$. {For $\ga\in G_i$ with the end point $p\in \Si_i$ and a point $q\in\Si_j,\,j>i$,  let us consider a sequence of different periodic orbits $p=p_0\prec p_1\prec\dots\prec p_k=q$ such that $\ga\cap W^u_{p_1}\neq\emptyset$, the length of the longest such chain is denoted by $beh(q|\ga)$.}

Let $\ga\in\check G_i$ be a separatrix of $p\in \Si_i$ {and let $N^t_\ga$ be the connected component of $N^t_p \smallsetminus W^u_p$ which contains $\ga$}. We endow {with} the index $\ga$ {(resp. $p$)  the preimages in $M$ (through the linearizing map $\mu_p$) of  all objects from the linear model $\mathcal N$ associated with the separatrix $\ga$ (resp. $p$); for being precise we decide that $\mu_p(\ga)=Ox_3^+$.} For a separatrix $\ga$ in $\check G_i$, let us fix a saddle point $q_\ga$ such that $beh(q_\ga|\ga)=1$. Notice that the intersection $\ga\cap W^u_{q_\ga}$ consists of a finite number {of} heteroclinic orbits. Let $T_p=W^u_p\cap W^s_{\Omega_1}$.

\begin{lemm}\label{order(n-1)} Let $n\geq 1$, {$i\in\{0,\dots,n-1\}$}. For every $\ga\in\check G_{i}$ there are  positive numbers $\rho$, $\de$ such that for every heteroclinic point {$Z^0_\ga\in (\ga\cap W^u_{q_\ga})$} with $z^0<\ep$ the following properties hold:

{\rm (1)} $U_{p}$ avoids all heteroclinic points;

{\rm (2)} either $c_p=U_p\cap T_p$ or the sets $\partial U_p$ and $T_p$ intersected transversally and each path-connected component of the intersection $U_p\cap T_p$ is a segment which intersects each from the both connected components of $\partial U_p$ at a unique point.

For chosen $c_p$ there is a positive number $\ep$  such that for every heteroclinic point {$Z^0_\ga\in (\ga\cap W^u_{q_\ga})$} with $z^0<\ep$ the following properties hold:

{\rm (3)} {$\vp(\tilde d_p)\subset {\phi}_{{\vp}^{u}_i,\psi^s_i}(N^{a_i}_i)$};

{\rm (4)} $\vp(\tilde c_p)\cap{\phi}_{{\vp}^{u}_i,\psi^s_i}(\tilde c^0_p)=\emptyset$, $\vp(\tilde c_p^1)\cap{\phi}_{{\vp}^{u}_i,\psi^s_i}(\tilde c_p^0)=\emptyset$ and $\vp(\tilde\beta_\ga)\subset{\phi}_{{\vp}^{u}_i,\psi^s_i}(\tilde V_\ga)$.
\end{lemm}    
\begin{demo} Let $\ga\in{\check G_{i}},\,i\in\{0,\dots,n-1\}$. If $W^u_p$ contains a compact heteroclinic curve which is non-contractible in $W^u_p\setminus p$ then we assign $c_p$ to be this heteroclinic curve. In the opposite case, due to Lemma \ref{Ge}, there is a generic  $\rho>0$ such that the  curve $c_p$ avoids all heteroclinic points. Since ${W^s_{l}}$ accumulates on $W^s_{k}$ for every $l<k$, then $K_{p}\cap W^s_{i-1}$ is made of a finite number of heteroclinic points $y_1,\dots,y_r$ which we can cover {by} closed 2-discs $b_1,\dots,b_r\subset int\,K_p$. In $K_{p}\smallsetminus int(b_1\cup\dots\cup b_r)$ there is a finite number of heteroclinic points from $W^s_{i-2}$ which we cover by the union of a finite number {of} closed 2-discs, and so on. Thus we get that all heteroclinic points in $K_{p}$ belong to the union of {finitely many} closed 2-discs avoiding $\partial K_{p}$. Therefore, there is $\delta\in(0,\frac{\rho}{4})$ such that  $U_{p}$ avoids heteroclinic points. This proves item (1). As the set $T_p$ is closed one-dimensional $C^{1,1}$-lamination then due to the theory of the general position there is a generic $U_p$ with the property (2).  

By assumption of Theorem \ref{t.invariant}, $\vp$ is defined on the complement of the stable manifolds and, by definition, ${\phi}_{{\vp}^{u}_{i},\psi^s_{i}}$ coincides with $\vp$ on $W^u_i\smallsetminus L^s$, and hence on $U_p$. As $\vp$ and ${\phi}_{{\vp}^{u}_{i},\psi^s_{i}}$ are continuous, we can choose $\ep>0$ sufficiently  small so that, if $Z^0_{\ga}$ is any  heteroclinic point in the intersection $\ga\cap W^u_{q_\ga}$ with $z^0<\ep$, the requirements of (3) and (4) are fulfilled.
\end{demo}

Let us fix $U_p$ satisfying  items (1)-(2) of Lemma \ref{order(n-1)} and let us define
$$U_i=\bigcup\limits_{p\in \Si_i}\left(\bigcup\limits_{k=0}^{per(p)-1}f^k(U_{p})\right),\quad K_i=\bigcup\limits_{p\in \Si_i}\left(\bigcup\limits_{k=0}^{per(p)-1}f^k(K_{p})\right).$$

\begin{lemm}\label{order} Let $n\geq 2$. For every $i\in\{0,\dots,{n-2}\}$ and $\ga\in\check G_i$, there is a heteroclinic point $Z^0_\ga\in \ga$ satisfying the conditions of Lemma \ref{order(n-1)} and in addition:
  $$\mathcal T_\ga\cap\tilde U_j=\emptyset~~~for~~~j\in\{i+1,\dots,{n-1}\}.$$
\end{lemm}
\begin{demo} {In this statement, {it is meant} that  $\tilde U_{n-1}$ is associated with the points $Z^0_{\ga}, \ga\in \check G_{n-1}$ given by Lemma \ref{order(n-1)} and $\tilde U_{i}$ is associated with the points $Z^0_{\ga}, \ga\in\check G_{j}$ given by Lemma \ref{order} for every $j>i$. Therefore, it makes sense to prove this Lemma by decreasing induction on $i$ from ${i=n-2}$ to $0$. That is what is done below. It is also worth noticing  that nothing is required with respect to $\Si_n$; the reason why is that the {stable separatrices of $\Si_n$} have no heteroclinic points.}

Let us first prove the lemma  for $i= n-2$. Let $\ga\in {\check G_{n-2}}$ and let $p$ be the saddle end point of $\ga$. Notice that the intersection $\ga\cap K_{n-1}$ consists of a finite number points $a_1,\dots,a_l$ avoiding $U_{n-1}$. Let $d_{1},\dots,d_{l}\subset K_{n-1}$ be compact discs with centres $a_1,\dots,a_l$ and radius $r_*$ (in linear coordinates of $N_p$) avoiding $U_{n-1}$. Let us choose a number $n^*\in\mathbb N$ such that $\frac{\rho}{2^{n^*}}<r^*$. Let $Z^*_\ga\subset\ga$ be a point such that the segment $[p,Z^*_\ga]$ of $\ga$ avoids $\tilde K_{n-1}$ and $\mu_p(Z^*_\ga)=Z^*=(0,0,z^*)$ where $z^*<\ep$. Then every heteroclinic point $z^0_{\ga}$ so that $z^0<\frac{z^*}{2^{n^*}}$ possesses the property: $\mathcal T_\ga\cap\tilde K_{n-1}$ avoids $\tilde U_{n-1}$. 

{For the induction, let us assume now that the construction of the desired heteroclinic points}
is done for $i+1,i+2,\dots,n-2$. Let us do it for $i$. Let $\ga\in  {\check G_i}$. By assumption of the induction $(\bigcup\limits_{k=i+1}^{j-1}\mathcal T_k)\cap\tilde U_j=\emptyset$ for $j\in\{i+2,\dots,n-1\}$. Since ${W^s_{k-1}}$ accumulates on $W^s_{k}$ for every $k\in\{0,\dots,n\}$, then $(\bigcup\limits_{k=i+1}^{j-1}\mathcal T_k)\cap K_j$ is a compact subset of $K_j$ and the intersection $(\ga\setminus(\bigcup\limits_{k=i+1}^{j-1}\mathcal T_k))\cap K_j$ consists of a finite number points $a_1,\dots,a_l$ avoiding $U_{j}$. Let $d_{1},\dots,d_{l}\subset K_{j}$ be compact discs with centres $a_1,\dots,a_l$ and radius $r_*$ (in linear coordinates of $N_p$) avoiding $U_{j}$ and such that $r_*$ is less than  {the distance between} $\partial(K_j\setminus U_j)$ and $(\bigcup\limits_{k=i+1}^{j-1}\mathcal T_k)\cap K_j$. Similar to the case $i=n-2$ it is possible to choose a heteroclinic point $Z^0_{\ga}$ sufficiently close to the saddle $p$ where $\ga$ ends such that the set $(\mathcal T_\ga\setminus(\bigcup\limits_{k=i+1}^{j-1}\mathcal T_k))\cap \tilde K_j$ avoids $\tilde U_j$. 
\end{demo}

Everywhere below, we assume that for every $\ga\subset\check G_i$ the neighborhood  
$\mathcal T_\ga$  {satisfies} to Lemmas \ref{order(n-1)}   {and} \ref{order}. Let $$\mathcal T_i=\bigcup\limits_{\ga\subset\check G_i}\left(\bigcup\limits_{k=0}^{per(\ga)-1}f^k(\mathcal T_{\ga})\right).$$ 

For $\ga\subset\check G_i,\, j>i$, let us denote by ${\mathcal{J}}_{\ga,j}$ the union of all connected components of ${W}^u_{j}\cap\mathcal T_\ga$ which do not lie in {$int\,\mathcal T_k$} with $i<k<j$. 
Let $\mathcal{J}_\ga=\bigcup\limits_{j=i+1}^{n}{\mathcal{J}}_{\ga,j}$ {and $$\mathcal J_i=\bigcup\limits_{\ga\subset\check G_i}\left(\bigcup\limits_{k=0}^{per(\ga)-1}f^k(\mathcal J_{\ga})\right).$$} 

Let {$\mathcal W_\ga$} be the fundamental domain of $f^{per(\ga)}\vert_{\mathcal T_\ga\setminus W^u_p}$
 {limited by the plaques of the two heteroclinic  points $Z^0_\ga$ and $f^{per(\ga)}Z^{0}_\ga$. Notice that}
 $\ga\cap{\mathcal W_\ga}$ is a fundamental domain of $f^{per(\ga)}\vert_\ga$. 
 
\begin{lemm}\label{nado} The set ${\mathcal{J}}_{\ga,j}\cap{\mathcal W_\ga}$ consists of a finite number  {of}  closed 2-discs.
\end{lemm} 
\begin{demo} Let $\hat\ga=p_i(\ga)$, $\hat{\mathcal T}_\ga=p_i(\mathcal T_\ga)$, $\hat{\mathcal T}_{i,j}=p_i(\mathcal T_j)$. Since ${W^u_{j}}$ accumulates on $W^u_{l}$ {only when} $l<j$, then the set
$\hat{W}^u_{j}\setminus\bigcup\limits_{l=i+1}^{j-1}\hat{\mathcal T}_{l,i}$ is a compact set. Due to Lemma \ref{order} the intersection of $\hat{\mathcal T}_\ga\cap\partial{\hat{\mathcal T}}_{l,i}$ consists of 2-discs which are projection with respect to $p_i$ of leaves of the foliation $F^u_i$. Thus the intersection $\hat\ga\cap(\hat{W}^u_{j}\setminus\bigcup\limits_{l=i+1}^{j-1}\hat{\mathcal T}_{l,i})$ consists of a finite number of closed 2-discs.
\end{demo}

Due to Lemma \ref{nado}, the set $\mathcal{J}_\ga\cap\ga\cap{\mathcal W_\ga}$ consists of a finite number {of heteroclinic points;} denote them $Z^2_\ga, \ldots, Z^{m}_\ga$ ($m$ depends on $\ga$).  {Finally, choose an} arbitrary point $Z^1_\ga\in\ga$ so that the arc $(z^0_\ga,z^1_\ga)\subset\ga$ does not contain heteroclinic points from $\mathcal J_\ga$. Let us construct $\mathcal H_\ga$ using the point $Z^1=\mu_p(Z^1_\ga)$. Without loss of generality we will assume that $\mu_p(Z^i_\ga)=Z^i=(0,0,z^i)$ for $z^0>z^1>\dots>z^m>\frac{z^0}{4}$ and $\mu_p(T_p)=T$. For $i=0, \ldots, n-1$ let $$\mathcal H_i=\bigcup\limits_{\ga\subset\check G_i}\left(\bigcup\limits_{k=0}^{per(\ga)-1}f^k(\mathcal H_{\ga})\right),\quad \mathcal {M}_{i}=V_f\bigcup\limits_{k=0}^{i}{(G_k\cup \Si_k)}{\quad{\rm and}\quad}\mathcal {M}'_{i}=V_{f'}\bigcup\limits_{k=0}^{i}{(G'_k\cup \Si'_k)}.$$

{\begin{lemm} \label{main0} There is an equivariant homeomorphism ${\varphi}_0:\mathcal{M}_0\to \mathcal{M}'_0$ with following properties:

{\rm 1)} ${\varphi}_0$ coincides with ${\varphi}$ out of $\mathcal{T}_0$;

{\rm 2)} ${\varphi}_0\vert_{\mathcal H_0}={\phi}_{\psi^u _{0},\psi^s_{0}}\vert_{\mathcal H_0}$,
 where $\psi^u_0={\varphi}\vert_{W^u_0}$;

{\rm 3)} $\vp_0(W^u_1)=W'^u_1$ and $\vp_0(W^u_k\setminus\bigcup\limits_{j=1}^{k-1}int\,\mathcal T_j)\subset W'^u_k$ for every $k\in\{2,\dots,n\}$;

{\rm 4)} $\vp_0(W^s_{\Om_1}\cap\mathcal M_0)=W^s_{\Om'_1}\cap\mathcal M'_0$.
\end{lemm}}
\begin{demo} The desired $\vp_0$ should be an interpolation between $\vp: V_f\smallsetminus \mathcal T_0 \to  {M'}$ and $\phi_{{\vp}^u _{0},\psi^s_{0}}\vert_{\mathcal H_0}$. Due to Lemma \ref{order(n-1)}  (3) {and the equivariance of the considered maps}, the embedding $${\xi_0}={\phi}^{-1}_{\psi^{u}_0,\psi^s_0} {\vp}:\mathcal T_0\setminus W^s_0\to M$$ is well-defined. Let $\ga\subset\check G_0$ be a separatrix ending at $p\in\Si_0$ and $\xi_\ga=\xi_0|_{\mathcal T_\ga}$. By construction, the topological embedding $\xi=\mu_p\xi_\ga\mu^{-1}_p:\mathcal T\to N$ satisfies to all conditions of Proposition \ref{F1}. Let $\zeta$ be the  embedding from the inclusion of that Lemma and  $\zeta_\ga=\mu_p^{-1}\zeta\mu_p$. Independently, one does the same for every separatrix $\ga\subset\check G_0$. Then, it is extend to all separatrices in $G_0$ by equivariance. {As a result}, we get a homeomorphism $\zeta_0$  {of $\mathcal T_0$ onto $\xi_0(\mathcal T_0)$ which coincides with $\xi_0$ on $\partial\mathcal T_0$. Now, define the homeomorphism $\vp_0:\mathcal M_0\to\mathcal M'_0$ to be equal to ${\phi}_{\psi^{u}_0,\psi^s_0}\zeta_0$ on $\mathcal T_0$ and to $\vp$ on $\mathcal M_0\setminus \mathcal T_0$. One checks the next properties:}

1) ${\varphi}_0$ coincides with ${\varphi}$ out of $\mathcal{T}_0$;

2) ${\varphi}_0\vert_{\mathcal H_0}={\phi}_{\psi^u _{0},\psi^s_{0}}\vert_{\mathcal H_0}$;

3) $\vp_0( {\mathcal J}_0)\subset L^u$;

4) $\vp_0(W^s_{\Om_1}\cap\mathcal M_0)=W^s_{\Om'_1}\cap\mathcal M'_0$.

Property 3) and the definition of the set $ {\mathcal J}_\ga$ imply that $\vp_0(W^u_1)=W'^u_1$ and $\vp_0(W^u_k\setminus\bigcup\limits_{j=1}^{k-1}int\,\mathcal T_j)\subset W'^u_k$ for every $k\in\{2,\dots,n\}$.
Thus $\vp_0$ satisfies  {all  required conditions of the} lemma. 
\end{demo}

\begin{lemm} \label{main} Assume $n\geq 2$, $i\in\{0,\dots,n-2\}${, and assume} there is an equivariant topological embedding ${\varphi}_i:\mathcal{M}_i\to{M'}$ with following properties:

{\rm 1)} ${\varphi}_i$ coincides with ${\varphi}_{i-1}$ out of $\mathcal{T}_i$;

{\rm 2)} ${\varphi}_i\vert_{\mathcal H_i}={\phi}_{\psi^u _{i},\psi^s_{i}}$,
 where $\psi^u_i={\varphi}_{i-1}\vert_{W^u_i}$ {and} $\vp_{-1}=\vp$; 

{\rm 3)} there is an $f$-invariant {union of tubes}
 $\mathcal B_i\subset(\mathcal T_i\cap\bigcup\limits_{j=0}^{i-1}\mathcal H_j)$ {containing
 $(\mathcal T_{i}\cap(\bigcup\limits_{j=0}^{i-1}W^s_j))$} where ${\varphi}_i$ coincides with $\vp_{i-1}$ 
 $($we assume $\mathcal B_0=\emptyset)$;
 
{\rm 4)} $\vp_i(W^u_{i+1})=W'^u_{i+1}$ and $\vp_i(W^u_k\setminus\bigcup\limits_{j=i+1}^{k-1}int\,\mathcal T_j)\subset W'^u_k$ for every $k\in\{i+2,\dots,n\}$;

{\rm 5)} $\vp_i(W^s_{\Om_1}\cap\mathcal M_i)=W^s_{\Om'_1}\cap\mathcal M'_i$.

Then there is a homeomorphism $\vp_{i+1}$ with the same properties {\rm 1)-5)}. 
\end{lemm}
\nd \begin{demo} The desired $\vp_{i+1}$ should be an interpolation between 
 $\vp_i: \mathcal M_{i+1}\smallsetminus \mathcal T_{i+1} \to M'$ and 
 $\phi_{{\psi}^u_{i+1},\psi^s_{i+1}}\vert_{\mathcal H_{i+1}}$ {where $\psi_{i+1}^u=\vp_i\vert_{W^u_{i+1}}$.}
Let $\ga\subset\check G_{i+1}$ be a separatrix ending at $p\in\Si_{i+1}$. {It follows from the definition of the set $\mathcal J_i$ and the choice of the point $q_\ga$ that $(W^u_{q_\ga}\cap \mathcal T_i)\subset \mathcal J_i$. Then, due to condition 4) for $\vp_i$ we have $\vp_i(W^u_{q_\ga}\cap \mathcal T_i)\subset W^u_{q'}$. By the property 1) of the homeomorphism $\vp_i$ and the properties of $\mathcal T_{i+1}$ from Lemmas \ref{order(n-1)} {(1)} and \ref{order}, we get that $\vp_i|_{\tilde U_p}=\vp|_{\tilde U_p}$. Then ${\phi}_{{\vp}^{u}_{i+1},\psi^s_{i+1}}|_{\tilde U_p}={\phi}_{{\psi}^{u}_{i+1},\psi^s_{i+1}}|_{\tilde U_p}$. Thus it follows from the property (3) in Lemma \ref{order(n-1)} that the following embedding is well-defined: $\xi_\ga={\phi}_{{\psi}^{u}_{i+1},\psi^s_{i+1}}^{-1}\vp_i:{\mathcal T}_\ga\setminus(\ga\cup p)\to M'$.}

By construction, the topological embedding $\xi=\mu_p\xi_\ga\mu^{-1}_p$ satisfies to all conditions of Proposition \ref{F1}. Let $\zeta$ be the  embedding {which is yielded by that proposition}. Define $\zeta_\ga=\mu_p^{-1}\zeta\mu_p$. Notice that by the property 3) of the homeomorphism $\psi^s$ in Lemma \ref{1dim} and by the properties $\psi^u_{i+1}=\vp_i|_{W^u_i}$, we have that $\zeta_{\ga}$ is {the} identity on a neighborhood $\tilde B_{\ga}\subset(\mathcal T_\ga\cap\bigcup\limits_{j=0}^i\mathcal H_j)$ of $\mathcal T_{\ga}\cap(\bigcup\limits_{j=0}^iW^s_j)$. Independently, one does the same for every separatrix $\ga\subset\check G_{i+1}$. {Assuming that $\zeta_{f(\ga)}=f'\zeta_{\ga}f^{-1}$  and $\tilde B_{i+1}=\bigcup\limits_{\ga\subset\check G_{i+1}}\left(\bigcup\limits_{k=0}^{per(\ga)-1}f^k(\tilde B_{\ga})\right)$ we get a homeomorphism $\zeta_{i+1}$ on $\mathcal T_{i+1}$.}  Thus the required homeomorphism coincides with ${\phi}_{\psi^{u}_{i+1},\psi^s_{i+1}}$ on $\mathcal H_{i+1}$ and with $\vp_i$ out of $\mathcal T_{i+1}$. 
\end{demo}

{Let $G$ be the union of all stable one-dimensional separatrices which do not contain heteroclinic points, $ N^t_{G}=\bigcup\limits_{\gamma\subset G} N^t_\gamma$ and  $$\mathcal {M}=\mathcal M_{n-1}\cup G.$$ Also we have similar objects with prime for $f'$.

{\begin{lemm} \label{main-} There are {numbers} 
$0<\tau_1<\tau_2<1$ and an equivariant embedding $h_{\mathcal M}:\mathcal M\to {M'}$ with the following properties: 

{\rm 1)} $h_{\mathcal M}$ coincides with {${\varphi}_{n-1}$} 
out of $N^{\tau_2}_{G}$;

{\rm 2)} $h_{\mathcal M}$ coincides with ${\phi}_{{\vp}_{n-1},{\psi}^s}$ on $\vert_{\mathcal N^{\tau_1}_{G}}$,
{where $\psi^s: L^s\to L'^s$ is yielded by 
Lemma \ref{1dim};}

{\rm 3)} there is an $f$-invariant neighborhood of the set {$N_{G}\cap(G_0\cup\dots\cup G_{n-1})$} where $h_{\mathcal M}$ coincides with {$\vp_{n-1}$};

{\rm 4)} $h_{\mathcal M}(W^s_{\Om_1}\cap\mathcal M)=W^s_{\Om'_1}\cap\mathcal M'$.
\end{lemm}}
\begin{demo} Let $\check G\subset G$ be a union of separatrices from $G$ such that $\ga_{2}\neq f^k(\ga_1)$ for every $\ga_1,\ga_2\subset\check G$, $k\in\mathbb Z\setminus\{0\}$ and $G=\bigcup\limits_{\ga\in\check G}orb(\ga)$. Let $i\in\{0,\dots,n\}$, $p\in\Si_i$ and $\gamma\subset G$. 

Notice that $\left(N_{\ga}\setminus(\ga\cup p)\right)/f^{per(\ga)}$ is homeomorphic to $\hat X\times[0,1]$ where $\hat X$ is 2-torus and the natural projection $\pi_{\ga }: N_{\ga }\setminus(\ga \cup p)\to \hat X\times[0,1]$ sends $\partial N^t_{\ga }$ to $\hat X\times\{t\}$ for each $t\in(0,1)$ and sends $W^u_p \setminus p$ to $\hat X\times\{0\}$. Let  $\xi_{\ga }={\phi}^{-1}_{{\vp}_{n-1}|_{W^u_i},{\psi}^s_i}\vp_{n-1}|_{N_{\ga }^{a_i}\setminus(\ga \cup p)}$ and  $\hat\xi_{\ga}=\pi_{\ga }\xi_{\ga }\pi^{-1}_{\ga }|_{\hat X\times[0,a_i]}$. Due to item 3) of Lemma \ref{main}, the homeomorphism $\hat\xi_{\ga }$ coincides with {the} identity in some neighborhood of $\pi_\ga(N^{a_i}_{\ga}\cap(G_0\cup\dots\cup G_{n-1}) )$. Let us choose this neighborhood of the form $B_\ga\times[0,a_i]$. Let $\hat T_\ga=\pi_\ga(T_p)$.  Let us choose numbers $0<\tau_{1,\ga }<\tau_{2,\gamma }<a_i$ such that $\hat \xi_{\ga }(\hat X\times[0,{\tau_{2,{\ga }}}])\subset\hat X\times[0,\tau_{1,\ga })$. {By the construction,} $\hat \xi_{\ga }:\hat X\times[0,{\tau_{2,{\ga }}}]\to {X\times[0,1]}$ {is} a topological embedding which is the identity on $\hat X\times\{0\}$,  
$\hat \xi_{\ga }|_{B_{\ga }\times[0,\tau_{2,\ga}]}=id|_{B_{\ga }\times[0,\tau_{2,\ga }]}$ and, due to item 4) of Lemma \ref{main}, $\hat\xi_{\ga }(\hat T_\ga\times[0,\tau_{2,\ga}])\subset\hat T_\ga\times[0,1]$. 
Then, due to Proposition \ref{F2}, 

1. there is a homeomorphism $\hat\zeta_{\ga }:X\times[0,\tau_{2,\ga }]\to \hat\xi(X\times[0,\tau_{2,\ga }])$ such that  $\hat\zeta_{\ga }$ is identity on $X\times[0,{\tau_{1,{\ga }}}]$ and is $\hat \xi_{\ga }$ on $X\times\{\tau_{2,\ga}\}$. 

2. $\hat \zeta_{\ga }|_{B_{\ga }\times[0,\tau_{2,\ga }]}=id|_{B_{\ga }\times[0,\tau_{2,\ga }]}$.

3. $\hat\zeta(\hat{T}_\ga\times[0,\tau_{2,\ga}])\subset{\hat{T}_\ga\times[0,1]}$. 

Let $\zeta_{\ga }$ be a lift of $\hat \zeta_{\ga }$ on $N^{\tau_{2,\ga}}_{\ga }$ which $\xi_\ga$ on $\partial N^{\tau_{2,\ga}}_{\ga }$. Thus $\vp_{\ga }={\phi}_{{\vp}_{n-1}|_{W^u_i},{\psi}^s_i}\zeta_{\ga }$ is the desired extension of $\vp_{n-1}$ to ${ N}_{\ga }$. Doing the same for every separatrix $\ga \subset \check G$ and {extending it to the other separatrices from $G$ by equivariance,}
we get the {required} homeomorphism $h_{\mathcal M}$ for $\tau_1=\min\limits_{\gamma \subset\check G}\{\tau_{1,\ga }\}$ and $\tau_2=\min\limits_{\gamma \subset\check G}\{\tau_{2,\ga }\}$.
\end{demo}

\section{Topological background}\label{IIII}

The following proposition is a corollary of Theorem 3.1 from \cite{GrMeZh}. In fact in paper \cite{GrMeZh} the objects are required to be smooth, but actually the results are true in the case when the objects are tame.

\begin{prop}\label{solid} Let $P$ be homeomorphic to $K\times[0,1]$, where $K=\mathbb S^1\times[0,1]$ and $Q\subset P$ is a tame embedded annulus such that $P\setminus Q$ is not connected and the annuli $K\times\{0\}$, $K\times\{1\}$ belong to the different connected components of $P\setminus Q$. Then the set $P\setminus Q$ consists of two connected components, the closure of each of which is homeomorphic to $P$.
\end{prop}

\begin{prop}\label{B2} Let 
\begin{itemize}
\item $C$ be a compact subset of $[0,1]$ including $0$ and $1$; 
\item $\mathcal L$ be a lamination $\{L_t=\mathbb R^2\times\{t\}\}_{t\in C}$;
\item there is a tame topological embedding $g:\partial(\mathbb D^2\times[0,1])\to\mathbb R^2\times[0,1]$ such that $g(\partial(\mathbb D^2\times\{t\}))\subset L_t$ and $g^{-1}(\partial(\mathbb D^2\times\{t\}))\subset L_t$ for any $t\in C$.
\end{itemize}
Then there is a homeomorphism $h:\mathbb R^2\times[0,1]\to\mathbb R^2\times[0,1]$ such that $h(L_t)=L_t$ for any $t\in C$ and $h=g$ on $\partial(\mathbb D^2\times[0,1])$.
\end{prop}
\begin{demo} Let us introduce the canonical projection $p:\mathbb R^2\times[0,1]\to\mathbb R^2$, where $p(r,t)=r$. Let us consider a homotopy $g_t:\partial\mathbb D^2\to\mathbb R^2,t\in[0,1]$ given by formula $g_t(x)=p(g(x,t))$. By \cite{Ep}, there is an isotopy $\bar g_t:\partial\mathbb D^2\to\mathbb R^2,t\in[0,1]$ such that $\bar g_t=g_t$ for $t\in C$ and $\lim\limits_{t\to C}||g_t(x)-\bar g_t(x)||=0$. Let us extend this isotopy up to an isotopy $\bar G_t:\mathbb R^2\to\mathbb R^2$. Let $\bar G:\mathbb R^2\times[0,1]\to\mathbb R^2\times[0,1]$ be a homeomorphism given by the formula $\bar G(r,t)=(\bar G_t(r),t)$. 

Let $Q=\partial(\mathbb D^2\times[0,1])$ and $\bar Q=\bar G^{-1}(g(\partial(\mathbb D^2\times[0,1])))$. Let us define $\psi:\partial Q\to\partial\bar Q$ by the formula $\psi=\bar G^{-1}g$. By the construction $\psi$ identity on $\mathcal L$ and $\lim\limits_{t\to C}||\psi(x,t)-(x,t)||=0$. Let us show that there is a map $\Psi:\mathbb R^2\times[0,1]\to\mathbb R^2\times[0,1]$ which coincides with $\psi$ on $\partial Q$, is the identity on $\mathcal L$ and such that $\lim\limits_{t\to C}||\Psi(x,t)-(x,t)||=0$, thus $h=\bar G\Psi$ will be  the required homeomorphism. 

For this aim let us denote by $A_{a,b}$ a connected component of $Q\setminus\mathcal L$ bounded by the leaves $L_a,L_b,a,b\in C$. If the set $C$ is finite then, by Alexander trick, there is a homeomorphism $\Psi_{a,b}:\mathbb R^2\times[a,b]\to\mathbb R^2\times[a,b]$ which is $\psi$ on $A_{a,b}$ and is identity on $L_a\cup L_b$. Then $\Psi$ is composed by $\Psi_{a,b}$. If $C$ is infinite then for a sequence of annuli $A_{a_n,b_n}$ such that  $\lim\limits_{n\to\infty}(b_n-a_n)=0$, below we construct homeomorphisms $\Psi_{n}:\mathbb R^2\times[a_n,b_n]\to\mathbb R^2\times[a_n,b_n]$ which is $\psi$ on $A_{a_n,b_n}$, is identity on $L_{a_n}\cup L_{b_n}$ and such that $\lim\limits_{n\to\infty}||\Psi_{n}(x,t)-(x,t)||=0$, what finishes the proof. 

Let $\psi_n=\psi|_{A_{a_n,b_n}}$. As $\lim\limits_{t\to C}||\psi(x,t)-(x,t)||=0$ then there is a sequence $\delta_n$ which tends to $0$ as $n\to\infty$ and such that $b_n-a_n<\delta_n$ and $\psi_n$ moves no point more than $\delta_n$. Let $U_n$ be a solid torus which is the one-sided $\delta_n$-neighborhood of $A_{a_n,b_n}$ and $\Si_n=\partial U_n$. Let $\phi_n:\Si_n\to\mathbb R^2\times[a_n,b_n]$ be a topological embedding which is $\psi_n$ on $A_{a_n,b_n}$ and identity on the other part of $\Si_n$. 

Let $\Si'_n=\phi_n(\Si_n)$ and $U'_n$ is a solid torus bounded by $\Si'_n$ (see Proposition \ref{solid}). Let us choose in $U_n$ an even number of vertical meridian discs $D^1_n,\dots,D^{2k_n}_n$ with distance between of them less than $3\de_n$ and such that $\phi_n(\partial D^{2i-1}_n)$ avoids $\bigcup\limits_{i=1}^{k_n}D^{2i}_n$. The closed curve $\phi_n(\partial D^{2i-1}_n)$ is a meridian in the torus $\Si'_n$, hence, it is the boundary of a disc $D'^{2i-1}_n$ in $U'_n$ whose  interior avoids $\Si'_n$. By the standard pushing procedure we can get that $D'^{2i-1}_n$ avoids $\bigcup\limits_{i=1}^{k_n}D^{2i}$. Since  every connected component of the sets  $U_n\setminus\bigcup\limits_{i=1}^{k_n}D^{2i-1}_n$ and $U'_n\setminus\bigcup\limits_{i=1}^{k_n}D'^{2i-1}_n$ is a 3-ball then there is a homeomorphism $\Phi_n:U_n\to U'_n$ such that is $\phi_n$ on $\Si_n$ and sends $D^{2i-1}_n$ to $D'^{2i-1}_n$. By the construction it moves no point more than $4\delta_n$.  

Doing the same for the other one-sided $\delta_n$-neighborhood of $A_{a_n,b_n}$ and extending by identity out of the $\delta_n$-neighborhood of $A_{a_n,b_n}$, we get $\Psi_n$.  
\end{demo}

\begin{rema}\label{B1} A similar proposition obviously true for a similar one-dimensional lamination in $\mathbb R^1\times[0,1]$.  
\end{rema}

\begin{prop} \label{F2} Let $\hat X$ be a compact topological space, $0<{\tau_1}<{\tau_2}<1$ and $\hat\xi:\hat X\times[0,{\tau_2}]\to\hat X\times[0,1]$ be a topological embedding which is the identity on $\hat X\times\{0\}$ and $\hat X\times[0,\tau_1]\subset\hat \xi(\hat X\times[0,{\tau_2}])$. Then

1. there is a homeomorphism $\hat\zeta:\hat X\times[0,\tau_2]\to \xi(\hat X\times[0,\tau_2])$ such that  $\hat \zeta$ is identity on $\hat X\times[0,{\tau_1}]$ and is $\hat\xi$ on $\hat X\times\{\tau_2\}$. 

2. if for a set $\hat{B}\subset\hat X$ the equality $\hat\xi|_{\hat{ B}\times[0,\tau_2]}=id|_{\hat{ B}\times[0,\tau_2]}$ is true then $\hat \zeta|_{\hat{ B}\times[0,\tau_2]}=id|_{\hat{ B}\times[0,\tau_2]}$.

3. if for a set $\hat T\subset\hat X$ the inclusion $\hat\xi(\hat{T}\times[0,\tau_2])\subset{\hat{T}\times[0,1]}$ is true then $\hat\zeta(\hat{T}\times[0,\tau_2])\subset{\hat{T}\times[0,1]}$.
\end{prop} 
\begin{demo} Let us choose $l\in(\tau_1,\tau_2)$ such that $\hat X\times[0,l]\subset\hat \xi(\hat X\times[0,{\tau_2}])$. Define a homeomorphism $\kappa:[\tau_1,1]\to[0,1]$ by the formula $$\kappa(t)=\cases{(x,\frac{l(t-\tau_1)}{l-\tau_1}),~t\in[\tau_1,l];\cr
(x,t),~t\in[l,1].\cr}$$ Let $\mathcal K(x,t)=(x,\kappa(t))$ on $\hat X\times[\tau_1,1]$. Then the required homeomorphism can be defined by the formula $$\hat \zeta(x,t)=\cases{(x,t),~t\in[0,\tau_1];\cr
\mathcal K^{-1}\xi(\mathcal K((x,s)))),~s\in[\tau_1,\tau_2].\cr}$$ Properties 2 and 3 automatically follows from this formula.
\end{demo}

\end{document}